\documentclass[a4paper]{article}
\usepackage[all]{xy}
\usepackage{latexsym,amssymb,amsfonts,amsmath}
 \newtheorem{definition}{Def.}[section]
 \newtheorem{theorem}[definition]{Theorem}
 \newtheorem{lemma}[definition]{Lemma}
 \newtheorem{cor}[definition]{Corollary}
 \newtheorem{prop}[definition]{Proposition}
 \newtheorem{remark}[definition]{Remark}
\newcommand{\orig}[1]{\mathsf{#1}}

\newcommand{\nat}{\mbox{\cal N}}

\newcommand{\mtt}{\mbox{\rm mTT}}
\newcommand{\mttdp}{\mbox{\rm mTT$^{dp}$}}
\newcommand{\cmtt}{\mbox{\rm cmTT}}
\newcommand{\mttps}{\mbox{\rm MLTT}}
\newcommand{\qmtt}{\mbox{\rm qmTT}}
\newcommand{\emtt}{\mbox{\rm emTT}}
\newcommand{\emttdp}{\mbox{\rm emTT$^{dp}$}}

\newcommand{\emtts}{\mbox{{\rm emTT}$_{set}$}}
\newcommand{\mtts}{\mbox{\rm mTT$_{set}$}}
\newcommand{\cemtt}{\mbox{\rm cemTT}}
\newcommand{\cq}{\mbox{\rm Q(mTT)}}
\newcommand{\cqset}{\mbox{\rm Q(mTT)$_{set}$}}
\newcommand{\cqdp}{\mbox{\rm Q(mTT$^{dp}$)}}
\newcommand{\cqis}{\mbox{\rm Q(mTT)/$\simeq$}}
\newcommand{\cqdpis}{\mbox{\rm Q(mTT$^{dp}$)/$\simeq$}}

\newcommand{\cqac}{\mbox{\rm Q(MLTT)}}
\newcommand{\cset}{\mbox{\rm {\cal C}(MLTT)}}

\CompileMatrices

\begin{document}

\title{
A minimalist two-level foundation for constructive mathematics
}
\author{Maria Emilia Maietti \\
 Dipartimento di Matematica Pura ed Applicata\\
 University of Padova\\
  via Trieste n. 63 - 35121  Padova,
 Italy \\
maietti@math.unipd.it
}
\date{February 6, 2009}
\maketitle

\begin{abstract}
We present a two-level theory to formalize
constructive mathematics as advocated in a previous paper
with G. Sambin.

 One level is given
by an intensional type theory, 
called Minimal type theory. This theory
 extends a previous version  with collections.

 The other level is given
 by an extensional set theory that  is interpreted
in the first one by means of a quotient model.

This two-level theory has two main features: it is minimal
among the most relevant foundations for constructive mathematics;
it is constructive thanks to the way the extensional level
is linked to the intensional one which 
 fulfills the
``proofs-as-programs'' paradigm and
 acts as a programming language.
\end{abstract}




\section{Introduction}
In a previous paper  with G. Sambin  \cite{mtt} we argued about the necessity of 
building a foundation for constructive mathematics to be taken as a common core
among relevant existing foundations in axiomatic set theory, such as Aczel-Myhill's CZF theory~\cite{czf},
or in category theory, such as the internal theory of a topos (for example
in \cite{tumscs}), 
or in type theory, such  as Martin-L{\"o}f's type theory~\cite{PMTT} and
 Coquand's Calculus of Inductive Constructions~\cite{TC90}.

There we also argued what it means for a foundation   to be constructive.
The idea is that a foundation to develop  mathematics is constructive if it satisfies
 the ``proofs-as-programs'' paradigm, namely if it enjoys a realizability model
 where to extract
programs from proofs. If such a semantics is defined in terms of
 Kleene realizability~\cite{DT88},
 then the foundation
turns out to be consistent with the formal Church thesis (for short CT) 
and the axiom of choice (for short AC).
In \cite{mtt} we took such a consistency property
as our formal notion of  ``proofs-as-programs'' foundation.
This notion appears very technical in comparison
with the intuitive  proofs-as-programs paradigm.
Actually, one of our aims is  to explore which commonly conceived
proofs-as-programs theories satisfy our formal notion of consistency
with CT and AC together.

To this purpose in   \cite{mtt}, we
first noticed that theories
satisfying  extensional properties, like extensionality of functions, can not
satisfy
our proofs-as-programs requirement.
This is due to the well known result by Troelstra~\cite{tr} 
that in intuitionistic arithmetics on finite types
extensionality of functions  is inconsistent 
with CT+AC  (see \cite{constrII} or 
\cite{Beeson}).  At present, only theories presented
in terms of an intensional type theory, such as  Martin-L{\"o}f's one in \cite{PMTT}, seem to fit into our paradigm.

This led us to conclude that 
in a proofs-as-programs theory one can only represent extensional concepts
by modelling them via intensional ones in a suitable way, for example
as done in \cite{toolbox,altlics,altob,disttheshof}.

Therefore in \cite{mtt} we ended up  
in saying that a  {\it constructive}  foundation for mathematics must be  equipped with
{\it two levels} of the following kind: with an {\it intensional} level that acts as a programming language
and is the actual proofs-as-programs theory; with an {\it extensional}
level that acts as the set theory where to formalize mathematical proofs.
Then, the constructivity of the whole foundation relies on the fact that
 the extensional level must be implemented over the intensional level, but not only this.
Indeed, in \cite{mtt} following Sambin's forget-restore principle in \cite{toolbox}
 we required that  extensional concepts must be abstractions
of intensional ones as result of forgetting
 irrelevant computational information. Such information
is then  restored when we implement  extensional concepts
back at the intensional level. 
Here we push forward this by saying that the needed extensional
concepts can be obtained by just abstracting over equalities of intensional
ones. Hence,   it is sufficient {\it to build a quotient model
over the intensional theory} in order to be able {\it to implement 
the extensional level in the intensional one}.
In this paper we present an example of such a two-level constructive foundation
that is also minimal among the most relevant constructive ones.

{\bf Intensional level.}
The intensional level of our two-level foundation is essentially
 obtained by extending the typed calculus  introduced in \cite{mtt},
and  called  Minimal Type Theory (for short \mtt),
with  the notion of ``collection'' and corresponding constructors
needed to implement the  power collection of a set as a quotient.
 We still call \mtt\ such an extension.

 Then, the original version of  Minimal Type Theory in \cite{mtt}
essentially represents
the set-theoretic part  of the \mtt\ version presented here, called \mtts.
We say ``essentially''
because here we adopt a version of type theory {\it  with explicit substitution rules and without
the $\xi$-rule for lambda-terms as in \cite{modmar}}.

The reason to adopt such a modified version
is due to the fact that, as suggested to us by P. Martin-L{\"o}f and
T. Streicher, its set-theoretic fragment
\mtts\ 
directly enjoys  Kleene's realizability interpretation of
intuitionistic connectives~\cite{DT88}, and
hence it turns out to satisfy our proofs-as-programs requirement.
We expect that  this is also the case for the whole 
\mtt\  by further extending the realizability
interpretation to support collections. 

Instead, we still do not know whether  the \mtt\ version in \cite{mtt} -
that includes   the $\xi$-rule - 
satisfies our proofs-as-programs requirement,
namely whether it is consistent with the formal Church thesis and the axiom of choice. This problem
can be reduced to asking whether intensional Martin-L{\"o}f's type theory
 in \cite{PMTT}, even with no universes,
is consistent with the formal Church thesis

{\small $$(CT)\qquad\, \forall f \in Nat\rightarrow Nat \quad
\exists e\in Nat\quad\quad
\ (\, \forall x\in Nat\ \ \exists y\in Nat\ \ T(e,x,y) \ \& \ 
U(y)=_{\nat}f(x)\, )$$}

\noindent
where {\small $T(e,x,y)$} is the Kleene predicate expressing that $y$ is
  the computation executed by the program numbered $e$ on the input $x$
and $U(y)$ is output of the computation $y$.
This technical problem seems to be
still open. Luckily, we realized that  we
 do not need to solve such a problem: indeed, we can take 
 a version of \mtt\ with explicit substitutions and without 
 the  $\xi$-rule, to satisfy
our proofs-as-programs requirement in an easier way, because
{\it the absence of the $\xi$-rule does not affect
the properties of the quotient model we will build over it.}
In particular, the quotient model will validate extensionality of 
functions as lambda-terms
anyway. This was first noticed in categorical terms in \cite{Ros,bcsr}.

{\bf Extensional level.}
The extensional level of our two-level foundation is taken to be
 an {\it extensional} dependent type theory with quotients, called \emtt.
This extends that presented in \cite{mai07} with collections
and related constructors needed to represent
the power collection of
a set
with  $\varepsilon$-relation and
comprehension used in everyday mathematical practice.
The set-theoretic part of  \emtt\  includes the fragment without
universes of extensional  Martin-L{\"o}f's type theory in \cite{ML84}.

{\bf Quotient model.} We will  interpret 
 our extensional theory \emtt\ 
in a quotient model 
 built over \mtt.
This model is based on the well-known notion of
total setoid \`a la Bishop~\cite{Bishop}
and the interpretation 
shows that the design of \emtt\ over \mtt\ satisfies
Sambin's forget-restore principle in \cite{toolbox}.
Indeed,
the interpretation 
 represents
the process of restoring all
the irrelevant computational information missing at the extensional level.
Moreover, it turns judgements of \emtt, which
are undecidable as 
 those of extensional  Martin-L{\"o}f's type theory in \cite{ML84}, 
into judgements of \mtt\ that are all decidable.
This forget-restore process 
 is  very evident when looking at the design of \emtt-propositions
and  their interpretation into \mtt.
In fact, whilst in \emtt, as in \mtt,  all propositions 
are identified with collections of their
proofs, in \emtt, despite of \mtt, they  are inhabited by
 at most only one proof in order to express the fact that 
{\it \emtt-propositional
 proofs  are indeed  irrelevant}.
This allows us to introduce  a canonical proof-term $\mathsf{true}$
to express that $\phi$ is valid in \emtt\  
if and only if we can derive  
 $\mathsf{true}\in \phi$ in \emtt.
It is only when we interpret in \mtt\  a derived judgement
of the form $\mathsf{true}\in \phi$  
that we need {\it  to  restore a specific proof-term} containing
all the forgotten computational information about its derivation.  

{\bf Benefits of adopting a two-level theory.} We hope that making explicit an extensional level over an intensional
type theory, as we do here with \emtt\ over \mtt,
will be useful to formalize
  mathematics in intensional type theory.
Indeed, in the current practice of formalizing mathematics in intensional type theory,
 one ends up to work with setoid
constructions, and hence to work {\it within a quotient model}. Here we extract 
a {\it theory valid in
one of  such setoid models}. Therefore, one is then dispensed
to work directly in the model with all the heavy type-theoretic details regarding setoids and families of setoids.
He can develop and formalize his theorems in a simpler extensional theory like our \emtt.
 The interpretation of the extensional
level into the intensional one, given once and for all, 
 guarantees that a formalization of theorems at the extensional level
 is then inherited at the intensional one.

{\bf Open issues.}
Our extensional level \emtt\ does {\it not include} all the
type-theoretic constructors that our quotient model can support.
For example, our quotient model over \mtt\ supports effective
quotients on generic collections and not only effective quotients
over sets as in \emtt.
Moreover,  
in our quotient model every object, which is a quotient of
an intensional set over an arbitrary equivalence relation,
 is covered by a quotient copy of an intensional set, namely by a quotient
of an intensional set over the identity relation.
In the case we build our two-level foundation
by taking   Martin-L{\"o}f's type theory as our
intensional level, the above observation has a very important consequence:
the axiom of choice, which  is not valid over generic quotients, turns out to be valid
 over copies of intensional sets in the quotient model.
 This implies that  we can
consistently add  an axiom expressing that
 ``every extensional set is covered by a set satisfying the axiom of choice''
 to  
the extensional theory abstracted  over Martin-L{\"o}f's type theory.
This axiom  was first noticed by P.Aczel and expressed as the Presentation axiom in his CZF theory (see \cite{czf}).
We leave to future research how to formalize such an axiom in an extensional type
theory abstracted over  Martin-L{\"o}f's one, and more generally how to formalize  the precise {\it internal
language of our quotient model} over \mtt\ or extensions,
namely the theory that fully captures all the type-theoretic constructors
that can be modelled via quotients in it.

{\bf Minimality.}
The presence of two levels in our foundation  facilitates its comparison
with other foundations, given that we can choose the most appropriate level at which
to make the comparison.
To establish the minimality of our foundation we will compare  intensional theories,
such as  
 Martin-L{\"o}f's one or the Calculus of Inductive Constructions, with its intensional level,
while extensional theories, such as the internal calculus of a generic topos
(as devised, for example, in \cite{tumscs})
or  Aczel-Myhill's CZF theory, with its extensional level. Also  logic enriched type theory
in \cite{AG} can  be  compared with our \mtt. Indeed it appears
as a fragment of our \mtt\  except that, being just a  many-sorted logic
 on Martin-L{\"o}f's type theory, its
 propositions are not inhabited with proofs
and they are not seen as collections (or sets) of their proofs as in our \mtt.
In \mtt\ we use such a property  to represent useful constructions on subsets.

Two-level theories, where one level is related to the other
via a quotient completion,  already appeared in the literature.
One of this is Hyland's effective topos~\cite{eff}. There the underlying
theory is given by a tripos~\cite{tripos}, namely
a realizability model of many-sorted intuitionistic
 logic indexed
on classical set theory. Then the topos is obtained by freely adding quotients
to a 
 regular category associated to the tripos~\cite{Car}.
However, the effective topos can be seen as obtained by a quotient completion
on a lex category, too~\cite{Car}. This latter completion is closer to our
quotient completion. The precise correspondence between our quotient completion
and the ones existing in the literature of category theory
 is left
to future work with the study of a general notion
of quotient completion.

{\bf Summary. } The main contributions of this paper are the following:
\begin{list}{-}{ }
\item
We introduce a two-level foundation for constructive mathematics where
both levels are given by type theories  \`a la
 Martin-L{\"o}f: one called \mtt\ is {\it intensional}
and the other called \emtt\ is   {\it extensional}.
They both essentially extend with collections previously introduced
theories respectively in \cite{mtt} and in \cite{mai07}.
The foundation is {\it minimal} among the most relevant known constructive foundations in type theory, or in set theory, or in category theory,
to be compared with ours at the appropriate level.
\item
The extensional level \emtt\  is interpreted in a quotient model \`a la Bishop
built over the intensional one \mtt\  by means of {\it canonical isomorphisms}.
This is because equality of \emtt\ types 
 gets interpreted into an isomorphism of intensional types.
As a consequence, \emtt\ can be viewed as a language {\it to reason within our
quotient model} over \mtt. As an application,
we get that the \emtt-formulation of the axiom of choice 
turns out to be interpreted in \mtt\ as
exactly  Martin-L{\"o}f's extensional axiom of choice in \cite{MLac,Jes},
and hence it is not valid. Even the axiom of unique choice is expected not to be valid in \emtt\
 as advocated in \cite{mtt}.
\item
We adopt an intensional version of \mtt\  {\it without the $\xi$-rule for lambda-terms},
after noticing from  \cite{Ros,bcsr} that its presence is irrelevant
to interpret the extensional level via quotients.
As said in \cite{modmar}, the absence of the $\xi$-rule
opens the way to interpret \mtt\  via Kleene's realizability
interpretation~\cite{DT88}, and hence to show that it satisfies our proofs-as-programs 
requirement of consistency with CT+AC, as required to the intensional
level
according to our notion of {\it constructive two-level foundation}.
This lets us to avoid the problem of proving consistency with  CT+AC
for intensional theories with the $\xi$-rule as in \cite{mtt,PMTT}, which is still open.
\end{list}

\section{The  intensional level \mtt}
\label{intlev}
Here we briefly describe the intensional level of  our two-level constructive foundation.
It consists of an intensional type theory in the style of 
Martin-L{\"o}f's one in
\cite{PMTT}, which essentially extends that presented
in \cite{mtt} with the notion of ``collection''.
Indeed, the version  in \cite{mtt} called Minimal Type theory
(for short \mtt)
 essentially 
corresponds to the set-theoretic part of that presented here, which we
still call \mtt.
We says ``essentially''
since the set-theoretic part of our new \mtt, called here \mtts, has different
rules about equality from that in \cite{mtt}.

We thought of modifying the original calculus in \cite{mtt} for
the following reasons.
First,
we wanted to extend the calculus in \cite{mtt} with collections
and the necessary constructors to support power collections of sets
via quotients. This opens the way 
to formalize various mathematical
theorems where power collections are used, like,
for example, those about formal topology~\cite{Sam2002,gsmin,BP}.  
As a consequence we had to equip the modified calculus with propositions
closed under usual intuitionistic connectives and quantification
over generic collections. We then called {\it small propositions}
those closed only under quantification over sets beside  usual intuitionistic
connectives.
Lastly, we modified the equality rules in order to easily satisfy
the proofs-as-programs paradigm, namely the consistency of \mtt\
 with the axiom of choice
and the formal Church thesis.

More in detail,
the typed calculus \mtt\ is written 
in the style of Martin-L{\"o}f's type theory \cite{PMTT} 
 by means of the following
 four kinds of judgements:

{\small  $$A \ type\ [\Gamma] \hspace{.5cm} A=B\ type\ [\Gamma] 
\hspace{.5cm} a \in A\ 
 [\Gamma] \hspace{.5cm} a=b \in A\ [\Gamma] $$}

\noindent
that is the type judgement (expressing that something is a specific type), 
the type equality judgement (expressing when
two types are equal), the
term judgement (expressing that something is a term of a certain type)
 and the term equality judgement
(expressing the definitional equality between terms of the same type), respectively, all under a context  $\Gamma$.
The contexts $\Gamma$  of these judgements are formed
as in \cite{PMTT} and they are telescopic \cite{DEBRUIJN} since
types are {\it dependent}, namely they are 
 allowed to depend on variables ranging over other types.
The precise rules of \mtt\ are given in the appendix~\ref{mttsyn}.

Types include collections, sets, propositions and small propositions and hence
 the word {\it type} is only used as a meta-variable,
namely 
{\small $$type \in \{ col, set,prop,prop_s\, \}$$}
Therefore, in \mtt\ types are actually formed by using the following
judgements:

{\small  $$A \ set\ [\Gamma] \qquad A\ col\ [\Gamma] 
\qquad A\ prop\
 [\Gamma]\qquad A\ prop_s\
 [\Gamma] $$}

As in \cite{mtt}, the general idea is to define a many-sorted logic, but now sorts
include both sets and collections. 
The main difference between sets and collections is that
{\it sets are those collections that are inductively generated}, namely those
whose most external constructor is
equipped with introduction and elimination rules,
and all of their collection components are so. 
According to this view we will allow elimination rules of sets 
to act also towards collections.

Our sets will be closed under the empty set, the singleton set,
strong indexed sums, dependent products, disjoint sums,
lists. These constructors are formulated as in  Martin-L{\"o}f's type theory
with the modification that their elimination rules vary on all types.
In order to view sets as collections, we add the rule {\bf set-into-col}
{\small
$$\begin{array}{l}
 \displaystyle{ \frac
       {\displaystyle \ A
        \hspace{.1cm} set\  }
{ \displaystyle \ A
        \hspace{.1cm} col\  }}
\end{array}\qquad \qquad$$}
The logic of the theory is described by means of propositions and small propositions. Small propositions are those propositions
closed only under intuitionistic connectives and quantification over sets.
To express that
 a small proposition is also a proposition
we add the subtyping rule {\bf prop$_s$-into-prop}
 {\small
$$
\begin{array}{l}
\displaystyle{ \frac
       {\displaystyle\  A\ prop_s\ }
      { \displaystyle\  A\ prop\  }}
\end{array}\qquad\qquad
$$}

\noindent
As explained in \cite{mtt}, since we restrict our
consideration only to mathematical propositions, it makes
sense to identify a proposition with the collection of its proofs.
To this purpose we add the rule { \bf prop-into-col}
 {\small
$$
\begin{array}{l}
\displaystyle{ \frac
       {\displaystyle\  A\ prop \ }
      { \displaystyle\  A\ col\  }}
\end{array}\qquad\qquad
$$}

\noindent
However, proofs of small propositions
are  inductively generated. Hence,
small propositions, as propositions in  \cite{mtt},
are though of as sets of their proofs by means of  the rule
 {\bf prop$_s$-into-set}
 {\small
$$
\begin{array}{l}
\displaystyle{ \frac
       {\displaystyle\  A\ prop_s\  }
      { \displaystyle\  A\ set \ }}
\end{array}\qquad\qquad
$$}

\noindent
The rules {\bf prop$_s$-into-set} and {\bf prop-into-col} allow us to form
the strong indexed sum of a small propositional function
{\small $\phi(x)\ prop_s\ [x\in A]$}, or simply of a propositional function,
{\small $$\Sigma_{x\in A} \phi(x)$$}

\noindent
both on sets and on collections. Given that we will define a subset 
as the equivalence class of a small
propositional
function, then the {\bf prop$_s$-into-set} rule
is  relevant to turn a small propositional function on a set
into a set, and hence to represent
  functions between subsets as in \cite{toolbox} and
to represent families indexed on a subset as advocated in
  \cite{BP}. The same can be said about subcollections.
Moreover, the identification of a proposition with the collection (or set) of its proofs
allows also to derive all the induction principles for propositions
depending on a set, because set elimination rules
can  act towards all collections including propositions.

In order to interpret the power collection of a set as a  quotient
of propositional functions, to \mtt\
we add  the collection of
 small propositions  $\mathsf{prop_s}$ and the collection of functions
from a set
 towards $\mathsf{prop_s}$. 
Since such function collections towards  $\mathsf{prop_s}$ are instances of dependent product
collections, to easily show some meta-theoretic 
properties about \mtt\ we will consider an extension of \mtt, called \mttdp,
with generic dependent product collections (see the appendix~\ref{mttsyn} for its rules).
Finally, by still keeping the minimality of our
two-level foundation,
we   close \mtt\ collections under strong indexed sums  in order to give a simple categorical interpretation
of the extensional level in the model we will build over \mtt.

It is worth noting that
the subtyping relation of propositions into collections and that 
 of small propositions into sets, via the rules
 {\bf prop-into-col} and   {\bf prop$_s$-into-set}  respectively,
are very different from 
 the subtyping relation of sets into collections via the 
  {\bf set-into-col} rule, or that of small propositions into propositions
via the {\bf props-into-prop} rule.
Indeed, the subtyping rules
 {\bf prop-into-col} and  {\bf prop$_s$-into-set} 
do not affect the elimination rules of propositions and small propositions:
they express  a {\it merely inclusion}.
Instead the subtyping rules {\bf set-into-col} and {\bf props-into-prop}
take 
part in the definition of sets and in that of  small propositions,
because  elimination rules of sets and of 
small propositions  act respectively towards all collections
and towards all
propositions.

 There are important motivations, already explained in \cite{mtt},
 behind the fact that 
 elimination rules of propositions act only towards all
propositions and not towards all collections, as well
as  those of small propositions  do not act towards all sets.

First, as said in \cite{mtt},
 propositions have their own distinct origin,
and only {\it a posteriori} they are recognized as collections of their proofs.
Then, a more technical reason is that we want to prevent the validity of the axiom
of choice. Indeed, as described  in \cite{mtt}, the axiom of choice turns out to be valid
 if we allow elimination rules of propositions towards
all collections or of small propositions towards all sets, because the existential quantification
would be then equivalent to the corresponding strong indexed sum on the same constituents~\cite{luo}.
The reason to reject  the general
 validity of the axiom of choice in \mtt\  is to get a minimalist foundation compatible
with the existing ones, including the internal theory of a generic topos
where the axiom of choice is not always valid.
All this attention to avoiding the validity of the axiom of choice
was paid in \cite{mtt} because there we were trying
to get our minimal foundation by modifying 
 Martin-L{\"o}f's intensional
type theory in  \cite{PMTT} (here called MLTT).
MLTT is not minimal just because it validates the axiom of choice,
given that it
 follows the isomorphism
``propositions as sets'' and hence it identifies the existential quantifier with
the strong indexed sum. To discharge such an isomorphism, and hence the validity of the axiom of choice, it is sufficient
to introduce a primitive notion of propositions with the mentioned restrictions
on their elimination rules. As  result of this process
 our \mtt\ version with collections, as well as that in \cite{mtt}, can be  naturally embedded in
Martin-L{\"o}f's intensional type theory \cite{PMTT},
if we interpret   sets as sets in a fixed  universe,
for example in the first universe $U_0$ in \cite{PMTT}, and collections
as generic sets. Then, propositions are interpreted as sets, always by following
the isomorphism ``propositions as sets'', and,
 after identifying small propositions with sets in $U_0$,
the collection of small proposition will be
of course interpreted as $U_0$ itself.

In \cite{mtt} we chose to work with MLTT
  in order
 to get a minimalist proofs-as-programs foundation for its intensionality. Indeed, 
 extensional theories, such the internal calculus
of a topos or Aczel-Myhill's CZF, can not satisfy  our ``proofs-as-programs''
requirement (see \cite{mtt}).
Actually, in \cite{mtt} we designed a version of \mtt\  
of which we still do not know 
whether it satisfies our proofs-as-programs
requirement, namely whether it is consistent with  the axiom of choice
and the formal Church thesis. This problem can be reduced to asking
whether intensional Martin-L{\"o}f's type theory
in \cite{PMTT} is consistent with the formal Church thesis.
 Here, we do not solve such problems
but we adopt a version of \mtt\  that hopefully satisfies our
requirement by means of Kleene's realizability interpretation~\cite{DT88}.

We got to this version after a suggestion 
 by T. Streicher and P. Martin-L{\"o}f  already reported in \cite{modmar}.
There it is said that  Kleene's realizability interpretation of set-theoretic
constructors validates
the first order version of Martin-L{\"o}f's type theory
with explicit
substitution rule for terms
{\small $$\begin{array}{l}
      \mbox{sub)} \ \
\displaystyle{ \frac
         { \displaystyle 
\begin{array}{l}
 c(x_1,\dots, x_n)\in C(x_1,\dots,x_n)\ \
 [\, x_1\in A_1,\,  \dots,\,  x_n\in A_n(x_1,\dots,x_{n-1})\, ]   \\[2pt]
a_1=b_1\in A_1\ \dots \ a_n=b_n\in A_n(a_1,\dots,a_{n-1})
\end{array}}
         {\displaystyle c(a_1,\dots,a_n)=c(b_1,\dots, b_n)\in
 C(a_1,\dots,a_{n})  }}
\end{array}
$$}

\noindent
in place of  the  usual  term equality rules in \cite{PMTT}.
 Hence this modified first-order version of MLTT is 
consistent with  the formal Church thesis and, in turn, satisfies our
proofs-as-programs requirement, given that the axiom of choice
is a theorem there.
Instead,
the version with the original term equalities in  \cite{PMTT} is not validated
 by Kleene's realizability interpretation just because of
the presence of the $\xi$-rule
{\small $$
\mbox{\small  $\xi$} \
\displaystyle{\frac{ \displaystyle c=c'\in  C\ [x\in B]  }
{ \displaystyle \lambda x^{B}.c=\lambda x^{B}.c' \in \Pi_{x\in B} C}}
$$ }

Now, given that, as we say above,
we can interpret our version of  \mtt\ in Martin-L{\"o}f's type theory  with at least one universe,
 if we perform the same change of equality rules for \mtt\ then 
Kleene's realizability interpretation surely validates the set-theoretic version of \mtt \
and hopefully also the whole version~\footnote{The realizability 
interpretation of  the collection of small
propositions, as well as of the first universe in Martin-L{\"o}f's type theory in  \cite{PMTT},
is a delicate point. We expect to interpret it  as the subset  of
 codes of small propositions.},
by providing  a  proof of \mtt-consistency   with the axiom of choice and
 formal Church thesis  (the axiom of choice
holds because the realizability interpretation interprets
 the existential quantifier as the strong indexed sum).

\noindent
Luckily, we can take this modified version of \mtt\ without the  $\xi$-rule
as the {\it intensional level} of our desired two-level foundation
with no effect on the interpretation of its extensional level.
Indeed, the above change of term equality rules
 does  not  affect the properties of the quotient
model we are going to build over \mtt\ and where we will
interpret \emtt.
 This was first noticed in categorical terms
 in \cite{Ros,bcsr} (see remark~\ref{ros}).

Note that our \mtt\ does not include the boolean universe
$U_b$ used in  \cite{mai07} to derive disjointness of binary sums.
Indeed, in the presence of a collection of small propositions
this is derivable anyway (see the proof of theorem~\ref{mainth1}).

Our present version of \mtt, as that in \cite{mtt}, 
can be still though of as  a predicative version of the Calculus of Inductive Constructions based on \cite{TC90} with types and propositions: \mtt\ collections and sets are simply interpreted as types,
propositions as themselves, and the
 collection of small propositions as the type of all
propositions. 

Of course, it makes sense to compare \mtt\ just with  intensional
type theories as those already mentioned.

\section{The extensional level \emtt}
\label{ext}
Here we briefly describe the extensional level of our
desired two-level foundation. This is the level where we will actually formalize
constructive mathematics.

As well as the intensional level, it consists of a
 type theory, called \emtt,  written now in the style
of 
  Martin-L{\"o}f's {\it extensional} type theory in \cite{ML84}.

The main idea behind the design of \emtt\
goes back to that of toolbox in \cite{toolbox}, and it follows
the {\it forget-restore
principle} conceived by G. Sambin and introduced there.
According to this principle,
extensional concepts must be  abstractions of
 intensional concepts as result of forgetting
 irrelevant computational information, that can be restored
when implementing these extensional concepts back at the intensional level.
Here, we think it is enough to require that extensional constructors
must be obtained by abstracting only over equalities
 of corresponding constructors in \mtt. Hence, it seems sufficient
to add quotients to \mtt\ to be able to represent our extensional
level. This leads us to conclude that the extensional level
should be considered as a fragment of the internal theory of
 a model built over \mtt\ by just adding quotients.

 In the literature
it is well known how to add quotients to a type theory
by building  setoids (see \cite{disttheshof,ven}) on it.
The extensional theory  \emtt\ we propose here 
can be interpreted in a suitable model of total setoids \`a la Bishop
that we will describe in the following.

We need to warn that \emtt\ is not precisely the {\it internal language} of 
the quotient  model we adopt, namely it is not fully complete with it, or
in other words
it does not necessarily capture all the constructions valid in the model,
especially at the level of collections.
For example, in \emtt\  collections are not closed
under effective quotients while in the model they are, instead.
Hence, different extensional levels may be considered over \mtt,
and even over the same quotient model.
A criterion to decide what to put in \emtt\ 
at the level of collections is that of preserving
minimality, for example with respect
to Aczel's CZF.


\emtt\ extends the set-theoretic version introduced in \cite{mai07}, called here
\emtts,
with collections and related constructors needed to represent
the power collection of
a set
with  $\varepsilon$-relation and
comprehension used in everyday mathematical practice.

\emtt\ is essentially obtained as follows: we  first take the {\it extensional} version of \mtt,
 in the same way the so-called extensional
 Martin-L{\"o}f's type theory in   \cite{ML84} is obtained from the intensional version in \cite{PMTT},
with the warning of replacing the collection of small propositions 
{\small $\mathsf{prop_s}$} with
 its quotient  
under equiprovability  {\small ${\cal P}(1)\, \equiv\, \mathsf{prop_s}/\leftrightarrow$};
then we collapse propositions into {\it mono collections} according to the notion in \cite{tumscs};
and finally we add effective quotient sets as in \cite{tumscs}.
The precise rules of \emtt\  are given in the appendix~\ref{emtt}.
The form of judgements to describe \emtt\ are those of \mtt.


\noindent
One of the main differences between \emtt\ and \mtt\ amounts to being that
between the extensional version of Martin-L{\"o}f's type theory  in \cite{ML84}
 and the intensional one in \cite{PMTT}. It consists in the fact  that
while
  type judgements in the intensional version
are decidable, those in the extensional one are no longer so~\cite{disttheshof}.
Another difference is that in \emtt\
 propositions are {\it mono} as in  \cite{tumscs}, that is 
they
are inhabited by at most one proof by introducing in \emtt\  the following rule:
{\small
$$
\mbox{{\bf prop-mono}) }  \,\displaystyle{ \frac{ \displaystyle A\ prop\ [\Gamma]\qquad p\in A\  [\Gamma] \quad q\in A\  [\Gamma]}{\displaystyle p=q\in A\ \ [\Gamma]}}  $$}
Propositions are then {\it mono collections} and small propositions are 
{\it mono sets}.
This property allows us to forget proof-terms of propositions,
namely to make {\it proofs of propositions irrelevant},
by introducing a canonical proof-term called $\mathsf{true}$  for them:
{\small 
$$
\begin{array}{l}
\mbox{{\bf prop-true}) }  \,\displaystyle{ \frac{ \displaystyle A\ prop\ \qquad p\in A\  }
{\displaystyle \mathsf{true}\in A}} 
\end{array}
$$}
This canonical proof-term allows us to interpret true-judgements
 in  \cite{ML84,sienlec} directly in \emtt\  as follows:
{\small
$$A\ true\ [\, \Gamma; B_1\ true,\dots, B_m\ true\, ]\, \equiv\, \mathsf{true}\in A\ [\, \Gamma,\,  y_1\in B_1,\dots, y_m
\in B_m\, ]
$$}

\noindent
Then, we can prove that, according to this
interpretation, all true judgements of the logic in \cite{ML84,sienlec} are valid in \emtt.

A key feature of extensional type theory in \cite{ML84}
is the presence of {\it Extensional Propositional Equality}, written
{\small $\mathsf{Eq}(A,a,b)$} to express that $a$ is equal to $b$.
This  is stronger than  Propositional Equality {\small $\mathsf{Id}(A,a,b)$}
in intensional type theory~\cite{PMTT}, and also in
\mtt\ (see appendix~\ref{mttsyn}), because the
 validity of {\small $\mathsf{Eq}(A,a,b)$}
 is equivalent to the definitional equality of terms
{\small $a=b\in A$} (both under the same context). Furthermore,
it   is also mono.
We add {\small $\mathsf{Eq}(A,a,b)$} to \emtt\
as a proposition, which is small when $A$ is a set
(see  appendix~\ref{emtt}).

Another key difference between \emtt\ and \mtt\  is that in \emtt\
 we can form
 effective quotient sets (see \cite{tumscs}). Then, 
in the presence of quotient effectiveness it is crucial to require that
 propositions are mono collections, as well as that
small propositions are mono sets.
Indeed, if we identify small propositions with sets simply, 
or propositions with collections,
quotient effectiveness 
 may lead to classical logic (see \cite{Maieff}), because it yields
to a sort of  choice operator,
 and hence it is no longer a constructive rule.

Moreover, observe that  the
 set-theoretic part of our \emtt, called \emtts,
is a variation of the internal type theory of a list-arithmetic
locally cartesian closed pretopos, as devised in \cite{tumscs}.
Indeed, \emtts\ is not exactly that because in the internal
theory of a generic pretopos
{\it small propositions are identified with mono sets}, while in 
\emtt\ {\it small  propositions are only some primitive mono sets},
and it does not necessarily follow that
that all mono sets are small propositions.
In this way we avoid the validity of {\it the axiom of unique choice}, which
 would instead be valid under the identification of  small propositions with mono sets (see \cite{tumscs}).

Lastly, 
 in \emtt\ we add the necessary constructors
to represent 
the power collection of  a set as the
``power set''
 in the internal theory of a topos devised in
the style of   Martin-L{\"o}f's type theory in \cite{tumscs}.
To this purpose we put in \emtt\ the power collection  of the singleton set
 {\small ${\cal P}(1)$}. This is represented as 
the quotient of the collection of small propositions  
 under equiprovability. Therefore a subset is represented
as an equivalence class of small propositions.
Then, we add collections of functions towards  {\small ${\cal P}(1)$}
in order to represent the power collection of a set
$A$:
 {\small $${\cal P}(A)\, \equiv\, A\rightarrow {\cal P}(1) $$}

\noindent
 Therefore, a subset of $A$, being  an element of {\small  ${\cal P}(A)$},
is represented as a function from $A$ to  {\small ${\cal P}(1) $}.

\noindent
Moreover, we can represent  $\varepsilon$-relation and
comprehension used in everyday mathematical practice as follows.
To this purpose we need to assume the extensional equality
$\mathsf{Eq}( {\cal P}(1),\, U, \, V\, )$ on ${\cal P}(1)$
to be a small proposition given that  the equality on 
${\cal P}(1)$ is
 intended to be defined
as equiprovability of the small propositions characterizing the equivalence classes $U,V$,
and this is still a small proposition.
Then, given
{\small
$W\in  {\cal P}(A)$} and {\small $a\in A$} we  define the small proposition
{\small $$a\varepsilon W\, \equiv\, \mathsf{Eq}( {\cal P}(1),\, W(a), \, [\mathsf{tt}]\, )$$}

\noindent
where {\small $\mathsf{tt} \, \equiv\, \bot \rightarrow \bot$} is
the truth constant (it may be represented by
 any tautology),
 and  {\small $[\mathsf{tt}]$}  its equivalence class under equiprovability.

\noindent
Furthermore, for any derivable
  {\small  $B(x)\in \mathsf{prop_s}\ [x\in A]$ }
we define
 {\small  $$ \{\, x\in A\, \mid\, B(x)\, \}\in {\cal P} (A)\, \equiv\,
\lambda x^A. [B(x)]\in  {\cal P} (A)$$}

\noindent
Then, thanks to  the rules eq-${\cal P}$) and  eff-${\cal P}$) 
of  {\small${\cal P}(1)$} and the equality rules expressing
 extensionality of function collections  in  appendix~\ref{emtt},
 we can prove that the equality between subsets determined
by propositional functions is the usual extensional one as in \cite{toolbox}:
{\small $$ \{\, x\in A\, \mid\, B(x)\, \}\, = \,  \{\, x\in A\, \mid\, C(x)\, \}\, \in  {\cal P}(A) \mbox{ holds in \emtt\ }\mbox{  iff }\forall_{x\in A} \, B(x)\leftrightarrow C(x)
 \mbox{ holds in \emtt.}$$}

\noindent
Furthermore, also the comprehension axiom 
{\small $$a\varepsilon \{ x\in A\, \mid \, B(x)\, \} \leftrightarrow
B(a)$$}

\noindent
holds in \emtt\  because
{\small $a\varepsilon \{ x\in A\, \mid \, B(x)\, \}$} is equal to
{\small $\mathsf{Eq}( \, {\cal P}(1),\,[B(a)]\, , \, [\mathsf{tt}]\, )$},
which is valid if and only if $B(a)$ is valid, too, again by the  rules eq-{\cal P}) and  eff-{\cal P})
of {\small ${\cal P}(1)$} and those of $\mathsf{Eq}$ in  appendix~\ref{emtt}.

\noindent
Thanks to the rule $\eta$-{\cal P}) and those about function collections,
 we can prove that for any subset {\small$W\in {\cal P}(A)$} we can derive
{\small$$W=\{ \, x\in A\ \mid\ x\varepsilon W\, \}\in  {\cal P}(A)$$}

\noindent
Alternative  rules to form the power collection of a set
can be deduced from the analysis in \cite{tumscs}
on how to represent the subobject classifier
 in the internal type theory of a topos.

Finally, it is worth noting that,
 thanks to  proof irrelevance of  propositions,
 we can  implementing functions
between  subsets as in \cite{toolbox}
without running into the problem pointed out in \cite{carl}.

The desire of {\it representing power collections of sets} together with  {\it proof irrelevance
of propositions} is a key motivation
to work in a two-level foundation. Indeed, such constructions
 can not be directly represented in \mtt\
because \mtt\ has no quotients. On the other hand it is sufficient
to build a quotient model over \mtt\ to interpret them, and to interpret
the whole \emtt.

In order to better present the proofs about this quotient model
over \mtt,   
 it is more convenient to show them for an extended two-level 
theory. The extensional level of this extended two-level theory is taken to be
 an extension of \emtt, called \emttdp.
\emttdp\ is first obtained by extending \emtt\ 
 with dependent product collections, 
given that function collections towards {\small ${\cal P}(1)$} can be seen
as  instances of them 
(see the rules in  appendix~\ref{emtt}).
Consequently, the intensional level is necessarily taken to be \mttdp, namely
\mtt\  extended with
dependent product collections, in order
to support the interpretation of corresponding
collections in \emttdp.
Moreover,
in \emttdp\ we also  include effective quotients on collections,
(that we do not include in \emtt!)
in the attempt to capture  the largest extensional theory 
\`a la Martin-L{\"o}f  valid in our quotient model over \mttdp.
However \emttdp\ does  not seem to be the internal language
of our quotient model over  \mttdp\ yet, as well as \emtt\
is not at all that of our  quotient model over \mtt. 
  We leave 
to future work  how to determine the fully complete theory
of our quotient model over \mttdp, and eventually that over \mtt.

\noindent
It is worth noting that our quotient models over \mtt\ and  over \mttdp\ 
support also the interpretation of the collection
of small propositions {\small  $\mathsf{prop_s}$}, if added
to \emtt\ and to \emttdp\ respectively. Indeed,  {\small  $\mathsf{prop_s}$}
 turns out to be
interpreted as   {\small ${\cal P}(1)$},
  because equality of \emtt\ propositions
will be interpreted into equiprovability of \mtt\ propositions.
We chose of putting in \emtt\ only the quotient
collection of small propositions under equiprovability  {\small ${\cal P}(1)$},
and not the collection of small propositions {\small $\mathsf{prop_s}$}, 
in order
 to make \emtt\ easily interpretable in the internal theory
of a topos. Indeed, if we consider
the formulation ${\cal T}_{top}$
of the internal theory of a topos devised in \cite{tumscs}
in the style of  Martin-L{\"o}f's extensional type theory,
 then \emtt\ sets and collections are translated into ${\cal T}_{top}$ types,
\emtt\ small propositions and propositions into ${\cal T}_{top}$ mono types,
that represent
 ${\cal T}_{top}$ propositions,  and  {\small ${\cal P}(1)$ }
  is translated as the ${\cal T}_{top}$ type representing
the subobject classifier. 

Note that, as in \mtt\, also in \emtt\ we do not include the boolean universe
$U_b$ used in  \cite{mai07} to derive disjointness of binary sums.
Indeed, in the presence of  {\small  ${\cal P}(1)$ }
sum disjointness is derivable anyway~\footnote{The proof of disjointness
in \emtt\ is  similar to that
for \mtt\ mentioned in the proof of theorem~\ref{mainth1}.}.

Our \emtt\ is also compatible with the notion of 
a predicative topos~\cite{swopos}, being its set-theoretic part a fragment
of the internal theory of a locally cartesian closed pretopos.
In particular, if 
we perform over Martin-L{\"o}f's type theory with universes
the  quotient model we will build over \mtt\ to interpret \emtt,
then we will get 
a predicative topos (see \cite{mptt,swopos}).
Given that our intensional level
 \mtt\ can be interpreted in Martin-L{\"o}f's type theory
with universes,
this yields that \emtt\ can be also interpreted in the predicative topos
over it.

Finally, \emtt\ is certainly compatible with Aczel's CZF~\cite{czf} 
by interpreting sets a CZF sets, collections as classes,
 propositions as subclasses of the singleton and
small propositions as subsets of the singleton (in order to make the
rules {\bf prop-into-col} and {\bf prop$_s$-into-set}  valid).
In particular,
the power collection  {\small ${\cal P}(A)$} of a set $A$ is interpreted as the corresponding
power collection of subsets.

\section{The quotient model}
In what follows we are going to define a model of quotients over \mtt.
This is
based on the well-known notion
of setoid~\cite{Bishop,disttheshof}, namely a set with an equivalence relation over it, that we apply
to any type we consider.
We will then interpret \emtt\ in such a model.

The model will be presented in a categorical shape. 
The reason is that we know the categorical semantics
of \emtt\  and, if the quotient model turns out to be an instance of it,
then  we can interpret \emtt\ into the model, and hence into \mtt.
 This model represents also a way to freely add quotients
to \mtt\ in a suitable sense. The precise categorical formulation of its universal property
is left to future work.

\noindent
We first start with defining the  category of ``extensional collections''
namely collections equipped with an 
equivalence relation.
In the case where the collection is a set, this notion
is known as  ``total setoid''
 in the literature~\cite{disttheshof,ven}.

\begin{definition}
\label{setoid}
\rm
The category {\small {\bf \cq}} is defined as follows:\\
{\small \bf $\mathsf{Ob}$\cq}:  objects  are pairs {\small $(A, =_{A})$} where
$A$ is a collection in \mtt, called ``support'', and 
{\small $$x=_{A}y\ prop\ [x\in A, y\in A]$$}

\noindent
is an equivalence relation on the collection $A$. This means that in \mtt\ there exist proof-terms
witnessing reflexivity, symmetry and transitivity of the relation: 
{\small $$\begin{array}{l}
\mathsf{rfl}(x) \in  x=_{A}x\ \
  [\, x\in A\, ]
\\
\mathsf{sym}(x,y,u) \in  y=_{A}x \  \ [\, x\in A,\, y\in A,\  u\in x=_{A}y\, ]\\
\mathsf{tra}(x,y,z,u,v) \in  x=_{A}z\ \ [\, x\in A,\, y\in A,\,  z\in A,\  u\in 
  x=_{A}y,\  v\in  y=_{A}z\, ]\end{array}$$}

\noindent
We call {\small $(A, =_{A})$} {\it extensional collection}.

\noindent
{\small \bf {$\mathsf{Mor}$Q(\mtt)}}:  morphisms  from an object {\small $(A,=_{A})$} to {\small$(B,=_{B})$}
are \mtt\ terms $f(x)\in B\ [x\in A]$ preserving the corresponding equality, i.e. in \mtt\  there exists a proof-term
{\small $$\mathsf{pr}_1(x,y,z)\in  f(x)=_{B}f(y)\ \  [\, x\in A,\, y\in A, \ z\in x=_{A}y\, ]     $$}

\noindent
Moreover, two morphisms {\small $f,g: (A,=_{A})\rightarrow  (B,=_{B})$}
are equal if and only if in \mtt\  there exists  a proof-term
{\small $$\mathsf{pr}_2(x)\in   f(x)=_{B}g(x)\ \  [\, x\in A\, ]     $$
}
\end{definition}


The category 
{\small \cq}
comes naturally equipped with an indexed category (or split fibration) satisfying the universal property of comprehension (see
\cite{jacobbook} for its definition) thanks to the closure of \mtt\ collections under strong indexed sums:
\begin{definition}\em
The indexed category: 
{\small $$\mathcal{P}_q: \mbox{\rm Q(\mtt)}^{\mathsf{OP}}\rightarrow {\mathsf{Cat}}$$}

\noindent
is defined as follows.  For each object {\small$(A,=_{A})$} in {\small\mbox{\rm Q(\mtt)}\ }then
{\small$\mathcal{P}_q (\, (A,=_{A})\, ) $} is the following category:
{\small $\mathsf{Ob} {\mathcal P}_q (\, (A,=_{A})\, )$} are the
propositions {\small $P(x)\ prop\ [x\in A]$} depending on $A$ and
preserving the equality on $A$, namely for such propositions there exists a proof-term:
{\small $$\mathsf{ps}(x,y,d)\in P(x)\rightarrow P(y)\ [x\in A,\, y\in
A,\, d\in x=_{A}y]~\footnote{ Indeed, from this, by using the symmetry
of $ x=_{A}y$ it follows that $ P(x)$ is equivalent to $P(y)$ if
$x=_{A}y$ holds.}$$}

\noindent
Two objects  {\small $P(x)\ prop\ [x\in A]$} and {\small $Q(x)\ prop\ [x\in A]$}
are equal if a proof of {\small $P(x)\, \leftrightarrow \, Q(x)\ prop\ [x\in A]$} can be derived in \mtt.

\noindent
Morphisms in {\small $\mathsf{Mor} {\mathcal P}_q (\, (A,=_{A})\, )$ }
are given by the following order
 {\small
$$\begin{array}{l} {\mathcal P}_q (\, (A,=_{A})\, ) (\, P(x)\, ,
\,Q(x)\, ) \, \equiv\, P(x)\leq Q(x)\\[10pt] \qquad \mbox{ iff }\mbox{
there exists a proof-term } \mathsf{pt}(x) \in P(x)\rightarrow Q(x)\
[x\in A]
\end{array}
$$}

\noindent
 Moreover, for every morphism {\small $f:(A,=_{A})\rightarrow
(B,=_{B})$} in {\small \cq } given by {\small $f(x)\in B\ [x\in A]$}
then $\mathcal{P}_q (f)$ is the substitution functor, i.e.  {\small
$\mathcal{P}_q (f) (\, P(y)\, )\, \equiv\, P(f(x)) $} for any
proposition {\small $P(y)\ prop\ [y\in B]$} (recall that ${\mathcal
P}_q$ is contravariant).
\end{definition} 

\begin{lemma}
\label{ind}
${\mathcal{P}_q}$ is an indexed category satisfying the universal property of comprehension.
\end{lemma}
{\bf Proof.}  To describe the universal property of comprehension, we consider
the Grothendieck completion {\small $ Gr(\mathcal{P}_q)$} of
$\mathcal{P}_q$ (see \cite{jacobbook} for its definition) and the
functor {\small $T:\cq\rightarrow Gr(\mathcal{P}_q)$ } defined as
follows: {\small $$T(\, (A,=_{A})\, )\, \equiv \, ( \, (A,=_{A})\, ,
\, \mathsf{tt}\, )\qquad \qquad T (f)\, \equiv\, (f,
id_{\mathsf{tt}})$$} 

\noindent
where $\mathsf{tt}$ is the truth constant. $T$ satisfies the following universal property of comprehension
(cfr. \cite{jacobbook}):
for every  $Gr(\mathcal{P}_q)$-object {\small $ (  \, (A,=_{A})\, , \, P \, )$} there exists a \mbox{\rm Q(\mtt)}-object
{\small $$Cm (\,  (  \, (A,=_{A})\, , \, P \, )\, )\, \equiv \, (\Sigma_{x\in A}P(x)\, ,\, =_{Cm}\, )$$}

\noindent 
where {\small $z_1=_{Cm}z_2  \, \equiv \, \pi_1(z_1)\, =_{A} \pi_1(z_2)$}
for {\small $z_1,z_2\in \Sigma_{x\in A}P(x)$}, such that, for each $Gr(\mathcal{P}_q)$-morphism {\small $ (f,p):  T(\, (C,=_{C}) \, )\, \longrightarrow\,  (  \, (A,=_{A})\, , \, P \, ) $}  there is a {\it unique  morphism} {\small $ [f,p]: 
 (C,=_{C}) \,\longrightarrow \,Cm (\,  (  \, (A,=_{A})\, , \, P \, )\, )$} in
 \mbox{\rm Q(\mtt)}
such that in $Gr(\mathcal{P}_q)$
{\small $$ (\eta_1,\eta_2)\, \cdot\,  T(\, [f,p]\, ) =(f,p)$$}

\noindent
where {\small $\eta_1\, \equiv\, \pi_1^P$} and {\small $\eta_2\, \equiv\, \pi_2^P$}
are the first and second projections of the indexed sum {\small $\Sigma_{x\in A}P(x)$}.

Note that $Cm(-)$ on objects does not define an operation from  $Gr(\mathcal{P}_q)$-objects to $\mbox{\rm Q(\mtt)}$-objects unless
 the Axiom of Choice is used.

\medskip

\begin{definition} [$\mathcal{P}_q$-proposition]
\label{catprop}\em
Given a proposition {\small $P\in \mathsf{Ob}\mathcal{P}_q( \, (A,=_{A})\, )$}
then the first component  of the comprehension adjunction counit
{\small 
$(\eta_1,\eta_2): \, ( \, (\Sigma_{x\in A}P(x)\, ,\, =_{Cm}\, )\, , \, \mathsf{tt}\, )\rightarrow (\, (A,=_A)\, ,\, P \, )$}
given by the first projection 
{\small $$\eta_1\, \equiv\, \pi_1^P:  (\Sigma_{x\in A}P(x)\, ,\, =_{Cm}\, )\rightarrow  (A,=_A)$$}
 is a monic morphism in $\mbox{\rm \cq}$ and it is called a 
{\it ${\mathcal{P}_q}$-proposition}.
\end{definition}

\begin{definition} [$\mathcal{P}_q$-small proposition]
\label{catprops}\em
A $\mathcal{P}_q$-proposition 
{\small $\eta_1\, \equiv\,\pi_1^P:  (\Sigma_{x\in A}P(x)\, ,\, =_{Cm}\,
  )\rightarrow  (A,=_A)$} is called a {\it ${\mathcal{P}_q}$-small proposition}
if   the  proposition 
{\small $P\in \mathcal{P}_q( \, (A,=_{A})\, )$} is small.
\end{definition}



\noindent
In order to make clear how the set-theoretic part of {\small \cq\ }
can interpret the set-theoretic fragment
 \emtts\
of \emtt, we  single out the category
of extensional sets from that of extensional collections. 

\begin{definition}\em
The category  {\small \cqset} is defined as the full subcategory of {\small \cq\ }
equipped with extensional sets {\small $(A,=_{A})$} where the support
is a set and the equivalence relation is small, namely
{\small $A \ set$} 
and {\small $x=_{A}y \ prop_s\ [x\in A,y\in A] $}
are derivable in \mtt.
\end{definition}

\begin{definition}\em
We define the functor
{\small ${\mathcal{P}_{q}}_{set}:
 \cq_{\mathsf{set}}^{\mathsf{OP}}\rightarrow {\mathsf{Cat}}$}
 as the restriction of {\small ${\mathcal{P}_q}$} on {\small  \cqset}
 to small propositions, namely for every extensional set {\small $(A,=_{A})$}
we have that
{\small  ${\mathcal{P}_{q}}_{set}(\, (A,=_A)\,)$} is the full subcategory of
 {\small $\mathcal{P}_{q}(\, (A,=_A)\,)$} containing only small propositions.
\end{definition}
Then, in an analogous way to lemma~\ref{ind} we can prove:
\begin{lemma}
${\mathcal{P}_{q}}_{set}$ is an indexed category 
satisfying the universal property of comprehension.
\end{lemma}
Analogously to definition~\ref{catprop},
 we can define a {\it ${\mathcal{P}_{q}}_{set}$-proposition},
which is indeed also small.
\\


\noindent
In order to prove that the category {\small \cq\ } has all the necessary structure to interpret \emtt,
in particular that it is closed under certain dependent products,
 it will be useful to know that
{\small \cq}-morphisms correspond to  extensional
dependent collections
defined in an analogous way to  dependent sets in 
\cite{Bishop,notepal} (see also \cite{dyb}) as follows:
\begin{definition}[extensional dependent collection]
\label{ds}
\rm
Given an object {\small $(A,=_{A})$ }of \cq, abbreviated with $A_{=}$,
we define an {\it extensional dependent collection} on {\small $(A,=_{A})$} written 
{\small $$B_{=}(x)\ [\, x\in A_{=}\, ]$$}

\noindent
as a dependent collection
{\small $B(x)\ col\ [\, x\in A\, ]$}, called ``dependent support'',
 together with an equivalence relation
{\small  $$y=_{B(x)}y'\ prop\ [x\in A, \, y\in B(x),\, y'\in B(x)].$$}

\noindent
Moreover, for any $x_1,x_2\in A$ such that
$x_1=_{A}x_2$ holds
 there must exist a {\em substitution morphism}
{\small $$\sigma_{x_1}^{x_2}(d,y)\in B(x_2)\ [x_1\in A,\,
 x_2\in A,\, d\in x_1=_{A}x_2,\,  y\in B(x_1)]$$}

\noindent
 preserving the equality on $B(x_1)$, namely there exists a proof of
{\small $$\begin{array}{l}
\sigma_{x_1}^{x_2}(d,y)=_{B(x_2)}\sigma_{x_1}^{x_2}(d,y')\  prop\ [x_1\in A,\,  x_2\in A,\, d\in x_1=_{A}x_2,\\[3pt]
\qquad \qquad \qquad \qquad \qquad \qquad \qquad \qquad y\in B(x_1),\,  y'\in B(x_1),\,  w\in y=_{B(x_1)}y' ].\end{array}$$}

\noindent
and
{\it not depending} on $d\in x_1=_{A}x_2$
in the sense that we can derive a proof of

{\small $$\sigma_{x_1}^{x_2}(d_1,y)=_{B(x_2)}\sigma_{x_1}^{x_2}(d_2,y)
\ prop\ [x_1\in A,\,
 x_2\in A,\, d_1\in x_1=_{A}x_2,\,  d_2\in x_1=_{A}x_2,\,y\in B(x_1)]$$}

\noindent
Furthermore, $\sigma_{x}^{x}(\mathsf{rfl}(x),-)$ is the identity, namely there exists a proof of
{\small $$\sigma_x^x(\mathsf{rfl}(x),y)=_{B(x)}y\  prop\ [x\in A, \,
y\in B(x)]$$}

\noindent
and the $\sigma_{x_1}^{x_2}$'s are closed under composition,
 namely there exists a proof of

{\small $$\begin{array}{l}
\sigma_{x_2}^{x_3}(\, d_2\, ,\, \sigma_{x_1}^{x_2}(d_1,y)\, )=_{B(x_3)}\sigma_{x_1}^{x_3}
(\, \mathsf{tra}(x_1,x_2,x_3,d_1,d_2)\, ,\,  y\, )\ prop\ \\[3pt]
\qquad \qquad
 [x_1\in A,\,  x_2\in A,\, x_3\in A,\, y\in B(x_1),\, d_1\in x_1=_{A}x_2,\,
d_2\in x_2=_{A}x_3].\end{array}$$}

\noindent
Note that the no dependency on the proofs of the extensional collection
equality $=_A$ allows to use the following  abbreviations:
{\small $$\sigma_{x_1}^{x_2}(y)\, \equiv \sigma_{x_1}^{x_2}(d,y)\qquad
$$} 

\noindent
for {\small $x_1\in A,\,
 x_2\in A,\, d\in x_1=_{A}x_2,\, y\in B(x_1)$}.

Categorically speaking,
we can see that
an extensional dependent collection is given by a functor from
a suitable groupoid category to {\small \cq.}
Indeed, for any extensional type {\small $(A,=_A)$} we can define
the category 
 {\small ${\cal G}(\, (A,=_A)\, )$}  as follows: its objects
are the elements of $A$ (with their definitional equality as equality)
and for {\small $a_1,a_2\in A$} then {\small $a_1\leq a_2$} holds if  and only if we can derive {\small $p\in a_1=_A a_2$} in \mtt.
Then, any extensional dependent collection on  {\small $(A,=_A)$}
is given by a functor 
{\small 
$$\sigma_B: {\cal G}(\, (A,=_A)\, )\rightarrow\cq\qquad \mbox { such that }\qquad
\sigma_B(x)\, \equiv\, B_=(x)\qquad \sigma_B(x_1\leq x_2)\,\equiv\,
\sigma_{x_1}^{x_2}$$}
 
\noindent
where $\sigma_{x_1}^{x_2}$ stands for the
extensional morphism given by
{\small $$\sigma_{x_1}^{x_2}(y)\in B(x_2)\ [x_1\in A,\,
 x_2\in A,\, d\in x_1=_{A}x_2,\,  y\in B(x_1)]$$}
Note that {\small ${\cal G}(A,=_A)$ } is a groupoid category and hence
its image along $\sigma_B$ is a groupoid category, too.
In particular we get that
 every $\sigma_{x_1}^{x_2}$ actually gives rise to an {\it isomorphism}
between  {\small $(\, B(x_1)\, , \, =_{B(x_1)}\, )$} and {\small 
$(\, B(x_2)\, , \, =_{B(x_2)}\, )$} for given {\small $x_1,x_2\in A$} and 
{\small $d\in a_1=_A a_2$}.
\end{definition}

\noindent
Analogously, we can give the definition of extensional dependent set:

\begin{definition}\label{dset}
\rm
An extensional dependent collection
{\small $B_{=}(x)\ [\, x\in A_{=}\, ]$} is an {\it extensional dependent set}
if its support is a dependent set and its equivalence relation is small, namely we can derive
{\small $$B(x)\ set\ [\, x\in A\, ]\qquad\qquad y=_{B(x)}y'\ prop_s\ [x\in A, \, y\in B(x),\, y'\in B(x)]$$}
\end{definition}
Analogously, we can give the definition of 
{\it extensional dependent proposition} and of {\it extensional
dependent small proposition}, where the latter must be also
an extensional dependent set.
In the following we will often speak of extensional dependent type to include any of them.

\begin{definition}
\em \label{extcat}
Let us call {\small $Dep_{\cq}(\, (A, =_{A})\, )$} the category whose objects are extensional dependent
 collections {\small $B_{=}(x)\ [x\in A_{=}]$}
on  the extensional collection {\small $(A,=_{A})$}, and whose
morphisms are 
{\it extensional terms}
{\small $$b(x,y)\in B_{=}(x)\ [\, x\in A_{=},\, y\in C_=(x)\, ]$$}

\noindent
that are dependent terms
{\small $b(x,y)\in B(x)\ [\, x\in A,\, y\in C(x)\, ]$}
preserving the equality on $A$ and that on $C(x)$, namely in \mtt\  there exists a proof of
{\small $$\begin{array}{l}
\sigma_{x_1}^{x_2}(\, b(x_1,y_1)\, )=_{B(x_2)}
\,  b(x_2,y_2)\ prop\  [\, x_1\in A,\, 
x_2\in A,\, w\in x_1=_{A}x_2,\\
\qquad \qquad \qquad\qquad \qquad \qquad \qquad  \qquad \qquad y_1\in C(x_1),\, y_2\in C (x_2),\, z\in 
\sigma_{x_1}^{x_2}(y_1)=_{C(x_2)}y_2\,  ]\end{array} $$}

\noindent
and two extensional terms  {\small $b(x,y)\in B_=(x)\ [\, x\in A_=,\, y\in C_=(x)\, ]$} and {\small $b'(x,y)\in B_=(x)\ [\, x\in A_=,\, y\in C_=(x)\, ]$} are equal if and only if 
in \mtt\  there exists a proof of
{\small $$\begin{array}{l}
 b(x_1,y_1)=_{B(x_1)}
\,  b'(x_1,y_1)\ prop\  [\, x_1\in A,\,  y_1\in C(x_1)\,  ]\end{array} $$}
\end{definition}

Now, we are ready to prove that 
{\small \cq}-morphisms correspond to  extensional
dependent collections. Categorically, this can be expressed
by saying that the slice category of {\small \cq\ } over
an extensional collection {\small $(A,=_{A})$}
(see \cite{jacobbook} for the definition of slice category) is equivalent
to the category of extensional dependent collections on  {\small $(A,=_{A})$}:
\begin{prop}
\label{dep}
The category {\small $\mbox{\rm \cq}/(A,=_{A})$}
 is equivalent~\footnote{Equivalence is here formulated as the existence of a dense fullly  faithful functor from {\small $ Dep_{\cq}(\, (A,=_{A})\, )$} to  {\small $\mbox{\rm \cq}/(A,=_{A})$ } unless the axiom of choice is used in the metatheory.} to {\small $ Dep_{\cq}(\, (A,=_{A})\, )$}.
\end{prop}
{\bf Proof.}
We just describe how to associate an extensional dependent collection
to a {\small \cq}-morphism and conversely.

\noindent
Given a {\small \cq}-morphism {\small $f: (C,=_C)\rightarrow (A,=_{A})$} then
the support of the extensional dependent collection associated with it is
{\small $$\Sigma_{y\in C}\,  f(y)=_A x\ col \ [x\in A]$$}

\noindent
with equality 
{\small $$z=_{\Sigma f }z'\, \equiv\, 
\pi_1(z)=_C \pi_1(z')$$}

\noindent
for {\small $x\in A$} and {\small  $z,z'\in \Sigma_{y\in C}\,  f(y)=_A x$}.
Then, we define
{\small $$\sigma_{x_1}^{x_2}(z)\, \equiv\, \langle \pi_1(z), \mathsf{tra}(\pi_2(z),
d)\rangle$$}

\noindent
for {\small $d\in x_1=_A x_2$}.
Clearly $\sigma_{x_1}^{x_1}(z)=_{\Sigma f}z$ and
also the transitivity property holds.

Conversely,
given an extensional dependent type {\small $B_{=}(x)\ [x\in A_{=}]$},
we consider the extensional type {\small  $ (\, \Sigma_{x\in A} B(x)\, ,
\, =_{\Sigma}\, )$}
where for {\small $z,z'\in \Sigma_{x\in A} B(x)$}
{\small $$z=_{\Sigma}z'\, \equiv\, 
\exists{d\in \pi_1(z)=_A \pi_1(z')}\ \ \  \sigma_{\pi_1(z)}^{\pi_1(z')}( \pi_2(z))=_{B(\pi_1(z'))}\pi_2(z') $$}

\noindent
 and we build  the {\small \cq}-morphism
{\small $$\pi^B: (\Sigma_{x\in A} B(x), =_{\Sigma})\rightarrow (A, =_A).$$}

\noindent
called {\it comprehension of the extensional dependent collection
 {\small $B_{=}(x)\ [x\in A_{=}]$}} and given by {\small $\pi^B(z)\, \equiv\, \pi_1(z)$} for {\small $z\in \Sigma_{x\in A} B(x)$}.
\medskip

\begin{definition}\em
We call {\it dset-morphism}
a  {\small \cq}-morphism 
{\small $$\pi^B: (\Sigma_{x\in A} B(x), =_{\Sigma})\rightarrow (A, =_A).$$}

\noindent
that is  the comprehension of an extensional dependent
set {\small $B_{=}(x)\ [x\in A_{=}]$}, namely of an extensional dependent
collection whose support is a dependent set and whose equivalence relation
is small.
\end{definition}

\noindent
In an analogous way, we can define the category of extensional collections
{\small \cqdp}\ on {\small\mttdp}, and we can prove
the same properties shown so far for {\small \cq}.

\noindent
Then, in order to describe the categorical structure of {\small \cq}, as well as
that of {\small \cqdp},
it is also useful to know that these models are closed under finite products
and equalizers (for their definition see, for example, \cite{M71}), 
namely that they are lex 
categories:
\begin{lemma}
\label{lex}
The category {\small \cq}, as well as  {\small \cqdp}, is lex
(i.e. with terminal object, binary products,
equalizers).
\end{lemma}
 {\bf Proof.}
In the following we indicate with {\small $c: (D,=_D)\rightarrow (C,=_C)$}
an arrow given by {\small $c(x)\in C\ [x\in D]$} in {\small \cq\ }.

\noindent
The {\it terminal} object is {\small $(\, \mathsf{N_1}, \mathsf{Id}(\mathsf{N_1}, x,y)\,)$}.
Then, for any object {\small  $(A,=_A)$}, the unique arrow to the terminal object
is {\small  $\star\in \mathsf{N_1}\ [x\in A]$}. The uniqueness can be proved
thanks to the fact that for any  {\small $d\in \mathsf{N_1}$} we can prove
that  {\small $\mathsf{Id}(\mathsf{N_1}, d,\star)$} holds by elimination rule on $\mathsf{N_1}$.

\noindent
The {\it binary product} of  {\small $(A,=_A)$} an {\small  $(B,=_B)$ }is
 {\small $(A\times B, =_{\times})$} where 
{\small 
$$z=_{\times }z'\, \equiv\quad  \pi_1(z)=_{A}\pi_1(z')\ \wedge\
\pi_2(z)=_{B}\pi_2(z')$$}

\noindent
The projections are {\small $\pi_1(z)\in A\ [z\in A\times B]$} and
{\small $\pi_2(z)\in B\ [z\in A\times B]$}.
The pairing of two arrows {\small $a: (C,=_C)\rightarrow (A,=_A)$}
and  {\small $b: (C,=_C)\rightarrow (B,=_B)$} is {\small $\langle a(z),b(z)
\rangle\in A\times B\
[z\in C]$}. 

\noindent
An {\it equalizer} of {\small $b_1,b_2: (A,=_A)\rightarrow (B,=_B)$} is
{\small  $(\, \Sigma_{y\in A}\, b_1(y)=_B b_2 (y)\, ,\, =_{eq}\,)$}
where {\small  $$ z=_{eq}z'\, \equiv\quad  \pi_1(z)=_{A}\pi_1(z')\ $$}

\noindent
for {\small  $z,z'\in \Sigma_{y\in  A}\, b_1(y)=_B b_2 (y)$.}
The embedding morphism is {\small $\pi_1(z)\in A\ [\, z\in \Sigma_{y\in A}\,
 b_1(y)=_B b_2(y) \,]$.}

\noindent
Moreover,
for any arrow {\small $c: (D,=_D)\rightarrow (A,=_A)$} equalizing $b_1,b_2$
and hence yielding to a proof  
{\small $$p(z)\in  b_1(c(z))=_B b_2 (c(z))\ [z\in D]$$}

\noindent
the unique arrow towards the equalizer is {\small  $\langle c(z),p(z) \rangle\in
\Sigma_{y\in A}\, b_1(y)=_B b_2 (y)\ [z\in D]$.}

\noindent
{\it Note that here we have  used the existence of proof-terms witnessing the
equality between morphisms in {\small \cq }.}
\medskip

Now, after knowing how binary products are defined in {\small \cq},
 we give the definition of 
 {\it categorical} equivalence relation induced
by an  equivalence relation of \mtt\ (also small):

\begin{definition} [$\mathcal{P}_q$-equivalence relation]
\label{eqrel}\em
Given an equivalence relation {\small $R\in \mathcal{P}_q( \, (A,=_{A})\times (A,=_{A})\, )$}, namely a proposition 
{\small $R(x,y)\ \mathsf{prop}\ [x\in A, y\in A]$} that preserves
 $=_{A}$ on both dependencies and is also an equivalence relation, 
we call {\it $\mathcal{P}_q$-equivalence relation}  the corresponding ${\mathcal{P}_q}$-proposition, namely  
{\small $$\eta_1\, \equiv\, \pi_1^R:  (\, \Sigma_{z\in A\times A}\,
R(\, \pi_1z\, ,\, \pi_2z\, )\, ,\, =_{Cm}\, )\rightarrow  (A\times A,
=_{\times})$$}
\end{definition}
Analogously, we can define a  {\it ${\mathcal{P}_{q}}$-small equivalence
relation} and a {\it ${\mathcal{P}_{q}}_{set}$-equivalence
relation},
which is indeed also small.

\noindent
In order to describe the categorical structure of {\small \cq},
 as well as that of {\small \cqdp}, it will be useful to know that pullbacks preserve dset-morphisms and
${\mathcal{P}_q}$-propositions, which in turn follows from the fact that
a pullback along the comprehension of an extensional dependent collection
is isomorphic to one of this form.
\begin{lemma}
\label{pul}
In {\small \cq},  as well as in {\small \cqdp}, the pullback of 
the comprehension of an extensional dependent collection {\small 
$B(x)_=\ [x\in A_=]$} 
{\small $$\pi_1^B:\, (\, \Sigma_{x\in A} B(x)\, ,\,  =_{\Sigma}\, )\rightarrow (A,=_A)$$}

\noindent
along a generic morphism {\small $\delta: (D,=_D)\rightarrow (A,=_A)$}
is isomorphic to the comprehension of an extensional dependent
collection {\small $B'(x)_=\ [x\in D_=]$}
where {\small $B'(x)\ [x\in D]\, \equiv\, B(\delta(x))\ [x\in D]$}
 {\small $$\pi_1^{B'}:\, (\Sigma_{x\in D}\ B(\, \delta(x)\, )\, ,\,  =_{B'(x)}\, )\rightarrow (D,=_D).$$}

\noindent
whose
equality is {\small $w=_{B'(x)}w'\, \equiv\, w=_{B(\delta(x))}w'$} for {\small $x\in D,
w,w'\in B'(x)$} and its substitution morphism is defined by using that of
{\small $B(x)$} as follows: for {\small $x_1\in D,  x_2\in D,\, w\in x_1=_{D}x_2, y\in B'(x_1)$}
{\small $$
\sigma^{x_2}_{x_2}(y)\, \equiv\, \sigma_{\delta(x_1)}^{\delta(x_2)}(y)
.$$} 
\end{lemma}

\noindent
From this we deduce the following corollary:
\begin{cor}
\label{pulp}
In {\small \cq},  as well as in {\small \cqdp}, the pullback of a 
dset-morphism
 {\small $$\pi_1^B:\, (\, \Sigma_{x\in A} \, B(x)\, , =_{Cm})\rightarrow
 (A,=_A)$$}

\noindent
along a morphism {\small $\delta: (D,=_D)\rightarrow (A,=_A)$}
is isomorphic to a  dset-morphism.
Analogously, the pullback of a 
${\mathcal{P}_q}$-proposition is isomorphic to a $\mathcal{P}_q$-proposition,
and that of a
$\mathcal{P}_q$-equivalence relation to a $\mathcal{P}_q$-equivalence relation.
The same holds with respect to ${\mathcal{P}_{q}}$-small propositions and 
${\mathcal{P}_{q}}$-small equivalence relations, as well as with respect to
${\mathcal{P}_{q}}_{set}$-propositions and 
${\mathcal{P}_{q}}_{set}$-equivalence relations in  {\small \cqset}.
\end{cor}
{\bf Proof.}
This follows immediately from lemma~\ref{pul} after noticing
that {\small $B'(x)\ [x\in D]$} is of the same kind of type as 
{\small $B(x)\ [x\in A]$},
and it is an equivalence relation  if  {\small $B(x)$} is so. In particular, to apply
lemma~\ref{pul} properly on  $\mathcal{P}_q$-propositions, note that  
any {\small $P\in \mathsf{Ob}\mathcal{P}_q( \, (A,=_{A})\, )$}  induces
an extensional dependent proposition {\small $P(x)_=\ [x\in A_=]$}
whose support is {\small $P(x)\ prop\ [x\in A]$},
whose equality is the trivial one: 
{\small $w=_{P(x)}w'\, \equiv\, \mathsf{tt}$},
and whose substitution morphisms are defined as {\small 
$\sigma_{x_1}^{x_2}(y)\, \equiv\,
\mathsf{Ap}(\, \mathsf{ps}(x_1,x_2, d)\, , y)  $} for {\small $x_1,x_2\in A$} and
{\small $d\in x_1=_A x_2, y\in P(x_1)$}.
\medskip

Now, we are going to recall the categorical notion of {\it dependent product}
and its stability under pullback (see \cite{Seel,handJohn1}).
To this purpose we fix the notation about pullback: in a lex category
 {\small$\cal C$}
we indicate the projections of a pullback of {\small $b:B\rightarrow A$}
along  {\small$d: D\rightarrow A$} as follows
{\small $$
\def\objectstyle{\scriptstyle}
\def\labelstyle{\scriptstyle}
{\xymatrix @-0pc{
 D\times_A B\ar[rr]^{b^{\ast}(d)} \ar[d]_{d^{\ast}(b)}
& & B \ar[d]^{ b}      \\
D\ar[rr]_{d}&& A
}}
$$}
Furthermore, for any  morphism {\small $m: b'\rightarrow b$} in
 {\small${\cal C}/A$}
we indicate with
{\small $d^{\ast}(m): d^{\ast}(b')\rightarrow d^{\ast}(b)$}  the unique morphism 
 in {\small ${\cal C}/D$}
such that {\small $d^{\ast}(b)\cdot d^{\ast}(m)=d^{\ast}(b')$} and
{\small $b^{\ast}(d)\cdot d^{\ast}(m)=m\cdot ( b')^{\ast}(d)$} hold.

Then, we are ready to
 recall the definition of stable dependent product:
\begin{definition}\em
Given a lex category   {\small ${\cal C}$}, we say that {\small $\cal C$ } is closed
under  {\it the dependent product of a morphism {\small $c: C\rightarrow B$} over
a morphism {\small $b: B\rightarrow A$}}, if there exists a {\small ${\cal C}$}-morphism
{\small $\Pi_bc: \Pi_BC\rightarrow A$} with a {\small ${\cal C}/B$}-morphism
{\small $Ap:   b^{\ast} ( \Pi_bc)\rightarrow c$} 
 such that, for every  {\small ${\cal C}$}-morphism 
 {\small $d: D\rightarrow A$}
and any {\small${\cal C}/B$}-morphism
 {\small $m: b^{\ast}( d) \rightarrow c$}, there exists a unique
{\small ${\cal C}/A$}-morphism {\small $\widehat{m}:  d\rightarrow \Pi_bc$} in  ${\cal C}/A$ such
that  {\small $Ap\cdot b^{\ast}(\widehat{m})= m$} in 
 {\small${\cal C}/B$}.

Moreover, a {\it dependent product of a morphism {\small
 $c: C\rightarrow B$ } over
a morphism {\small $b: B\rightarrow A$}} is {\it stable
under pullback} if for every morphism
{\small $q: Q\rightarrow A $}  then 
{\small $q^{\ast}({\Pi_bc}): Q\times_A\Pi_BC\rightarrow Q$}
together with {\small $(b^{\ast}(q))^{\ast}(Ap)$} is a dependent product of
{\small $(b^{\ast}(q))^{\ast}(c):(Q\times_AB)\times_B C \rightarrow Q\times_AB$}
over {\small $q^{\ast}(b): Q\times_AB\rightarrow Q$}.
\end{definition}

Then, we recall the definition of exponential in a slice
category and its stability
under pullback:
\begin{definition}\em
Given a lex category {\small ${\cal C}$}, for any object $A$ in  
{\small ${\cal C}$}
 we say that {\small  ${\cal C}/A$} is closed
under  {\it the exponential of 
{\small  $ c: C\rightarrow A$}
to {\small $b: B\rightarrow A$} },
also called the {\it function space from $b: B\rightarrow A$ to
$ c: C\rightarrow A$}, 
if {\small  ${\cal C}$} is closed under the dependent product
 of {\small $b^{\ast}(c): C\times_AB\rightarrow
B$} over {\small  $b: B\rightarrow A$}.

Moreover, {\it the exponential of {\small $ c: C\rightarrow A$}
 to {\small $b: B\rightarrow A$}} is {\it stable under pullback}, if the corresponding
dependent product is stable under pullback.
\end{definition}

\noindent
Now, we are ready to describe the categorical structure
of  {\small \cq\ } and  {\small \cqdp\ }  
sufficient  to interpret
\emtt\ and \emttdp\ respectively (for the involved categorical definitions
not presented here
 see loc. cit. in \cite{tumscs}):
\begin{theorem}
\label{mainth1}
The following hold:
\begin{list}{-}{}
\item
The category {\small \mbox{\rm \cq}}
 is lex  with  parameterized lists of extensional sets,
 with stable finite propositionally disjoint sums~\footnote{The sums are only propositionally disjoint in the sense that the  pullback vertex of
   the sum injections  is the comprehension of the falsum, which is an
   initial object only in the fibres of ${\cal P}_q$.}  of extensional sets and 
with stable effective quotients
with respect to $\mathcal{P}_{q}$-equivalence relations.

Furthermore,
{\small \mbox{\rm \cq}} is also closed under stable dependent products
of any dset-morphism over another dset-morphism.

Moreover, in {\small \mbox{\rm \cq}}
 there is an object {\small $(\, \mathsf{prop_s}, \leftrightarrow\, )$},
where {\small $\leftrightarrow (p,q)\, \equiv\, ( p \rightarrow q)
\wedge  ( q \rightarrow p)$} for {\small $p,q\in \mathsf{prop_s}$},
 classifying  
${\mathcal{P}_{q}}$-small propositions, namely for every
object {\small $(A,=_{A})$} there is a bijection 
{\small $$\cq(\,  (A,=_{A})\, , \, (\, \mathsf{prop_s}, \leftrightarrow\, )\, )\simeq
 Sub_{{\mathcal{P}_{q}}\mbox{\rm -small}}(\, (A, =_{A})\,)$$ }

\noindent
between morphisms from {\small $(A,=_{A})$}
 to it and the collection of subobjects represented by ${\mathcal{P}_q}$-small
propositions on {\small $(A,=_{A})$}.

And, in  {\small \mbox{\rm \cq}} there exist local stable exponentials of
  $(\, \mathsf{prop_s}, \leftrightarrow\, )$ to dset-morphisms,
namely, 
for any extensional collection  $(A,=_A)$,
the slice category of {\small \mbox{\rm \cq}} over  $(A,=_A)$
is closed  under exponentials
of 
{\small
 $\pi_1:(A,=_A)\times (\, \mathsf{prop_s}, \leftrightarrow\, )\rightarrow (A,=_A)$} to any dset-morphism, and these exponentials are stable under pullback.

Lastly, the indexed category ${\cal P}_q$ validates first-order intuitionistic logic with equality,
namely it is an intuitionistic hyperdoctrine in the sense of \cite{hyper}
(see also \cite{eqhyp,dialectica,Pitts}).

\item
The category {\small \cqdp\ } enjoys the same properties as {\small \cq}, but in addition
is also
locally cartesian closed (i.e. with dependent products).
\end{list}

\end{theorem}

 {\bf Proof.}
Thanks to lemma~\ref{lex}
we already know that  {\small \cq\ } and  {\small \cqdp\ } are lex, and
we proceed to prove all the other properties.

\noindent
The {\it list object} on an extensional set  {\small $(C,=_C)$} is {\small  $(  List(C), =_{List(C)})$}
with
 {\small $$z =_{List(C)}z'\, \equiv \,
\exists_{l\in List(R)} \ \
\mathsf{Id}(\, List(C),\, \overline{\pi_1}(l),\,  z) \ \ \wedge \ \ 
\mathsf{Id}(\, List(C),\, \overline{\pi_2}(l),\,  z')$$}

\noindent
for  {\small $z,z'\in List(C) $}
 where {\small $R\, \equiv\, \Sigma_{x\in C}\, \Sigma_{y\in C} \, x=_{C}y$} and
 {\small $\overline{\pi_i}\equiv List(\pi_i)$} is the lifting
on lists of the $i$-th projection for $i=1,2$.
The empty list arrow is {\small $\epsilon\in List(C)\ [w\in D]$}
and the list constructor is {\small $\mathsf{cons}(z,y)\in  List(C)\ [z\in List(C),y\in
C]$}. Given  the {\small \cq}-morphisms
 {\small $a: (D,=_D)\rightarrow (M,=_M) $} and
 {\small$l: (M,=_M)\times (C,=_C)\rightarrow (M,=_M) $}, the recursor map is 
 {\small${\it El}_{list}(\, u,\,  a(w)\, ,\, (x,y,z).l(\langle z,y\rangle)\,)\in
M\ [w\in D, u\in List(C)] $}.
The uniqueness of the recursor map  follows by elimination rule
on lists.

\noindent
The {\it initial} object is {\small  $(\,\mathsf{N_0}, \mathsf{Id}(\mathsf{N_0}, x,y)\,)$}
and, for any object  {\small  $(A,=_A)$}, the unique arrow from the initial object
to it is  {\small $ \orig{emp_{o}}(x)\in A\ [\, x\in \mathsf{N_0}\, ]$}.
The uniqueness of such an arrow follows by the elimination rule on $\mathsf{N_0}$.

\noindent
The {\it binary coproduct} of extensional sets
 {\small  $(B,=_B)$} and  {\small $(C,=_C)$ } is
 {\small $(B+ C, =_{B+C})$~\footnote{More formally, thanks to disjointness of sums in \mtt\ we can define
$$\begin{array}{ll}
z=_{B+C}z'\, \equiv\,& (\ \exists_{x\in B}\,\exists_{x'\in B}\
\mathsf{Id}(B+C, z, \mathsf{inl}(x))\, \wedge\,\mathsf{Id}(B+C, z',
 \mathsf{inl}(x'))\ \wedge\ x=_Bx'\ )\ \\
&\ \ \vee\ 
 (\ \exists_{y\in C}\,\exists_{y'\in C}\
\mathsf{Id}(B+C, z, \mathsf{inr}(y))\, \wedge\,\mathsf{Id}(B+C, z',
 \mathsf{inr}(y'))\ \wedge\ y=_C y'\ ) \end{array}  $$}} 

\noindent
where 
{\small $$z=_{B+C}z'\, \equiv\,
 \left\{\begin{tabular}{ll}
$b=_{B}b'$ & \qquad\mbox{ if }
 $ z=\orig{inl}(b)\mbox{ and }  z'=\orig{inl}(b')$ \mbox{ for }
$b,b'\in B$\\
 $c=_{C}c'$ &
 \qquad \mbox{ if } $ z=\orig{inr}(c)\mbox{ and }  z'=\orig{inr}(c')$
\mbox{ for }$c,c'\in C$\\
$\bot$ &\qquad \mbox{ otherwise }\end{tabular}\right.$$}

\noindent
for {\small $z,z'\in  B+C$}. 
The injections are {\small $\orig{inl}(z)\in B+C\ [z\in B]$} and
 {\small $  \orig{inr}(z)\in B+C\ [z\in C]$} and the coproduct morphism
of   {\small $b: (B,=_B)\rightarrow (D,=_D)$} and
 {\small $c: (C,=_C)\rightarrow(D,=_D)$}
 is  {\small ${\it El}(z,\,  (y_1). b(y_1),\, (y_2). c(y_2)\, )
\in D\ [z\in B+ C]$.}
Uniqueness of the coproduct morphism follows by elimination rule
on the binary sum $B+C$. 
Sums are disjoint thanks to the  disjointness of sums in \mtt\ (this can be proved easily thanks to the fact that we can eliminate on disjoint sums toward
{\small $\mathsf{prop_s}$} by using 
{\small ${\it El}_+ (z, (x). \bot, (y) .\mathsf{tt}\, )\in \mathsf{prop_s}\
[z\in A+B]$}).
Stability of coproducts under pullbacks
follows by elimination on binary sums and from the fact that in \mtt\
we can prove  injectivity of sum injections~\footnote{For example,
injectivity of $\mathsf{inl}$ can be proved as follows.
Consider the term $p(u,z)\in A\ [u\in A, z\in A+B]$ where
$p(u,z)\, \equiv\, {\it El}_+ (z, (x).x,(y).u)$. Then, if 
$\mathsf{Id}(\, A+B,\, \mathsf{inl}(a),\, \mathsf{inl}(a')\, )$ holds for $a,a'\in A$,
by preservation of propositional equality we get that $\mathsf{Id}(\, A,\,
p( a,\,  \mathsf{inl}(a)\, )\, ,\, p( a,\mathsf{inl}(a'))$ holds, too. Hence,
from this
we conclude that $\mathsf{Id}(\, A,\,
a\, ,\,a')$ holds.}
and sum disjointness.

\noindent
The {\it classifying bijection} 
between morphisms from 
 {\small $(A,=_{A})$}
 to  {\small $(\mathsf{prop_s}, \leftrightarrow)$}
 and subobjects given by
 ${\mathcal{P}_{q}}$-small propositions on  {\small $(A,=_{A})$}
follows by construction of  ${{\cal P}_{q}}$-small propositions.

\noindent
The fact that ${\cal P}_{q}$ validates
  first-order intuitionistic logic with equality follows by construction
of the indexed functor.

\noindent 
Proving  that {\small \cq} is closed
under dependent products of dset-morphisms over
dset-morphisms, as well as proving that
 local stable exponentials of
  $(\, \mathsf{prop_s}, \leftrightarrow\, )$ to dset-morphisms
exist in {\small \cq},
can be seen as a particular case of proving  {\it local closure} of
{\small \cqdp}. Indeed, the latter is equivalent
 to proving that  {\small \cqdp}
is closed under  {\it dependent products} of generic morphisms (usually expressed in the form  that a
right adjoint to the pullback functor induced by any morphism
exists, see \cite{jacobbook,handJohn1}).
We proceed, then, to prove that {\small \cqdp}
is closed under dependent products of any morphism over any other.
 To this purpose,
thanks to the fact, shown in proposition~\ref{dep},
that
 any {\small \cqdp}-morphism {\small $f: (B,=_B)\rightarrow (A,=_A)$ }
 is isomorphic in 
{\small  $\mbox{\rm \cqdp}/(A,=_{A})$}
 to one of the form
{\small $$\pi_1^B:(\Sigma_{z\in A}\, B(z), =_{\Sigma})\rightarrow (A,=_A)$$}

\noindent
with {\small $B_{=}(z)\ [z\in A_{=}]$} extensional dependent collection,
it is enough to show the existence of dependent products for
morphisms of the form $\pi_1^B$.
In the next we will use the abbreviation  {\small $B_\Sigma\, \equiv\,\Sigma_{z\in A}\, B(z)$}.

\noindent
Then, a {\it dependent product} of 
{\small $$\pi_1^C: 
(\Sigma_{z\in B_{\Sigma}}C(z), =_{C_{\Sigma}}\, )\rightarrow (B_{\Sigma}, =_{B_{\Sigma}})$$}

\noindent
over {\small $\pi_1^B:\, (B_{\Sigma}, =_{\Sigma})\rightarrow (A,=_A)$},
is   {\small $\pi_1^{\Pi_B^C}:\,  (\, (\Pi_{B}C)_{\Sigma} \, ,\,  =_{\Pi})\rightarrow (A,=_A)$}
that is the comprehension (in the sense of proposition~\ref{dep}) of
the  extensional dependent collection
{\small $\Pi_{B}\, C(x)_=\ [x\in A_=]$} defined as follows:
 for {\small $x\in A$}
{\small $$\Pi_{B}\, C(x)\, \equiv\, 
 \Sigma_{h\in \Pi_{y\in B(x)}\, C (\langle x,y\rangle)}
\quad\ 
  \forall_{y_1\in B(x)}\ \forall_{y_2\in B(x)}\ \forall_{d\in y_1=_{B(x)}y_2}\ \ \ 
\sigma_{\langle x,y_1\rangle}^{\langle x,y_2\rangle }\, (
\, \mathsf{Ap}(h,y_1)\, )=_{C(\langle x,y_2\rangle )} \mathsf{Ap}(h,y_2) $$}
and 
 {\small $$z=_{\Pi_B^C(x)}z'\, 
\equiv\, 
 \forall_{y\in B(x)}\ \  \mathsf{Ap}(\, \pi_1(z)\, ,\,  y\, )
=_{C(\langle x, y\rangle) } 
\mathsf{Ap}(\, \pi_1(z')\, ,\,  y\, )\ \qquad
\mbox{ for }z,z'\in \Pi_{B} C(x) $$}

\noindent
and for {\small $x_1\in A,
x_2\in A, w\in x_1=_{A}x_2, u\in  \Pi_B^C(x)$} 
{\small $$\sigma_{x_1}^{x_2}(u)\, \equiv\, 
\langle \lambda y^{B(x_2)}. \, \sigma_{\langle x_1, \sigma_{x_2}^{x_1}(y)\rangle}^{
\langle x_2, y\rangle}(\,\mathsf{Ap}( h, \,  \sigma_{x_2}^{x_1}(y)\, )\, )\, , \, t\rangle$$}

\noindent
for {\small $u\, \equiv\, \langle h,p\rangle$} and suitable $t$
built by using $p$ and proof-terms witnessing properties of
  substitution isomorphisms.

\noindent
Then, supposing  that the pullback domain
 {\small $ (B_{\Sigma}, =_{B_{\Sigma}})\times_{(A,=_A)}
 (\, (\Pi_{B}C)_{\Sigma}\, ,\,  =_{\Pi}\, )$}
 is represented by
{\small $$
(\,\,  \Sigma_{w\in B_{\Sigma}}\, \Sigma_{u\in (\Pi_{B}C)_{\Sigma}}\,
 \pi_1(w)=_{A}\pi_1(u)\, , \,\,  =_{\times_{(A,=A)}}\,\,  )$$} 

\noindent
where
{\small $z=_{\times_{(A,=A)}}z'\, \equiv\,  v=_{B_{\Sigma}}v'\, 
\wedge\,  u =_{\Pi}u'
$} assuming that
{\small $z\, \equiv\, \langle v,\,  \langle u, p\rangle \, \rangle$
and  $z'\, \equiv\, \langle v',\, \langle u', p'\rangle\, \rangle$},
then the {\it application map}
{\small $$App: (B_{\Sigma}, =_{B_{\Sigma}})\, \times_{(A,=_A)}\,
(\, (\Pi_{B}C)_{\Sigma}, =_{\Pi}) 
\rightarrow (C_{\Sigma}, =_{C_{\Sigma}})$$}

\noindent
 is given by
{\small$$App(z)\, \equiv\, \langle \  \langle a_2,\ \sigma_{a_1}^{a_2}(b)
\, \rangle\  ,\,    
   \mathsf{Ap}( h,
\sigma_{a_1}^{a_2}(b))\ \rangle\in C_{\Sigma} \qquad \mbox{ for }
  z\in \,\Sigma_{w\in B_{\Sigma}}\,
   \Sigma_{u\in (\Pi_{B}\, C)_{\Sigma}}\, \pi_1(w)=_{A}\pi_1(u)$$}

\noindent
supposing  {\small $\pi_1(z)\, 
\equiv\, \langle a_1,b\rangle $}
and {\small $\pi_1(\pi_2(z))\, \equiv\,
\langle a_2, \langle h,p\rangle\, \rangle $}.
Note that $\sigma_{a_1}^{a_2}$ is well defined
because  {\small $q\, \equiv\, \pi_2(\pi_2(z))\in
a_1=_{A}a_2 $}.

Now, given a morphism  {\small $m: 
(D,=_D)\times_{(A,=_A)}(B_{\Sigma}, =_{B_{\Sigma}}) \rightarrow
(C_\Sigma, =_{C_{\Sigma}})$} in {\small $\cqdp/ (B_{\Sigma},=_{B_{\Sigma}})$ }
where  we assume the support
of the pullback domain
 {\small $ (D,=_D)\times_{(A,=_A)}(B_{\Sigma}, =_{B_{\Sigma}}) 
$} to be {\small $
\Sigma_{z\in D}\ \Sigma_{w\in B_{\Sigma}}\  \delta(z)=_{A}\pi_1(w)$}
for  {\small $\delta: (D,=_D)\rightarrow(A,=_A)$} as above,
then, the abstraction map {\small 
$\widehat{m}(d)\in  (\Pi_{B}\, C)_{\Sigma}\ [d\in D]$}
 is given by
 {\small $$\widehat{m}(d)\, \equiv\, \langle \delta(d), \, \langle\,  \lambda y^{B(\delta(d))}.\,
\sigma_{\pi_1(m(z_d))}^{\langle \delta(d),y\rangle}(\, \pi_2(\,  m(z_d )\, )\, )\ ,p \rangle\,  \rangle$$}

\noindent
where {\small $z_d\, \equiv\, \langle\, d\, ,\, \langle\,  \langle \delta(d), y\rangle\, , \,
\mathsf{ref}(   \delta(d))\, \rangle\ \rangle$} and suitable $p$ built
 by using proof-terms witnessing  properties of substitution isomorphisms
and equality preservation of $m$.

 
\noindent
In particular,
 the exponential of  {\small
 $\pi_1:(A,=_A)\times (\, \mathsf{prop_s}, \leftrightarrow\, )\rightarrow (A,=_A)$} 
to the dset-morphism
{\small $\pi_1^B: (B_\Sigma,=_\Sigma)\rightarrow (A,=_A)$} in {\small  $\cq/(A,=_A)$} is
the comprehension of the extensional dependent collection
{\small 
${\cal P}(B(x))_=\ [x\in A_=]$}
where {\small $${\cal P}( B(x))\, 
\equiv\, \Sigma_{h\in  B(x)\rightarrow \mathsf{prop_s}}\ \ 
  \forall_{y_1\in B(x)}\  \forall_{y_2\in B(x)}\ \  y_1=_{B(x)}y_2 \ \rightarrow\ (\, 
\mathsf{Ap}(h, y_1)\leftrightarrow  \mathsf{Ap}(h,y_2)\, ) $$}

\noindent
and its corresponding equality is
{\small $$z=_{{\cal P}(B(x))}z'\,
\equiv\, \forall_{y\in B(x)}\ \,  \mathsf{Ap} (\, \pi_1(z)\, ,y)\, \leftrightarrow\, 
\mathsf{Ap}(\, \pi_2(z)\, ,y)\ \qquad
\mbox{ for }z,z'\in {\cal P}(B(x)).$$}

\noindent
Stability under pullback of
dependent products of dset-morphisms over dset-morphisms, or
of  local exponentials of {\small $(\mathsf{prop_s}, \leftrightarrow)$},
follows easily
thanks to corollary~\ref{pulp}.

\noindent
The {\it  quotient of a $\mathcal{P}_{q}$-equivalence relation in \cqdp} 
  {\small$$r: (\, \Sigma_{z\in A\times A}\,
R(\, \pi_1 z\, ,\, \pi_2  z\, )\, ,\, =_{Cm}\, )\rightarrow (\, A\times A,=_\times\, )$$}
 
\noindent
is  {\small$(\, A, R\,)$.}
The quotient map from {\small $(A,=_A)$} to {\small $(\, A,\, R\,)$ }
is {\small $z\in A\ [z\in A]$.}
Given a map  {\small$a: (A,=_A)\rightarrow (D,=_D)$} coequalizing
the projections along $r$, namely $\pi_1\cdot r$
and $\pi_2\cdot r$, the unique map from  {\small$(\, A, R\,)$} to {\small $(D,=_D)$}
factoring $a$ is {\small $a(z)\in D\ [z\in A]$ } itself. 
The quotient map satisfies effectiveness by construction.
 The quotient stability under pullback
 follows easily thanks to corollary~\ref{pulp}. 
 \medskip

\noindent
Furthermore, in an analogous way to
theorem~\ref{mainth1},
we can prove that the category {\small \cq$_{\mathsf{set}}$}  enjoys all
 the categorical properties necessary to interpret \emtt-sets
(for their definitions see, for example, \cite{tumscs}).
In particular, to prove local closure we will make use of 
 a proposition analogous to \ref{dep} that can be proved in an
 analogous way as well:
\begin{prop}
\label{depset}
The category of extensional sets
{\small $Dep_{\cqset}(\, (A, =_{A})\, )$}, defined analogously to
definition~\ref{extcat},
on an extensional set \mbox{\small $ (A, =_{A})$} is equivalent to 
{\small $\mbox{\rm \cqset}/(A,=_{A})$}.
\end{prop}

\begin{theorem}
\label{mainth2}
The category {\small \mbox{\rm \cq$_{\mathsf{set}}$}}
 is lextensive (i.e. with terminal object, binary products,
equalizers and stable finite disjoint coproducts), list-arithmetic
(i.e. with parameterized list objects)  and 
locally cartesian closed  with stable effective quotients
with respect to ${{\cal P}_{q}}_{set}$-equivalence relations, and the
embedding {\small$I:  \mbox{\rm \cq$_{\mathsf{set}}$}\longrightarrow \cq$} preserves all such
a structure.
Moreover, the indexed category ${{\cal P}_{q}}_{set}$ validates first-order intuitionistic logic
with equality, namely it is an intuitionistic hyperdoctrine in the sense of \cite{hyper},
 and the natural embedding {\small  $i: {{\cal P}_{q}}_{set}\Longrightarrow
 {\cal P}_{q}\cdot I$} preserves such a structure.
\end{theorem}

\begin{remark}\em
\label{ros}
Note that to prove theorems~\ref{mainth1},~\ref{mainth2} we do not
need to use
any preservation of definitional  equality of discharged premises
in the elimination constructors, which is indeed absent in \mtt. For example, we did not need to use
rule E-eq list) in the appendix~\ref{emtt} about $l$ in {\small ${\it El}_{list}(s,\,  a,\, l\,)$}.
 Indeed,
 uniqueness of  lists
does not need such an equality preservation. The same can be said for binary products and coproducts with respect to elimination rules of strong indexed sums and
disjoint sums of \mtt.

Moreover, the proof of local closure, or existence of suitable dependent products,
 went through without any use of
the $\xi$-rule mentioned in section~\ref{intlev}.
This point was already noticed  
 in categorical terms
in \cite{Ros,bcsr} with the theorem stating 
that the exact completion of a lex category with weak dependent
products has dependent products.
In our case this can be read as follows.
Let us define the category {\small ${\cal C}(\mtt)_{set}$} in this way:
its objects
are \mtt-sets, its  
 morphisms from $A$ to $B$ are terms  {\small $b(x)\in B\ [x\in A]$}, and  two
morphisms  {\small$b_1(x)\in B\ [x\in A]$} and  {\small$b_2(x)\in B\ [x\in A]$} are equal if there exists a proof
of  {\small$\mathsf{Id}(\, B,\,  b_1(x),\,  b_2(x)\, )\ prop\ [x\in A]$} in 
\mtt.
The identity and composition are defined as 
in the syntactic categories in \cite{tumscs}.

Then, to prove the local closure in {\small \cqset\ }
 it is enough that in {\small${\cal C}(\mtt)_{set}$} 
the natural transformation
{\small $$App\cdot ({\pi_1^B})^\ast(-) : {\cal C}(\mtt)_{set}/(A,=_A)(\, \delta \, , \, \pi_1^{\Pi_{B}C}\, )\rightarrow 
{\cal C}(\mtt)_{set}/(B_{\Sigma},=_{B_{\Sigma}})(\,({\pi_1^B})^\ast (\delta)\, ,
 \pi_1^C\, )$$}

\noindent
is surjective  without necessarily enjoying
a retraction.
In turn, to prove such a surjectiveness
it is enough that $\lambda$-terms exist without necessarily satisfying the
$\xi$-rule, which is instead necessary to define a proper retraction.

\end{remark}

\noindent
\subsection{The interpretation of \mbox{\emtt}}
After theorem~\ref{mainth1}, in order to interpret \emtt\ in  {\small \cq\ }, and \emttdp\ in {\small \cqdp},
at a first glance it seems that we could  simply
use
 the interpretation in \cite{tumscs} given
by fibred functors (we remind from \cite{tumscs}
 that this overcomes the problem, first solved in \cite{Hofmann},
of interpreting substitution {\it correctly}
when following the {\em informal} interpretation first given 
by Seely in \cite{Seel} and recalled in \cite{handJohn2}).
But this interpretation requires a choice of   {\small \cq}-structure  in order
to interpret the various constructors. Unfortunately, we are not able
  to fix such a choice, if we take type equality as object equality.  In particular
we are not able to fix a choice of {\small \cq}-equalizers because they depend on the representatives of {\small \cq}-arrows.
Indeed, the  equalizer {\small $Eq(f,g)$} of two \cq-arrows {\small $f,g: (A,=_A)\rightarrow (B,=_B)$}
may be defined as  {\small $\Sigma_{x\in A}\ f(x)=_A g(x) $}, which depends on the chosen
$f,g$ and it is not intensionally equal to that built on  another choice $f',g'$ of the same two {\small \cq}-morphisms,
namely for which
 {\small $f=_{\cq}f'$}
and {\small $g=_{\cq}g'$} hold.
Now, given that we are not able to fix a choice of 
 {\small \cq}-arrows, we are not able to fix a choice of  {\small \cq}-structure.
Luckily, this problem can be avoided
if  we  work in {\small\cq\ } up to isomorphisms.
Indeed, if we take the category {\small \cqis\ } obtained from {\small \cq}\
by quotienting it over isomorphisms, then this category enjoys a unique choice
of the structure needed to give the interpretation in
\cite{tumscs}. In this case, also
the informal way of interpreting a dependent typed calculus  in \cite{Seel}
can be used, because it turns out to be correct. Our solution is to interpret
\emtt\ in a category {\small \cqis}\ obtained by quotienting  {\small \cq}\ only over suitable
isomorphisms, called {\it  canonical isomorphisms}, to make the proof of the validity theorem go through.
In particular, we  interpret \emtt -signatures of types and terms in {\small \cqis}\ 
by just using an interpretation of them  into {\small \cq}, 
or better into extensional dependent types and terms of \mtt.
This interpretation allows us
to determine  {\small \cqis}-objects and morphisms 
where to interpret  \emtt-signatures by selecting their representatives
up to canonical isomorphisms.

\noindent
To this purpose, we need to generalize  definition~\ref{ds} 
of extensional dependent collection and that of extensional dependent set, as well as that
of extensional dependent proposition and small proposition.
To include all cases we speak of extensional dependent type.
\begin{definition}\em
\label{depqt}
An {\it extensional dependent type}
{\small $$B_=(x_1,\dots x_n)\ [\, x_1\in {A_1}_=, \dots x_n\in {A_n}_=\, ]$$}

\noindent
is given by a dependent type
{\small $$B(x_1,\dots x_n)\ type\  [\, x_1\in A_1, \dots x_n\in A_n\, ]$$}

\noindent
together with the following isomorphisms:
{\small $$\begin{array}{l}
\sigma_{x_1,..x_{n}}^{x'_1,..x'_n}(p_1,\dots p_n, z)\in B(x'_1,\dots x'_n)\
[\, x_1\in A_1, \dots x_n\in A_n\, , x'_1\in A_1, \dots x'_n\in A_n\, ,\\
\qquad\qquad\qquad \qquad  p_1\in x_1=_{A_1} x'_1 \, ,\,  p_2\in 
\sigma_{x_1}^{x_1'}(x_2)=_{A_2} x'_2\, \dots \, 
p_n\in \sigma_{x_1,\dots x_{n-1}}^{x'_1,\dots , x'_{n-1}}(x_n)=_{A_n}x'_n ]
\end{array}$$}

\noindent
not depending on the proof-terms $p_i$ for $i=1,\dots , n$  and preserving the equality of the various $A_i$ 
for $i=1,\dots , n$  in the sense of
definition~\ref{ds}.
Such isomorphisms are also closed under identity and composition
as in definition~\ref{ds}. Analogously to definition~\ref{ds} we will use the abbreviation
{\small $$
\sigma_{x_1,..x_{n}}^{x'_1,..x'_n}( z)\, \equiv\,
\sigma_{x_1,..x_{n}}^{x'_1,..x'_n}(p_1,\dots p_n, z)$$}
\end{definition}

\noindent
Then we define the notion of extensional isomorphism between two extensional
dependent types:
\begin{definition}\em 
\label{iso}
Given two extensional terms
{\small $\tau_{B}^{C}(y)\in C_=\ [\Gamma_= , y\in B_=]\mbox{ and }
\tau_{C}^{B}(z)\in B_=\ [ \Gamma_= , y\in C_=]$}, we say
that {\small $\tau_{B}^{C}$} and {\small $\tau_{C}^{B}$}
provide an {\it extensional isomorphism} between the extensional types
{\small $ B_=\ [\Gamma_=]$} and
{\small $ C_=\ [\Gamma_=]$},
if  we can derive  proofs of 
{\small $$\tau_{B}^{C}(\tau_{C}^{B}(x))=_{C(x)} x\ prop\  [\Gamma_= , x\in C_=]
\quad \mbox { and }\quad  \tau_{C}^{B}(\tau_{B}^{C} (x))=_{B(x)} x\ prop\  [\Gamma_= , x\in B_=]$$}
\end{definition}
In the following,  we simply indicate an extensional isomorphism
with one of its parts, given that the inverse is uniquely determined up
to extensional equality.

In the following, when we say that a proposition $P$ implies another proposition $Q$ in \mtt, we mean that in \mtt\ 
there a proof of $P\rightarrow Q$ under the context defining the propositions.

In the interpretation we will need
to use suitable canonical isomorphisms between extensional dependent types defined as follows: 
\begin{definition}[canonical isomorphism]\label{isocan}
\em
Let  {\small $\tau_{B}^{C}(y)\in C_=\ [\Gamma_= , y\in B_=]$} be part of
 an extensional
isomorphism as defined in def.~\ref{iso}.
We say when $\tau_{B}^{C}$ with its inverse $\tau_{C}^{B}$
 is a {\it canonical isomorphism} by induction
on the derivation of  $B_= \ [\Gamma_=]$ and $C_= \ [\Gamma_=]$
as follows.  For easiness, 
here we suppose to work in the larger calculus {\it \mttdp}:
\begin{list}{-}{}
\item
If {\small $B= C\, \equiv\, \mathsf{prop_s}$}, or  
{\small $B\,= C\, \equiv\, \mathsf{N_0}$}, or {\small  $B=C\, \equiv\, \mathsf{N_1}$}, 
  then {\small $\tau_{B}^{C}(y)\,=_{B}\, y$} holds
under the appropriate context in \mtt.
\item
If {\small $B\, \equiv\, \Sigma_{x\in A}\, D(x)$} and
{\small $C\, \equiv\, \Sigma_{x\in A'}\, D'(x)$} where in \mtt\  both  $=_{B}$ and $=_{C}$ imply $=_{\Sigma}$  instantiated respectively for $B$ and $C$
as in proposition~\ref{dep},   with {\small $D_=\ [\Gamma_=, x\in A_=]$ } and {\small $D'_=\ [\Gamma_=, x\in A'_=]$}
and {\small $A_=\ [\Gamma_=]$} and {\small $A'_=\ [\Gamma_=]$} extensional
dependent types having canonical substitution isomorphisms, then
{\small $\tau_{B}^{C}(y)\, =_C\, {\it El}_{\Sigma}(\, y,\,
(w_1,w_2). \langle \tau_{A}^{A'}(w_1), \tau_{D}^{D'}(w_2)
\rangle\ )$} holds under the appropriate context  with  {\small $\tau_{A}^{A'}$} and {\small $\tau_{D}^{D'}$}
canonical isomorphisms where in particular {\small $\tau_{D}^{D'}(y)\in D'(\tau_A^{A'}(x))\ [\Gamma_=, x\in A_=,y\in D_=(x)]$} (note that we still use the notation {\small $\tau_{D}^{D'}$}
even if the types are on different contexts).
Moreover, also the substitution isomorphisms of  $B_=$ and $C_=$ and
  {\small $\tau_{C}^{B}$} are of this form.
\item
If  {\small $$B\, \equiv\,\Sigma_{h\in \Pi_{x\in A}\, D(x)}\quad\ 
  \forall_{x_1\in  A}\ \forall_{x_2\in  A}\ \forall_{d\in  x_1=_{A}x_2}\ \ \
\sigma_{\overline{x},x_1}^{\overline{x},x_2}\, (\mathsf{Ap}(h,x_1))=_{D(x_2)}
\mathsf{Ap}( h,x_2) $$}

\noindent
 with $\overline{x}$ the variables
in $\Gamma$ and analogously
{\small $$C\, 
\equiv\, \Sigma_{h\in \Pi_{x\in A'}\, D'(x)}\quad\ 
  \forall_{x_1\in A'}\ \forall_{x_2\in A'}\ 
\forall_{d\in  x_1=_{A'}x_2}\ \ \
\sigma_{\overline{x},x_1}^{\overline{x},x_2}\, (\mathsf{Ap}(\, h,\, x_1\, ))=_{D'(x_2)}
\mathsf{Ap}(\,  h,\, x_2\, )        $$ }

\noindent
where  $=_{B}$ implies in \mtt\  the following equivalence relation for {\small $ z,z'\in B$}
{\small $$\forall_{x\in A}\ \  \mathsf{Ap}(\, \pi_1(z)\, ,\, x\, )\, =_{D(x)}\, 
\mathsf{Ap}( \, \pi_1(z')\, ,\, x\, )\ $$}

\noindent
and $=_C$ is of analogous form,
with {\small $D_=\ [\Gamma_=, x\in A_=]$} and  {\small $D'_=\ [\Gamma_=, x\in A'_=]$}
and {\small $A_=\ [\Gamma_=]$} and {\small $A'_=\ [\Gamma_=]$ }
extensional dependent types having canonical substitution isomorphisms, then we can derive a proof of
{\small$$\tau_{B}^{C}(y)\, =_C\, \langle\,  \lambda w^{A'} .\,
\sigma_{\tau_{A}^{A'}(\tau_{A'}^{A}(w))}^{w} (\, 
\tau_D^{D'}(\,\mathsf{Ap} (\, \pi_1(y)\, ,\, \tau_{A'}^A(w))\, )\, )\, , \, p\rangle$$}

\noindent
for some proof-term $p$   under the appropriate context 
 with  {\small$\tau_{A}^{A'}$} and {\small$\tau_{D}^{D'}$}
canonical isomorphisms  where in particular {\small$\tau_{D}^{D'}(y)\in D'(\tau_A^{A'}(x))\ [\Gamma_=, x\in A_=,y\in D_=(x)]$.}
Moreover, also the substitution isomorphisms of  $B_=$ and $C_=$ and
  {\small $\tau_{C}^{B}$} are of this form.
\item
If {\small $B\, \equiv\, List(D)$ } and
{\small$C\, \equiv\, List(D')$ }  where in \mtt\ both $=_{B}$ implies $=_{List(D)}$  and $=_C$ implies $=_{List(D)}$  defined as    in the proof of theorem~\ref{mainth1}
with {\small$D_=\ [\Gamma_=]$ } and {\small $D'_=\ [\Gamma_=]$ } 
extensional dependent types having canonical substitution isomorphisms,
then
{\small$\tau_{B}^{C}(y)\, =_C\, 
{\it El}_{List}(y,\, \epsilon ,\,  (y_1,y_2,z).\, \mathsf{cons}( \,
z\, ,\, 
\tau_D^{D'}(y_2)\, )\,  )$} holds under the appropriate context
 with   $\tau_{D}^{D'}$
canonical isomorphism. Moreover, also the substitution isomorphisms of  $B_=$ and $C_=$ and  {\small $\tau_{C}^{B}$} are of this form.
\item
If {\small $B\, \equiv\, A+D$} and
{\small$C\, \equiv\, A'+D'$} where in \mtt\   both $=_{B}$  implies  $=_{A+D}$ and  $=_C$ implies $=_{A'+D'}$  defined as in the proof of theorem~\ref{mainth1}
with {\small $A_=\ [\Gamma_=]$} and {\small$A'_=\ [\Gamma_=]$} and {\small $ D_=\ [\Gamma_=]$} and {\small  $D'_=\ [\Gamma_=]$}
 extensional
dependent types having canonical substitution isomorphisms,
then
{\small$\tau_{B}^{C}(y)\, =_C\, 
{\it El}_{+}(y, (y_1).\, \mathsf{inl} (\, \tau_{A}^{A'}(y_1)\, )
\, ,\, (y_2).\, \mathsf{inr}(\, \tau_{D}^{D'} (y_2)\,)\, )$}  holds under the appropriate context
 with {\small $\tau_{A}^{A'}$} and {\small$\tau_{D}^{D'}$}
canonical isomorphisms. Moreover, also the substitution isomorphisms of  
$B_=$ and $C_=$ and
  {\small $\tau_{C}^{B}$} are of this form.
\item
If {\small $B\ prop\ [\Gamma]$} and {\small $C\ prop\ [\Gamma]$} are
derivable and  $=_B$ and $=_C$ are the trivial relation equating all proofs, then any {\small $\tau_B^C$} with {\small $\tau_C^B$} is canonical.

\end{list}

\end{definition}

\begin{remark}\em
Note that the above definition of canonical isomorphisms
 can be adapted to work for {\small \mtt}
by simply adding a case analogous to the third one where
{\small $$D(x)\, \equiv\, \mathsf{prop_s}\qquad  z=_{D(x)}z'\, \equiv\, z\leftrightarrow z'\ \mbox{ for } z,z'\in\mathsf{prop_s} $$}
with {\small  $\sigma_{\overline{x},x_1}^{\overline{x},x_2}(w)\, \equiv\, w\ \mbox{ for }
w\in D(x_1)$} and in particular {\small 
$$B\, \equiv\,\Sigma_{h\in A\rightarrow  \mathsf{prop_s} }\quad\ 
  \forall_{x_1\in A}\, \forall_{x_2\in A}\ \ x_1=_{A}x_2\, \rightarrow \, (\, \mathsf{Ap}(h,x_1)
\leftrightarrow \mathsf{Ap}(h,x_2)\, ) $$} 

\noindent
and furthermore {\small $C$} is defined analogously.
\end{remark}

Note that canonical isomorphisms are closed under composition, and more importantly
 between two extensional dependent types there is at most only one canonical isomorphism extensionally:
\begin{prop}
\label{coer}
Canonical isomorphisms enjoy the following properties:
\begin{list}{-}{ }
\item
If {\small$\tau_{A}^{B}(y)\in B_=\ [\Gamma_= , y\in A_=]$}
and   {\small$\tau_{B}^{C}(y)\in C_=\ [\Gamma_= , y\in B_=]$} are canonical isomorphisms,
then  {\small$\tau_{B}^{C}(\tau_{A}^{B}(y))\in  C_=\ [\Gamma_= , y\in A_=]$} is also a canonical isomorphism.

\item If  {\small$\tau_{A}^{B}(y)\in B_=\ [\Gamma_= , y\in A_=]$} and
 {\small ${\tau'}_{A}^{B}(y)\in B_=\ [\Gamma_= , y\in A_=]$} are both canonical isomorphisms, then we can prove
that they are equal as extensional terms, namely that we can derive a proof of
 {\small$$ \tau_{A}^{B}(y)=_{B}{\tau'}_{A}^{B}(y)\ prop\ [\Gamma_= , y\in A_=]$$ }
 
 \item If  {\small$\tau_{A}^{B}(y)\in B_=\ [\Gamma_= , y\in A_=]$}  is a canonical isomorphism, then 
 {\small$\tau_{A}^{B}(y)\in B_='\ [\Gamma_= , y\in A_=']$}   is also a canonical isomorphism provided that it is  an extensional isomorphism with $='_{A}$ implying $=_A$
 and $='_{B}$ implying $=_B$.
\end{list}
\end{prop}
{\bf Proof.}
The proof is by induction on the definition of canonical isomorphism by showing 
all  the points of   the main statement in the same time together with  the facts
that for any canonical isomorphism  {\small$\tau_{A}^{B}(y)\in B_=\ [\Gamma_= , y\in A_=]$}   then all all the substitutions morphisms of $A_=$ and $B_=$ are canonical
and that  for any extensional dependent type  $B_=\ [\Gamma_= ]$ with canonical substitution morphisms then the identity $\tau_B^B(w)\ \equiv\ w \in B_=\ [\Gamma_= ]$  is also
a canonical isomorphism, and finally canonical isomorphisms have canonical inverses.

\medskip

Then we define the quotient category of {\small \cq}\ over  canonical isomorphisms in which we will interpret \emtt. The equivalence relation generated by canonical isomorphisms coincides with extending the previous  collection of canonical isomorphisms by including all identity morphisms between objects of {\small \cq}~\footnote{We still get a well defined quotient category, if we make canonical also  all the morphisms between objects  $ {\small \cq}$  which include the identity morphism between their supports as a representative (and hence the equivalence relations  of the two objects are equivalent).}. \begin{definition}\em
We call {\small \cqis} the category obtained by quotienting  {\small \cq}\ 
over  canonical isomorphisms:
namely an object of {\small \cqis}\ is given by the equivalence class of a \cq-object {\small $(A,=_A)$} 
{\small $$[\, (A,=_A)\, ]$$}

\noindent
 where two objects are defined to be equivalent if 
they are isomorphic via a canonical isomorphism,
\noindent
and a morphism from {\small $[\, (A,=_A)\, ]$} to{ \small  $[\, (B,=_B)\, ]$}
is the equivalence class given by a \cq-morphism {\small $f:  (C,=_C)\rightarrow (D,=_D)$}
{\small $$[f]:  [(A,=_A)]\rightarrow [(B,=_B)]$$}

\noindent
such that {\small $(C,=_C)$ } is canonically isomorphic to  {\small $(A,=_A)$}, i.e. via a canonical isomorphism, and
{\small $(D,=_D)$} canonically isomorphic to {\small$(B,=_B)$},
where two such  morphisms {\small $f,g$}   are
defined to be equivalent,
if, supposing {\small  $f:  (C,=_C)\rightarrow (D,=_D)$} and {\small$g: (M,=_M)\rightarrow (N,=_N)$}
then
{\small $$g\cdot \tau_C^M=_{\cq}\tau_D^N\cdot f$$}

\noindent
for canonical isomorphisms {\small$ \tau_C^M$} and {\small$\tau_D^N$}, given
that {\small $(C,=_C)$ } turns out to be  canonically isomorphic to 
 {\small $(M,=_M)$} and
{\small $(D,=_D)$}  to {\small$(N,=_N)$}.
The equality between two morphisms coincides with the equality of their
 equivalence
classes.
The unit and composition are those inherited from  {\small \cq}.
\end{definition}
The category  is well defined
thanks to the first two  properties of canonical isomorphisms in prop.~\ref{coer}, which continue to hold with the addition of all identity morphisms.

We define the category {\small \cqdpis} analogously and we will use it to
interpret \emttdp.
\\

\noindent
{\bf Why  canonical isomorphisms.}
We will use canonical isomorphisms to interpret equality between
\emtt (\emttdp)-types.
The reason  is the following.
The underlying assumption is that,  in order to be able to interpret
quotient types, we 
 interpret \emtt-dependent types as \mtt-extensional dependent types
and \emtt-terms as \mtt-extensional terms.
 Then, it follows
that the definitional
equality between \emtt-terms must be interpreted in the existence
of a proof that the two terms
are equal according to the equality associated with their type interpretation.
 Now,
suppose that the
 \mtt-extensional dependent collection {\small $ B^I_=(x) \ [x\in A^I_=]$}  interprets
the \emtt\ collection  {\small $B(x)\ col\  [x\in A]$}, 
and that the judgement {\small$a_1=a_2\in A$} holds in \emtt\
and  is valid in \cq,
namely that  {\small$a_1^I=_{A^I} a_2^I$ } holds in \mtt. 
Now, in \emtt\ we get also
that  {\small$B(a_1)=B(a_2)$} holds, too, but in \mtt\ we just know that
  {\small$$B^I(a_1)_= \mbox{ is isomorphic to } B^I(a_2)_=$$}

\noindent
{\it via a substitution isomorphism}.
Therefore, we are forced to interpret type equality as isomorphism of types,
and to this purpose we introduce in \mtt\ the
judgement {\small $A_= =_{ext}B_=\ [\Gamma_=]$} for saying that 
{\small $A_=\ [\Gamma_=]$}
 is isomorphic to  {\small $B_=\ [\Gamma_=]$} via an \mtt-extensional
isomorphism.
But then, in order to make  the rule conv)
in appendix~\ref{mttsyn} valid as well, we may need to correct the interpretation $b^I$ of
an \emtt-term $b$ via an isomorphism.
Indeed, supposing that {\small $a^I\in {A_1^I}_=\ [\Gamma^I_=]$} in \mtt\  interprets
the derived \emtt-judgement {\small  $a\in A_1\ [\Gamma]$ }
 and that the \mtt-judgement
 {\small  ${A_1^I}_= =_{ext} {A_2^I}_=\ [\Gamma^I_=]$}
interprets the derived \emtt-judgement   {\small $A_1=A_2\ col \ [\Gamma]$},  in order
to make the interpretation of the \emtt-judgement {\small  $a\in A_2\ [\Gamma]$}
derived by conv) valid, given that in \mtt\ we only know  that   {\small 
${A_1^I}_= =_{ext} {A_2^I}_=\ [\Gamma^I_=]$} holds,
  we need to introduce in \mtt\  the judgement
 {\small $$a^I\in_{ext}  {A_2^I}_=\ [\Gamma^I_=]$$}

\noindent
to express that $a^I$ is {\it extensionally } of type $A_2^I$
if it belongs to a type isomorphic to it, namely if in \mtt\ we can derive
 {\small $$\tau_{A_1^I}^{A_2^I}(a^I)\in  {A_2^I}_=\ [\Gamma^I_=]$$}

\noindent
{\it via a canonical isomorphism $\tau_{A_1^I}^{A_2^I}$}.
But now, given that the interpretation of \emtt-terms in \mtt\
depends on isomorphisms, then  the interpretation
of \emtt-types depends on  isomorphisms, too, because types may depend on terms.
Luckily, we are able to give the interpretation by making use of {\it  canonical
isomorphisms only}, and hence we require the isomorphisms used so far to be
canonical.
The fact  that canonical isomorphisms
 between two interpreted types are at most one extensionally will allows us
  to prove the validity theorem.

Before interpreting  \emtt-signatures,
 we say when two extensional dependent types under a common context
and on different contexts  are isomorphic:

\begin{definition}\em 
\label{isoext} In \mtt\
we say that the extensional dependent type
{\small $ B_=\ [\Gamma_=]$} is {\it extensionally equal} to the extensional dependent type 
{\small $ C_=\ [\Gamma_=]$}
with the judgement
{\small $$ B_==_{ext} C_=\ [\Gamma_=]$$}

\noindent 
 if they are
 isomorphic via a canonical extensional
isomorphism as in  definition~\ref{isocan}.

\noindent
We will generally call the isomorphism components
{\small $\tau_{B}^{C}(y)\in C_=\ [\Gamma_= , y\in B_=]\mbox{ and }
\tau_{C}^{B}(z)\in B_=\ [ \Gamma_= , z\in C_=]$}.
\end{definition}

\begin{definition}\em 
\label{contis}
Given the \mtt-extensional dependent types
{\small $B_= \ [\Gamma_{=}]$} and {\small $D_=\ [\Delta_{=}]$}
where {\small $\Gamma_=\, \equiv\,
 x_1\in {A_1}_=, \dots ,x_{n}\in {A_{n}}_=$ } and
{\small $\Delta_=\, \equiv\, y_1\in {C_1}_=, \dots ,y_{n}\in {C_{n}}_=$}
we introduce the judgement
{\small $$B_= \ [\Gamma_{=}]\, =_{ext}\, D_=\ [\Delta_{=}]$$}

\noindent
to express that
 in \mtt\ we can derive 
{\small $A_1=_{ext}C_1$} and
{\small $ {A_i}_= \, =_{ext}\,
 D_i(\, \tau_{A_1}^{C_1}(x_1),\dots , \tau_{A_{i-1}}^{C_{i-1}}(x_{i-1})\, )_=\
[\, x_1\in {A_1}_=, \dots ,x_{i-1}\in {A_{i-1}}_=\, ]$}
 for canonical isomorphisms {\small $\tau_{A_i}^{C_i}$} for {\small
$i=2,\dots, n$} if 
{\small $n\geq 2$},
 and also
{\small 
$$B_=\, =_{ext}\, 
\widetilde{D}_=\ [\Gamma_=]$$}

\noindent
via a canonical isomorphism {\small $\tau_B^{D}$} 
where {\small  $\widetilde{D}\, \equiv\,D\,[y_1/ \tau_{A_1}^{C_1}(x_1),\dots ,
y_n/ \tau_{A_{n}}^{C_{n}}(x_{n})\,]  $} whose equality is
 {\small $w=_{\widetilde{D}}w' \, \equiv\,
z=_Dz'\, [y_1/ \tau_{A_1}^{C_1}(x_1),\dots ,
y_n/ \tau_{A_{n}}^{C_{n}}(x_{n})\, , z/w,\, z'/w]$}. 
\end{definition}

\noindent
Then, we are ready to define  when a term belongs to a type in an extensional way:
\begin{definition}\em
\label{termis}
Given an  \mtt-extensional term
{\small  $b\in D_=\ [\Delta_{=}]$} and an \mtt-extensional dependent type  
{\small  $B_=\ [\Gamma_=]$} where {\small $\Gamma_=\, \equiv\,
 x_1\in {A_1}_=, \dots ,x_{n}\in {A_{n}}_=$ } and
{\small $\Delta_=\, \equiv\, y_1\in {C_1}_=, \dots ,y_{n}\in {C_{n}}_=$},
and supposing that  {\small $B_=\ [\Gamma_{=}]\, =_{ext}\, D_=\ [\Delta_{=}]$}
with  canonical isomorphisms {\small  $\tau_{A_i}^{C_i}$} for {\small
$i=1,\dots, n$}, if 
{\small $n\geq 1$}, and {\small $\tau_B^D$},
then  in \mtt\ we say  that {\it $b$ is extensionally of type {\small $B_=\ [\Gamma_{=}]$}} with the judgement
{\small $$b \in_{ext}  B_{=}\ [\Gamma_{=}]$$}

\noindent
if in \mtt\ we can derive 
{\small $$  \widetilde{b}\in
\widetilde{B}_=\ [\Gamma_=]$$}

\noindent
where {\small $\widetilde{b}\, \equiv\, \tau_D^B(\, b (\, 
\tau_{A_1}^{C_1} (x_1),\dots ,\tau_{A_n}^{C_n}(x_n) \, )\, )$}.
\end{definition}

\noindent
Then, we define  extensional equality between terms:
\begin{definition}\em
\label{ist}
Given two \mtt-extensional terms
{\small $$b\in_{ext} B_=\ [\Gamma_{=}]\qquad c\in_{ext} B_=\ [\Gamma_{=}]$$}
in \mtt\ we say that they are {\it extensionally equal terms} with the judgement 
{\small $$b \, =_{ext}\, c \in_{ext}  B_=\ [\Gamma_{=}]$$}
if and only if in \mtt\ we can derive a proof
{\small $$p\in \widetilde{b}=_{\widetilde{B}}
\widetilde{c}\quad [\, \Gamma \,]$$}

\noindent
where {\small $\widetilde{b}, \widetilde{c}$} and {\small $\widetilde{B}$} are defined
respectively as in definitions~\ref{termis} and \ref{contis}.
\end{definition}

{\bf Interpretation of \emtt\ types and terms.}
Now we are ready to define the interpretation of \emtt\ type and term signatures  as 
 \mtt-extensional dependent types, with canonical substitution morphisms, and terms, respectively.
Then,  a type judgement
will be interpreted  as an extensional dependent type
{\small $$(\, B\ type\  [\,\Gamma\, ]\, )^I
\ \, \equiv\, B^I_=\ [\,\Gamma^I_{=}\, ]$$}

\noindent
and a term judgement as
 an extensional term which belongs {\it only extensionally} to the interpretation
of the assigned type as follows:
{\small $$ (\, b\in B\ [\Gamma])^I
\, \equiv\, b^I\in_{ext} B^I_= \ 
[\,\Gamma^I_{=}\, ]$$}

\noindent
In order to give an idea on how the interpretation is defined,
suppose to interpret {\small $\orig{Eq}(B, b_1, b_2)\ prop\
[\Gamma]$} assuming that {\small $(\, B\ type\  [\,\Gamma\, ]\, )^I
\ \, \equiv\, B^I_=\ [\,\Gamma^I_{=}\, ]$}. Then, the term signatures
$b_1$ and $b_2$ under context $\Gamma$  are  assumed to be interpreted as \mtt-extensional terms
{\small $b_1^I\in C_=\ [\Gamma^I_=]$} and {\small $b_2^I\in M_=\ [\Gamma^I_=]$}.
Now, to give the interpretation we do not require as usual that
  {\small $C_=$ } is equal in \mtt\ to  {\small $M_=$} and to {\small $B^I_=$},
but only that it is isomorphic to them via canonical isomorphisms. Hence we put
{\small $(\, \orig{Eq}(B, b_1, b_2)\ prop\ [\Gamma]\, )^I\, \equiv\,  \widetilde{b_1^I}=_{B^I} \widetilde{b_2^I}\ prop\ [\Gamma^I]$}
where we have corrected the interpretation of  $b_1^I$ and $b_2^I$ to match the type {\small $B^I$} as in definition~\ref{termis} by means of canonical isomorphisms.

This explains why we  give an interpretation
$(-)^I: \emtt\rightarrow \mtt $
 of \emtt-type and term signatures as
extensional dependent types with canonical substitution morphisms and extensional terms in \mtt,
that is not only partial as usual interpretations
of dependent type theories (see first paragraph in appendix~\ref{emcq}),
but also uses canonical isomorphisms.

Analogously, we define an  interpretation
$ (-)^I: \emttdp \rightarrow \mttdp $
of \emttdp-type and term signatures as \mttdp\  
extensional dependent types with canonical substitution morphisms and  extensional terms.

{\it The interpretations are properly defined in appendix~\ref{emcq}}.
Here, we just show the interpretation of the power collection of
the singleton  {\small ${\cal P}(1)$} and 
of dependent product sets, quotient sets and function collections
 towards  {\small${\cal P}(1)$}
 with their terms.
Therefore, to interpret these \emtt-types as extensional dependent types of \mtt,
we need to specify the support  of their interpretation
with  related
equality and substitution morphisms. Note that in the case we are  interpreting a type
or a term that requires to have already
interpreted more than one term,   we need to match the types of such terms and we assume
to correct them via {\it canonical isomorphisms}.
In the following,
we simply write
 $b^{\widetilde{I}}$ instead of $\widetilde{b^I}$.

\noindent
The power collection of the singleton is interpreted as the extensional collection classifying  ${\mathcal{P}_{q}}$-small propositions
in theorem~\ref{mainth1}, namely as the \mtt-collection
of small
propositions equipped with equiprovability as equality:

\noindent
{\small
$\begin{array}{l}
\mbox{\bf Power collection of the singleton}:\\ 
{\cal P}(1)^I\ col\  [\Gamma^I]\, \equiv\, 
\mathsf{prop_s}\ [\Gamma^I]\\[3pt]
\mbox{and }z=_{\mathsf{prop_s}^I}z'\, \equiv\,
 ( z \rightarrow z')
\wedge  ( z' \rightarrow z)$ for $z,z'\in \mathsf{prop_s}\end{array}$
\\
$\sigma_{\overline{x}}^{\overline{x}'}(w)\, \equiv\, w$ for
$\overline{x}, \overline{x'}\in \Gamma^I, w\in \mathsf{prop_s}$.\\
$([A])^I\, \equiv\, A^I$} for $A$ small proposition.
\\

\noindent
The dependent product set is interpreted similarly to the
extensional collection behind the dependent product construction 
in theorem~\ref{mainth1}, namely
as the strong indexed sum of functions  preserving the corresponding 
equalities. Two elements of this indexed sum are considered equal
if their first components, which are  lambda-functions, send a given
element to equal elements. This interpretation validates both
 $\beta$ and $\eta$ equalities
for functions and also the $\xi$-rule.

\noindent
{\small
$\begin{array}{l}
\mbox{\bf Dependent Product set}:\\
(\, \Pi_{y\in B} C(y)\, )^I\ set \ [\Gamma^I]\, \equiv\,\ 
 \Sigma_{h\in \Pi_{y\in B^I} \ C^I(y)}\quad\ 
  \forall_{y_1\in B^I}\ \ \forall_{y_2\in B^I}\ \ \forall_{d\in
  y_1=_{B^I}y_2}\ \
\sigma_{\overline{x},y_1}^{\overline{x},y_2}\, (\,\mathsf{Ap}(h,y_1)\,)
=_{C^I (y_2)}\mathsf{Ap}( h,y_2) \\
\mbox{and }z=_{\Pi}z'\, 
\equiv\, \forall_{y\in B^{I}}\ \  \mathsf{Ap}(\, \pi_1(z)\, ,y)=_{C^{I}(y)}\mathsf{Ap}(\,
\pi_1(z')\, ,y)\ 
\mbox{ for }z,z'\in (\,  \Pi_{y\in B}\, C(y)\, )^{I}\end{array}$
\\
$(\,\, \lambda y^{B}.c\, )^I\, \equiv\, \langle\,
\lambda y^{B^I}. c^{\widetilde{I}}\, , \, p \, \rangle$ where $p\in
\forall_{y_1\in B^I}\ \ \forall_{y_2\in B^I}\ \,  \forall_{d\in y_1=_{B^I}y_2}\
\sigma_{\overline{x},y_1}^{\overline{x},y_2}(\, c^{\widetilde{I}}(y_1)\, )
=_{C^I(y_2)} c^{\widetilde{I}}(y_2)$\\
$(\, \mathsf{Ap}(f,b) \, )^I\, \equiv\, \mathsf{Ap}(\,
\pi_1(f^{\widetilde{I}})\, , b^{\widetilde{I}})$\\
$\sigma_{\overline{x}}^{\overline{x'}}(w)\, \equiv\,\langle\, 
\lambda y'^{B^{I}(\overline{x'})}.\,  \sigma_{\overline{x},
\sigma_{\overline{x'}}^{\overline{x}}(y')}^{ \overline{x'},y'}(\,
\mathsf{Ap}(\, \pi_1(w)\, ,\, \sigma_{\overline{x'}}^{\overline{x}}(y') \,  )\, ), p\rangle$ for
$\overline{x}, \overline{x'}\in \Gamma^I$ and
 $w\in (\, \Pi_{y\in B}\, C(y)\, )^{I}(\overline{x})$.\\}
where $p$ is the proof-term witnessing the preservation of equalities obtained from $\pi_2(w)$.
\\

\noindent
The quotient on a set whose interpretation has support $A^I$ is interpreted
as the extensional set with same support $A^I$, but whose equality
is the interpretation
of the quotient equivalence relation, in a similar way to the construction of quotients in theorem~\ref{mainth1}. Then the
interpretation of the quotient map is simply given by the identity and effectiveness becomes trivially validated:

{\small
\noindent
$\begin{array}{l}
\mbox{\bf Quotient set}:\\
 (\, A/R\ set\ [\Gamma] )^I\, \equiv\, A^I \ set\ [\Gamma^I] \\[3pt]
\mbox{and }z=_{A/R^I}z' \, \equiv\,  R^I(z,z')\mbox{ for } 
z,z'\in A^I\end{array}$
\\
$(\, [a]\, )\, \equiv\, a^I$ and
 ${\it El}_{Q}(p,l)^I\, \equiv\, l^{\widetilde{I}}(p^{\widetilde{B}}) $\\
$\sigma_{\overline{x}}^{\overline{x'}}(w)$} is defined
as the substitution isomorphism of {\small $A^I_=\, [\, \Gamma^I_=]$}.
\\

\noindent
A function collection towards {\small ${\cal P}(1)$}  is interpreted as 
the 
extensional collection behind the local exponential
 construction
of the ${\mathcal{P}_{q}}$-small proposition classifier
in theorem~\ref{mainth1}:\\

{\small
\noindent
$\begin{array}{l}
\mbox{\bf Function collection toward ${\cal P}(1)$}:\\
 (\, B\rightarrow {\cal P}(1)\  col\ [\Gamma]\, )^I\, \equiv\,
\Sigma_{h\in B^I\rightarrow \mathsf{prop_s}} \ 
\quad\ 
  \forall_{y_1\in B^I}\ \ \forall_{ y_2\in B^I}\ \
  y_1=_{B^I}y_2\, \rightarrow\, (\, \mathsf{Ap}(
h,y_1)\, \leftrightarrow\, \mathsf{Ap}(h, y_2)\, ) \\[2pt]
\mbox{and }z=_{{\cal P}}z'\, 
\equiv\, \forall_{y\in B^{I}}\ \ \mathsf{Ap}( \, \pi_1(z)\, ,y)\, \leftrightarrow 
 \, \ \mathsf{Ap}(\, \pi_1(z')\, ,y)\ 
\mbox{ for }z,z'\in (\, B\rightarrow {\cal P}(1)\, )^{I}\end{array}$\\
$(\, \lambda y^{B}.c\, )^I\, \equiv\, \langle\,
\lambda y^{\widetilde{B}}. c^{\widetilde{I}}\, , \, p \, \rangle$ where $p\in
\forall_{y_1\in B^I}\ \ \forall_{ y_2\in B^I}\ \ y_1=_{B^I}y_2\, \rightarrow\,
(\, c^{\widetilde{I}}(y_1)\, \leftrightarrow \, c^{\widetilde{I}}(y_2)\, )$\\
$(\, \mathsf{Ap}(f,b) \, )^I\, \equiv\, \mathsf{Ap}(\, 
\pi_1(f^{\widetilde{I}})\, , b^{\widetilde{I}})$\\[2pt]
$\sigma_{\overline{x}}^{\overline{x'}}(w)\, \equiv\,\langle\, 
\lambda y'^{B^{I}(\overline{x'})}.\,  \sigma_{\overline{x},
\sigma_{\overline{x'}}^{\overline{x}}(y')}^{ \overline{x'},y'}
(\, \mathsf{Ap}(\, \pi_1(w)\, ,\,   \sigma_{\overline{x'}}^{\overline{x}}(y') \,   )\, ),\, p\, \rangle$ for
$\overline{x}, \overline{x'}\in \Gamma^I$ and 
$w\in (\, B\rightarrow {\cal P}(1) \, )^I(\overline{x})$}
where {\small $p$} is a proof-term witnessing the preservation of equalities obtained from {\small $\pi_2(w)$.}
\\

\begin{remark}
\em
Note that  we need to close \mtt\ collections under strong
indexed sums  in order
  to interpret function collections  toward  {\small ${\cal P}(1)$},
and hence in turn to interpret power collections of sets as described in
section~\ref{ext}.
\end{remark}

After giving the interpretation of \emtt-signatures (and \emttdp-signatures)
 into \mtt-extensional dependent types and terms,
we can  interpret \emtt\ (and \emttdp) judgements in the category {\small \cqis} ({\small
\cqdp/$\simeq$}) by following the idea behind the naive interpretation of dependent types
in \cite{Seel} and
in the completeness theorem in \cite{tumscs}.
 To this purpose
we need first to transform extensional dependent types into
\cq-arrows as in proposition~\ref{dep}.
\begin{definition}\em
\label{ta}
Given an \mtt-context {\small $\Gamma_{=}\,$} we define its indexed closure
{\small $Sig( \Gamma_{=})$} as a \cq-object
together with  suitable projections $\pi^n_j(z)$
for {\small $z\in Sig( \Gamma_{=})$}  and $j=1,\dots, n$
by induction on the length $n$ of {\small $\Gamma_{=}\,$} as follows:
\begin{list}{-}{ }
\item
If {\small $\Gamma_=\, \equiv\, x\in A_=$} then
{\small $$Sig(\,  x\in A_=\, )\, \equiv (A, =_A)\qquad \mbox{ and }\qquad \pi_1^1(z)\, \equiv\, z$$}
\item
If {\small $\Gamma_=\, \equiv\, 
\Delta, x\in A$ } of $n+1$ length then
{\small $$Sig(\,\Delta_=,  x\in A_=\, )\, 
\equiv\, (\, \Sigma_{z\in Sig(\Delta_=)}\  A[x_1/ \pi_1^n(z),\dots
x_n/\pi_n^n(z)]\, ,\, =_{Sig}\, )$$}

\noindent
where {\small $w=_{Sig}w'\, \equiv\,  \exists_{d\in  \pi_1(w)=_{Sig}
\pi_1(w')}\ \
\sigma_{\pi_1(w)}^{\pi_1(w')}( \pi_2(w))=_{A_{w'}}
\pi_2(w')\\
\mbox{for }w,w'\in    \Sigma_{z\in Sig(\Delta_=)}\  A\,
  [x_1/ \pi_1^n(z),\dots
x_n/\pi_n^n(z)]\, $}  and 
{\small $\pi_j^{n+1}(w)\, \equiv\
\pi_j^n(\pi_1(w))$} for $j=1,\dots, n$ and
{\small $\pi_{n+1}^{n+1}(w)\, \equiv\, \pi_2(w)$}
and {\small $A_{w'}\, \equiv\, A\,
  [x_1/ \pi_1^n(\pi_1(w')),\dots
x_n/\pi_n^n(\pi_1(w'))]\, $}.
\end{list}
\end{definition}

\begin{definition}[Interpretation  in {\small \cqis}]
\label{intcat}\em
The interpretation of \emtt\ judgements in the category {\small \cqis}
$$Int: \emtt\rightarrow \cqis$$
is defined by using the interpretation of \emtt-signatures into \mtt-extensional dependent types and terms in appendix~\ref{emcq} as follows.

\noindent
 An \emtt-dependent type is interpreted as the projection in {\small \cqis}\  of its
interpretation as \mtt-extensional dependent type
according to the idea of turning a dependent collection
into an arrow in proposition~\ref{dep}:
{\small $$Int(\, B\ type\  [\,\Gamma\, ]\, )
\ \, \equiv\,[\pi_1]\, :\, [\, Sig(\, \Gamma^I_=, B^I_=\, )\, ]\rightarrow 
[\,Sig(\, \Gamma^I_=\, )\, ] $$}

\noindent
Then, an \emtt-type equality judgement is interpreted as the morphism
equality of  type interpretations 
in {\small \cqis}\
{\small
$$ Int(\,  A=B \ type\ [\Gamma]\, )\, \equiv\,
Int(\, A\ type\  [\,\Gamma\, ]\, ) =_{\cqis}Int(\, B\ type\  [\,\Gamma\, ]\, )
$$}

\noindent
which amounts to proving that their interpretations  as \mtt-extensional dependent types
are extensionally equal, namely that in \mtt\
 we can derive {\small $$ A^I_=\,=_{ext}\,
 B^I_= \  [\,\Gamma^I_=\, ]$$}

\noindent
An \emtt-term is interpreted as a section of  the corresponding type 
and it is built out of its interpretation as \mtt-extensional term:
{\small $$Int(\, b\in B\  [\Gamma]\, )
\, \equiv\,[\, \langle z, \overline{b^I}\rangle \, ]:[\,Sig(\, \Gamma^I_=\, )\, ]\rightarrow
[\, Sig(\, \Gamma^I_=, B^I_=\, )\, ]
$$
}

\noindent
where $\overline{b^I}$ is obtained by substituting its free variables
 {\small $x_1,\dots, x_n$} with $\pi_j^n(z)$ for {\small $j=1,\dots, n$} as in definition~\ref{ta}.  The given interpretation amounts
 to deriving in \mtt\ 
 {\small $$b^I\in_{ext} B^I\ [\Gamma^I]$$} 

\noindent
 Note that the interpretation
of a term is a section of its type interpretation because
 {\small $ [\pi_1]\cdot [ \langle z, \overline{b^I}\rangle ]=_{\cqis} \mathsf{id}$ } holds.

\noindent
Finally,
an \emtt-term equality judgement is interpreted as the equality of  {\small \cqis}-morphisms interpreting the terms:
{\small $$Int(\,  a=b\in B \ [\Gamma]\, ) \, \equiv\, 
Int(\,a\in B\ type\  [\,\Gamma\, ]\, ) =_{\cqis} Int(\,b\in B\ type\ 
 [\,\Gamma\, ]\, ) $$}

\noindent
which amounts to deriving in \mtt\
{\small $$a^I\, =_{ext}\, b^I\in_{ext} B^I\ [\Gamma^I]$$}

\noindent
Analogously, we define the interpretation of \emttdp-judgements in
{\small \cqdpis}.
\end{definition}

\noindent
In order to prove the validity theorem  we need to know how to interpret weakening and substitution.
For easiness we just show how to interpret substitution.

\noindent
Note that in the following, given
a context  {\small $\Gamma\, \equiv\, \Sigma, x_n\in A_n,\Delta $}
with {\small $\Delta\, \equiv\, x_{n+1}\in A_{n+1}, ..., x_{k}\in A_{k}$} then 
for every {\small $a\in A_{n}\ [\Sigma]$} and for any type
{\small $B\ type\ [ \Gamma] $} we simply write the substitution
 of $x_n$  with {\small $a$} in  {\small $B$}
in the form
{\small $ B[x_n/ a]\ type\  [\Sigma, \Delta_a]$} instead of the more correct form
 {\small $ B[x_{n}/a_{n}] [x_{i}/x'_{i}]_{i=n+1,...,k}\ 
type\  [\Sigma, \Delta_a]$}
where  {\small
$\Delta_a \  \equiv \ x'_{n+1}\in A'_{n+1}, ..., x'_{k}\in A'_{k}$}
and {\small  $A'_{j}\  \equiv \ A_{j}\ [x_{n}/a_{n}] 
[x_{i}/x'_{i}]_{i=n+1,...,j-1}$}
for {\small $j=n+2,... ,k$}, if {\small $n+2 \leq k$}, otherwise
{\small $A'_{n+1}\  \equiv \ A_{n+1}\ [x_{n}/a_{n}]$}.
If  {\small $\Delta $} is the empty context, then {\small $\Delta_a$} is empty, too.
 Similar abbreviations are used also for terms.
\begin{lemma}
\label{isocompl2}
For any \emtt\ judgement
{\small $B\ type\ [ \Gamma] $}  interpreted in {\small \cqis}\ as
{\small $$[\pi_1]\, :\, [\, Sig(\, \Gamma^I_{=}\,,\,  w\in B^I_=\,  )]
\rightarrow [\, Sig(\, \Gamma^I_{=}\,)\,]$$}

\noindent
and {\small $b\in  B\ [\Gamma]$} interpreted in {\small \cqis}\ as
{\small 
$$[\langle z, \overline{b^I} \rangle]\, :
 [\, Sig(\, \Gamma^I_{=}\, )\,]\rightarrow
 [\, Sig(\, \Gamma^I_{=}\, ,\,  w\in B^I_=\,  )]
$$}

\noindent
 {\em substitution} is interpreted as follows:
supposed {\small $\Gamma\, \equiv\, \Sigma, x_n\in A_n,\Delta $}
with {\small $\Delta\, \equiv\, x_{n+1}\in A_{n+1}, ..., x_{k}\in A_{k}$}
if not empty,  
for every \emtt\ judgement {\small $a\in A_{n}\ [\Sigma]$} interpreted as
{\small $$
[\, \langle z, \overline{a^I}\rangle\, ]\,:\, 
[\, Sig(\, {\Sigma^I}_=\, )\,]\rightarrow 
[\, Sig(\,   {\Sigma^I}_= ,
x_n\in {A^I_n}_=\, )\, ]
$$}

\noindent
then
 {\small $$\begin{array}{l}
I(\, B[x_n/ a]\ type\ [ \Sigma,\Delta_a]\, )\, =_{\cqis}\\
\qquad \qquad \qquad 
\,[\pi_1]\, :\,
[\, Sig(\,  {\Sigma^I}_=,{\Delta_a^I}_= \, ,\, w\in
B^I[x_n/\widetilde{a^I}]_=\, )\, ]\rightarrow 
[ \, Sig(\, {\Sigma^I}_=, {\Delta_a^I}_=\,)\, ]
\end{array}$$}
and

{\small $$\begin{array}{l}
 Int(\,b[x_n/ a]\in  B[x_n/a]\
 type\ [ \Sigma, \Delta_a]\, )\, =_{\cqis}\\
\qquad \qquad \qquad 
\,[\, \langle z, \overline{b^I[x_n/\widetilde{a^I}]}\rangle\,]:
[Sig(\, {\Sigma^I}_=,{\Delta_a^I}_=\, )\, ]\rightarrow 
[\, Sig(\, {\Sigma^I}_=,{\Delta_a^I}_=,w\in
B^I[x_n/\widetilde{a^I}]_=\, )\, ]
\end{array}$$}

\noindent
where the support of {\small $B^I[x_n/\widetilde{a^I}]_=$} is 
{\small $ B^I[x_n/\widetilde{a^I}][x_{i}/x'_{i}]_{i=n+1,...,k}$ } and
{\small $z=_{B^I [x_n/\widetilde{a^I}]}z'\,
 \equiv\, (z=_{B^I}z')[x_n/\widetilde{a^I}][x_{i}/x'_{i}]_{i=n+1,...,k}$}
defines its equality, if {\small $\Delta$} is  not empty.
 \end{lemma}
{\bf Proof.} By induction on the interpretation of the signature.
\medskip

\noindent
An analogous lemma holds for the interpretation of \emttdp-signatures in {\small
\cqdpis}.

\begin{theorem} [validity of \emtt\ into \cqis]
\label{intset}
The calculus \emtt\  is valid  with respect to the
 interpretation in definition~\ref{intcat} of
\emtt-signatures in {\small \cqis}:\\
If  {\small $  A\ type \ [\Gamma]$} is derivable in \emtt,
 then 
{\small  $Int(\,  A\ type\ [\Gamma]\, )$} is well defined.
If {\small  $ a\in A\ [\Gamma]$} is derivable in \emtt,  then 
{\small  $Int(\, a\in A\ [\Gamma]\, )$} is well defined.
Supposing that {\small $  A\ type \ [\Gamma]$} and {\small $ B\ type\ [\Gamma]$}
are derivable   in \emtt,
if {\small  $ A=B \ type\ [\Gamma]$} is derivable in \emtt,
then {\small $Int(  A=B \ type\ [\Gamma]\, ) $} is valid.\\
Supposing that {\small $ a\in A\ [\Gamma]$} and {\small  $ b\in A\ [\Gamma]$} are  
derivable in \emtt, if
 {\small $ a=b\in A \ [\Gamma]$} is derivable  in \emtt,
then  {\small $Int(\,  a=b\in A \ [\Gamma]\, )$} is also valid.
\end{theorem}
{\bf Proof.}
We can prove the statements by induction on the derivation of the judgements
by making use of theorem~\ref{mainth1}.

\noindent
Note that the rule \mbox{ \small conv)} is validated because of the presence of  canonical isomorphisms
witnessing that two types are extensionally equal. Moreover, conversion rules are valid
thanks to the properties of canonical isomorphisms in proposition~\ref{coer} and
thanks to the fact that they send canonical elements to canonical ones.
\medskip 

\noindent
Analogously, we can prove:
\begin{theorem} [validity of \emttdp\ into \cqdpis]
\label{intset}
The calculus \emttdp\  is valid  with respect to the
 interpretation in definition~\ref{intcat} 
of \emttdp-signatures in {\small \cqdpis}:\\
If  {\small $  A\ type \ [\Gamma]$} is derivable in \emttdp,
 then 
{\small  $Int(\,  A\ type\ [\Gamma]\, )$} is well defined.
If {\small  $ a\in A\ [\Gamma]$} is derivable in \emttdp,  then 
{\small  $Int(\, a\in A\ [\Gamma]\, )$} is well defined.
Supposing that {\small $  A\ type \ [\Gamma]$} and {\small $ B\ type\ [\Gamma]$}
are derivable   in \emttdp,
if {\small  $ A=B \ type\ [\Gamma]$} is derivable in \emttdp,
then {\small $Int(  A=B \ type\ [\Gamma]\, ) $} is valid.\\
Supposing that {\small $ a\in A\ [\Gamma]$} and {\small  $ b\in A\ [\Gamma]$} are  
derivable in \emttdp, if
 {\small $ a=b\in A \ [\Gamma]$} is derivable  in \emttdp,
then  {\small $Int(\,  a=b\in A \ [\Gamma]\, )$} is also valid.
\end{theorem}

\begin{remark}\em
Note that the power collection of a set is interpreted in \cq\ as the quotient
of suitable small propositional functions under equiprovability as in \cite{toolbox}, but
with the difference that here it is a particular construction on the top of \mtt,
while  in \cite{toolbox} it is declared to be a type in the underlying
intensional theory, which then loses the decidability of type judgements.
Hence, in our interpretation of \emtt\ over 
\mtt\ a subset is interpreted as the equivalence
class determined by a small propositional function, 
while in \cite{toolbox} (and also in \cite{finmin}) it is identified with it.
\end{remark}

\begin{remark}\rm {\bf Internal logic  of {\small \cqdp}}.
As already announced,
 {\small  \emtt\ } is not at all the internal language of {\small \cq}, because,
for example, 
it does not include generic quotient collections whilst they are supported
in the model. Even
{\small \emttdp\ } is not the internal language of {\small \cqdp}.
One reason is the following. If we restrict to the set-theoretic
fragment of {\small  \emtt}, called 
{\small  \qmtt\ } in \cite{mai07} and here {\small \emtts},
then
 the interpretation  of implication,
of universal quantification and of dependent product set do not seem 
to be  preserved by the functor 
{\small $\xi: \cqset\rightarrow  {\cal C}(\mbox{\emtts } )$}
sending an extensional type into its quotient (see also \cite{benphd}),
 where {\small ${\cal C}(\emtts)$} is the syntactic category of 
{\small \emtts\ } defined as in \cite{tumscs}.

However,   if we take the {\it set-theoretic coherent
fragment {\small \cemtt } } of {\small\emtts},
then we expect {\small\cemtt\ }
 to be an {\it internal language of the quotient
model} built out of {\it the corresponding set-theoretic coherent
 fragment {\small \cmtt\ } of  {\small \mtts}}.
The coherent fragment {\small \mbox{\rm \cemtt} } is obtained 
from the set-theoretic fragment {\small \emtts\ }
 by cutting out
implication, universal quantification and dependent product sets.
Also the fragment {\small \mbox{\rm \cmtt}} is obtained
from  {\small\mbox{\mtts\ } } in an analogous way. A proof of this is left to future work.
\end{remark}




\noindent
\begin{remark}\rm
{\bf Connection with the exact completion of a weakly lex category}.
The construction of total setoids on \mtt\ corresponds categorically to
an instance of a generalization of
 the exact completion construction~\cite{Carvit,excom}
of a weakly lex category. 

The connection with the construction of the exact completion in ~\cite{Carvit,excom}
is clearer if we build our quotient model over 
 Martin-L{\"o}f's type theory  {\small\mttps\ } 
in \cite{PMTT}.
Then, such a quotient model built over  {\small \mttps}, always
with total setoids as in definition~\ref{setoid} and called
{\small \cqac},
 happens to be equivalent to 
the  exact completion construction in \cite{Carvit,excom}
performed on the weakly lex syntactic category, called {\small \cset}, associated with 
 {\small \mttps\ }
as in remark~\ref{ros}.
In particular, {\small \cqac\ }  turns out to be a {\it list-arithmetic
locally cartesian closed pretopos}.

Now, it is important to note that the quotient model {\small \cqac\ }
does  not seem to validate well-behaved quotients, if we identify {\it propositions
as sets} as done in   {\small\mttps}. Indeed,  under this identification
\mbox{\small \cqac}, {\it  but also} {\small \mbox{\cq}}, 
  supports  first-order extensional Martin-L{\"o}f's type theory in
  \cite{ML84} and  hence it {\it validates
the axiom of choice}, as a consequence of the fact that
 universal and existential quantifiers are identified
with  dependent products $\Pi$ and strong indexed sums $\Sigma$ respectively.
Then,  effectiveness of quotients,
being generally incompatible with the axiom of choice (see \cite{Maieff}), 
does not seem to be validated.

To gain  well-behaved  quotients in {\small \cqac},  one possibility is to reason
by identifying {\it propositions with mono collections},
as well as {\it small propositions with mono sets} like in the logic of a pretopos
(see \cite{tumscs}). Instead, in {\small \cq}\ we get them by identifying 
propositions only with  some mono collections, as well as
 small propositions only with  some mono sets.  Categorically this means
that  {\small \cq}\ is closed under effective quotients of categorical equivalence relations obtained  via the properties of the
  comprehension adjunction as  ${\cal P}_q$-equivalence relations.
In fact, even if {\small \cq\ } supports quotients of all the monic equivalence relations,
these do not seem to enjoy effectiveness.
Always categorically speaking, this means that {\small \cq\ } does not seem to be
 a pretopos,
even if it has quotients for all monic equivalence relations.
Indeed, 
we are not able to prove that all the monic equivalence relations
are in bijection with ${\cal P}_q$-equivalence relations, for which
effective quotients exist (which explains
why we introduced the concept of ${\cal P}_q$-equivalence relation!).

Finally note that effective quotients
in {\small \cq\ } and {\small \cqset}, and also in {\small \cqdp},
    are enough to make these models
regular (see \cite{handJohn1} for the categorical definition): indeed
one can define the image of {\small $f: (A,=_A)\rightarrow
(B,=_B)$ } as the quotient of {\small $(A,=_A)$} over its kernel, namely as
 {\small $(\, A,\,  f(x)=_Bf(y)\, )$}, after noticing that monic arrows are indeed injective.

\begin{remark}\em
In order to interpret quotients in \mtt\ we also considered to mimic the exact completion on
a regular category in \cite{Carvit,eff}. But we ended up just in a list-arithmetic pretopos, for example not necessarily closed
under dependent products.
Indeed, given that in this completion arrows are identified with
functional relations, if the axiom of choice is not generally valid
as it happens in \mtt, in order to define exponentials and dependent products
 we would
 need to use an impredicative quantification on relations
that is not allowed in  
\mtt\ for its predicativity.
\end{remark}


\subsection{The axiom of choice is not valid in \emtt}
\label{acem}
The axiom of choice is not derivable in \emtt.
This is not surprising, if we consider that it may be incompatible
with effective quotients (see \cite{Maieff}).
To show this fact, 
just observe  that the {\it propositional axiom of choice} written in \emtt\ 
as follows
{\small $$(AC)\qquad\forall x\in A\  \exists y\in B\ R(x,y) \ \   \longrightarrow
\ \ \exists f\in A\rightarrow B\ \ \forall  x\in A\ R(x,\, \mathsf{Ap}(f,x)\, )$$}
is exactly interpreted in {\small  \cqac}, and also in {\small \cq}, as:  
{\small $$\forall x\in A^I\  \exists y\in B^I\ R^I(x,y) \ \   \longrightarrow
\ \ \exists f\in (A\rightarrow B)^I\ \ \forall  x\in A^I\ \ R^I(x,\mathsf{Ap}(f,x))$$}
where we recall that
{\small  $$ (A\rightarrow B)^I\, \equiv\, \Sigma_{h\in A^I\rightarrow B^I}\quad\   
\forall_{x_1\in A^I}\ \forall_{x_2\in A^I}\ \   x_1=_{A^I}x_2\ \rightarrow
\mathsf{Ap}(h,x_1)=_{B^I} \mathsf{Ap}(h,x_2)$$}

\noindent
This  interpretation in \cq\ amounts  to be exactly
 equivalent to the  {\it extensional axiom of choice} in \cite{MLac,Jes}
given that {\small $R^I$} satisfies the conditions:
\begin{list}{-} {}
\item If {\small $a=_{A^I}a'$} then {\small $\forall_{z\in B^I}\, R^I(a,z)\rightarrow R^I(a',z)$ }holds.
\item If {\small  $b=_{B^I}b'$} then {\small $\forall_{z\in A^I}\, R^I(z,b)\rightarrow R^I(z,b')$} holds.
\end{list}
Therefore, the arguments in \cite{MLac,Jes}
exactly  show that the propositional axiom of choice {\it fails} to
be valid in the  quotient models  {\small \cqac}, and even more in {\small\cq}\
or in an Heyting pretopos, and hence also in \emtt.
Indeed, we can prove that the validity of the axiom of choice
in \emtt\ yields that all propositions are decidable
as shown in  \cite{MLac,Jes, MV}, whose proof goes back to Goodman-Myhill's one in \cite{ac}.
To prove this, we use a choice property valid for effective quotients
thanks to the fact that  propositions are mono:
\begin{lemma}
\label{q}
In \emtt\  for any quotient set {\small $A/R\ set\  [\Gamma]$}
we can derive a proof of
{\small $$\forall_{z\in A/R}\ \exists_{y\in A}\ [y]=_{A/R} z$$}
\end{lemma}
{\bf Proof.}
Given  {\small $z\in A/R$}, by elimination of quotient sets we
get {\small $\mathsf{El}_{\cal Q}(z, \, (x).\mathsf{true}\, )\in
\exists_{y\in A}\ [y]=_{A/R} z$} because for {\small $x\in A$} then
{\small  $
\exists_{y\in A}\ [y]=_{A/R} [x]$} holds by reflexivity of Propositional
Equality
by taking $y$ as $x$ and it is well defined since propositions
are mono.
\medskip

\begin{prop}
In \emtt\ the validity of the axiom of choice  
{\small $$(AC)\qquad\forall x\in A\  \exists y\in B\ R(x,y) \ \   \longrightarrow
\ \ \exists f\in A\rightarrow B\ \forall  x\in A\ R(x,\mathsf{Ap}(f,x))$$}

\noindent
on all \emtt\ sets $A$ and $B$ and small relation {\small
 $R(x,y)\ prop_s\ [x\in A, y\in B]$} 
implies that all small propositions are decidable.
\end{prop}
{\bf Proof.} We follow the proof in \cite{Jes}.
Let us define the following equivalence relation
on the boolean set {\small $Bool\, \equiv\, \mathsf{N_1}+ \mathsf{N_1}$} whose elements are called
{\small $\mathsf{true}\, \equiv \mathsf{inl}(\star)$} and {\small $\mathsf{false}\, \equiv \mathsf{inr}(\star)$}:
given any proposition $P$ we put
{\small $$R(a,b)\, \equiv\, a=_{Bool}b \vee P$$}

\noindent
Then, thanks to lemma~\ref{q} in \emtt\ we can derive a proof of
{\small $$\forall_{z\in Bool/R}\ \exists_{y\in Bool}\  z=_{Bool/R} [y]$$}

\noindent
and by the validity of the axiom of choice
we get a proof of
{\small $$\exists f\in Bool/R\rightarrow Bool\ \  \forall  z\in Bool/R\ \ 
z=_{Bool/R}[\mathsf{Ap}(f, z)]$$}

\noindent
that amounts to being an injective arrow by definition.
Then, given that the equality in {\small  $Bool$} is decidable
(which follows by sum disjointness), we get that the equality in {\small$Bool/R$} is decidable, too.
Hence, {\small  $[\mathsf{true}]=_{Bool/R}[\mathsf{false}]\, \vee\, 
\neg [\mathsf{true}]=_{Bool/R}[\mathsf{false}]$} is derivable. Then,
by effectiveness also {\small  $R(\mathsf{true}, \mathsf{false})\, \vee\, 
\neg R(\mathsf{true}, \mathsf{false})$} is derivable, too.
Now, given
that {\small $R(\mathsf{true},\mathsf{false})\leftrightarrow P$}, then
{\small $P\, \vee \, \neg P$} is derivable, too,  namely $P$ is decidable.
The logic of small propositions is then classical.
\medskip

\noindent
Note that a similar argument holds also at the level of propositions.
Moreover, about pretopoi we can deduce the following:
\begin{cor}
A locally cartesian closed
 pretopos enjoying the validity of the propositional axiom of choice
is boolean.
\end{cor}
{\bf Proof.}
The internal logic of a locally cartesian closed pretopos devised in \cite{tumscs}
is an extension of the set-theoretic fragment of \emtt\ and validates
lemma~\ref{q}.
\medskip

In  {\small \cqac\ } the propositional axiom of choice survives at least for those quotients
whose equivalence relation is the propositional equality of {\small \mttps}.
Only the validity of the axiom of unique choice continues to hold
in its generality in {\small \cqac\ } (see also \cite{MLac}).

\end{remark}


\subsection{What links the two levels?}
A question we need to address in forming a two-level foundation
 is to decide what mathematical link should tie the two
levels. In \cite{mtt} we said that the extensional level
must be obtained from the intensional one
 by following Sambin's  forget-restore principle expressed in \cite{toolbox}.

Our example of two-level foundation fully satisfies such a principle.
In particular, the interpretation of \emtt\ into the quotient model over \mtt\
makes visible the validity of such a principle.
Indeed, in \emtt\ we work with undecidable judgements, 
while those of \mtt\ are decidable.
The interpretation of \emtt\ into \mtt\ restores the forgotten information
of \emtt\ undecidable judgements  by transforming them into decidable ones, once
the lost information has been recovered.
For example,  the \emtt-judgement {\small $B=C\ set\ [\Gamma]$} is interpreted as the existence of a canonical isomorphism 
{\small $$\begin{array}{cc}
\emtt& \mtt\\[5pt]
I(\, B=C\ set\  [\Gamma]\, )\, \equiv\, &
\mbox{ there exists a canonical isomorphism with components}\\[5pt]
&\quad \tau_{B^I}^{C^I}(y)\in C^I_=\ [\Gamma^I_= , y\in B^I_=]\qquad\mbox{ and }
\qquad \tau_{C^I}^{B^I}(z)\in B^I_=\ [ \Gamma^I_= , z\in
  C^I_=]
\end{array}$$}

\noindent
This says that in order to interpret \emtt-type equalities into \mtt\
we need to restore the missing canonical isomorphisms and hence
to prove suitable decidable \mtt\ judgements.

\noindent
Another
example, already studied in \cite{ML84,sienlec,Val97},  is the
interpretation
 of the validity of a proposition {\small $A \ prop\ [\Gamma]$},
 expressed by the \emtt\ judgement
{\small $\mathsf{true}\in A$}, as the existence of a proof-term:

{\small $$\begin{array}{cc}
\emtt& \mtt\\[5pt]
I(\,\mathsf{true}\in A\, )\, \equiv\, &
\mbox{ there exists }p\in A^I
\end{array}$$}
In other words to interpret the validity of a proposition we need to
restore a proof-term of its \mtt-interpretation.

Such considerations reveal
 that the link between our two levels is {\it  not} then
{\it  a   merely interpretation} of the extensional level into the intensional one. Instead,
 the extensional level is designed over the intensional one
 {\it only by forgetting information
about  equality between types and terms} that can be {\it restored},
and hence it is implemented  {\it
only via a quotient construction performed on the intensional level.}

Therefore, this kind of two-level foundation rules out examples of  two-level foundations
where the extensional
level is governed by a classical logic that is interpreted
in the intensional one, like \mtt,
  via a double-negation interpretation.

\section{Future work}
Given that \emtt\ is not at all the internal language of our quotient model {\small \cq}, we would like
to find out whether such an internal language exists in terms of an extensional dependent type
theory. In particular, it would be useful
to find   the internal language of the quotient model {\small \cqac}\
 built over Martin-L{\"o}f's type theory in an analogous way to
 {\small \cq}. An application of this would be to extend \emtt\ with
 occurrences of the axiom of choice that are constructively admissible.
Indeed, whilst the propositional axiom of choice is generally constructively incompatible with \emtt,
as recalled in section~\ref{acem},  there are extensional sets
on which we can apply the axiom of choice without losing 
constructivity, as advocated by Bishop in \cite{Bishop}.
For example, in the quotient model {\small \cqac\ }  over Martin-L{\"o}f's type theory the axiom of choice is valid on those extensional sets whose equivalence
relation is the identity relation, including, for example,
the extensional set of natural numbers. These extensional
sets are  actually copies of intensional sets at the extensional level.

We think that knowing the internal language of {\small \cqac}\ in terms of an extensional dependent type theory \`a la  Martin-L{\"o}f
would help to characterize such intensional sets in \emtt\ without stating their existence in a purely
axiomatic way in the style of Aczel's Presentation axiom in \cite{czf}.

\noindent
{\bf Acknowledgements:}
Our first thanks go to Giovanni Sambin for all his inspiring ideas and his
fruitful collaboration jointly with us on constructive foundations: in
fact this paper was originated within the
research developed with him starting from \cite{mtt}.

 We then thank  Per Martin-L{\"o}f and  Thomas Streicher for very useful discussions
about the existence of a realizability model for  Martin-L{\"o}f's intensional type theory
satisfying the formal Church thesis. We also thank Peter Aczel, Jesper Carlstr{\"o}m, Giovanni Curi, 
Pino Rosolini, Claudio Sacerdoti Coen and  Silvio Valentini for other useful discussions on the topics developed here. We thank Ferruccio Guidi and Olov Wilander for making comments on previous versions of this paper.
\bibliographystyle{alpha}
\bibliography{biblioexin}

\begin{thebibliography}{{Joh}02b}

\bibitem[Alt99]{altlics}
Thorsten Altenkirch.
\newblock Extensional equality in intensional type theory.
\newblock In {\em 14th Symposium on Logic in Computer Science}, pages 412 --
  420, 1999.

\bibitem[AMS07]{altob}
T.~{Altenkirch}, C.~{McBride}, and W.~{Swierstra}.
\newblock Observational equality, now!
\newblock In {\em PLPV '07: Proceedings of the 2007 workshop on Programming
  languages meets program verification}, pages 57--68, New York, NY, USA, 2007.
  ACM.

\bibitem[AR]{czf}
P.~{Aczel} and M.~{Rathjen}.
\newblock Notes on constructive set theory.
\newblock Mittag-Leffler Technical Report No.40, 2000/2001.

\bibitem[BCP03]{ven}
G.~{Barthes}, V.~{Capretta}, and O.~{Pons}.
\newblock Setoids in type theory.
\newblock {\em J. Funct. Programming}, 13(2):261--293, 2003.
\newblock Special issue on {"Logical frameworks and metalanguages"}.

\bibitem[BCRS98]{bcsr}
L.~{Birkedal}, A.~{Carboni}, G.~{Rosolini}, and D.~{Scott}.
\newblock Type theory via exact categories.
\newblock In {\em Thirteenth Annual IEEE Symposium on Logic in Computer Science
  (Indianapolis, IN, 1998)}, pages 188--198. IEEE Computer Soc., 1998.

\bibitem[{Bee}85]{Beeson}
M.J. {Beeson}.
\newblock {\em Foundations of constructive mathematics. Metamathematical
  studies.}, volume~6 of {\em Ergebnisse der Mathematik und ihrer Grenzgebiete
  (3)}.
\newblock Springer-Verlag, Berlin, 1985.

\bibitem[{Bis}67]{Bishop}
E.~{Bishop}.
\newblock {\em Foundations of {C}onstructive {A}nalysis}.
\newblock McGraw-Hill Book Co., 1967.

\bibitem[{Car}95]{Car}
A.~{Carboni}.
\newblock Some free constructions in realizability and proof theory.
\newblock {\em J. Pure Appl. Algebra}, 103:117--148, 1995.

\bibitem[{Car}03]{carl}
J.~{Carlstr{\"o}m}.
\newblock Subsets, quotients and partial functions in {M}artin-{L}{\"o}f's type
  theory.
\newblock In {\em Types for proofs and programs}, volume 2646 of {\em Lecture
  Notes in Comput. Sci}. Springer, Berlin, 2003.

\bibitem[{Car}04]{Jes}
J.~{Carlstr{\"o}m}.
\newblock {EM} + {Ext-} + {AC}int is equivalent to {AC}ext.
\newblock {\em Mathematical Logic Quarterly}, 50(3):236--240, 2004.

\bibitem[CC82]{excom}
A.~{Carboni} and R.~{Celia Magno}.
\newblock The free exact category on a left exact one.
\newblock {\em Journal of Australian Math. Soc.}, 33:295--301, 1982.

\bibitem[{Coq}90]{TC90}
T.~{Coquand}.
\newblock Metamathematical investigation of a calculus of constructions.
\newblock In P.~{Odifreddi}, editor, {\em Logic in Computer Science}, pages
  91--122. Academic Press, 1990.

\bibitem[CR00]{Ros}
A.~{Carboni} and G.~{Rosolini}.
\newblock Locally cartesian closed exact completions.
\newblock {\em J. Pure Appl. Algebra}, pages 103--116, 2000.
\newblock Category theory and its applications (Montreal, QC, 1997).

\bibitem[CS08]{finmin}
F.~{Ciraulo} and G.~{Sambin}.
\newblock Finiteness in a {M}inimalist {F}oundation.
\newblock In I.~Scagnetto M.~Miculan and F.~Honsell, editors, {\em Types for
  Proofs and Programs International Conference, TYPES 2007}, volume 4941 of
  {\em Lecture Notes in Computer Science}, pages 51--68. Springer Berlin, 2008.

\bibitem[CV98]{Carvit}
A.~{C}arboni and E.M. {V}itale.
\newblock Regular and exact completions.
\newblock {\em Journal of Pure and Applied Algebra}, 125:79--116, 1998.

\bibitem[dB91]{DEBRUIJN}
N.G. de~{B}ruijn.
\newblock Telescopic mapping in typed lambda calculus.
\newblock {\em Information and Computation}, 91:189--204, 1991.

\bibitem[{Dyb}96]{dyb}
P.~{Dybjer}.
\newblock Internal type theory.
\newblock In {\em TYPES '95}, volume 1158 of {\em LNCS}, pages 120--134, 1996.

\bibitem[GA06]{AG}
N.~{Gambino} and P.~{Aczel}.
\newblock The generalized type-theoretic interpretation of constructive set
  theory.
\newblock {\em Journal of Symbolic Logic}, 71(1):67--103, 2006.

\bibitem[GM78]{ac}
N.~{G}oodman and J.~{M}yhill.
\newblock Choice implies excluded middle.
\newblock {\em Z. Math. Logik Grundlag. Math.}, 24:461, 1978.

\bibitem[{Gui}]{gui}
F.~{Guidi}.
\newblock {The Formal System lambda-delta}.
\newblock To appear in ACM Transactions on Computational Logic.

\bibitem[HJP80]{tripos}
J.M.E. {H}yland, P.T. {J}ohnstone, and A.M. {P}itts.
\newblock Tripos theory.
\newblock {\em Bull. Austral. Math. Soc.}, 88:205--232, 1980.

\bibitem[{Hof}94]{Hofmann}
M.~{Hofmann}.
\newblock On the interpretation of type theory in locally cartesian closed
  categories.
\newblock In Proceedings of CSL'94, September 1994.

\bibitem[{Hof}97]{disttheshof}
M.~{Hofmann}.
\newblock {\em Extensional Constructs in Intensional Type Theory.}
\newblock Distinguished Dissertations. Springer, 1997.

\bibitem[{Hyl}82]{eff}
J.~M.~E. {Hyland}.
\newblock The effective topos.
\newblock In {\em The L.E.J. Brouwer Centenary Symposium (Noordwijkerhout,
  1981)}, volume 110 of {\em Stud. Logic Foundations Math.}, pages 165--216.
  North-Holland, Amsterdam-New York,, 1982.

\bibitem[{Jac}99]{jacobbook}
B.~{Jacobs}.
\newblock {\em Categorical Logic and Type Theory}, volume 141 of {\em Studies
  in Logic}.
\newblock Elsevier, 1999.

\bibitem[{Joh}02a]{handJohn1}
P.~T. {Johnstone}.
\newblock {\em Sketches of an elephant: a topos theory compendium. Vol. 1,},
  volume~43 of {\em Oxford Logic Guides}.
\newblock The Clarendon Press, Oxford University Press, New York,, 2002.

\bibitem[{Joh}02b]{handJohn2}
P.~T. {Johnstone}.
\newblock {\em Sketches of an elephant: a topos theory compendium. Vol. 2.},
  volume~44 of {\em Oxford Logic Guides}.
\newblock The Clarendon Press, Oxford University Press, New York,, 2002.

\bibitem[{L}aw69]{dialectica}
F.W. {L}awvere.
\newblock Adjointness in foundations.
\newblock {\em Dialectica}, 23:281--296, 1969.

\bibitem[{L}aw70]{eqhyp}
F.W. {L}awvere.
\newblock Equality in hyperdoctrines and comprehension schema as an adjoint
  functor.
\newblock {\em Proc. Sympos. Pure Math.}, XVII:1--14, 1970.

\bibitem[{Luo}94]{luo}
Z.~{Luo}.
\newblock {\em Computation and reasoning. A type theory for computer science.},
  volume~11 of {\em International Series of Monographs on Computer Science}.
\newblock The Clarendon Press, Oxford University Press, New York, 1994.

\bibitem[{Mac}71]{M71}
S.~{Mac Lane}.
\newblock {\em Categories for the working mathematician.}, volume~5 of {\em
  Graduate text in Mathematics}.
\newblock Springer, 1971.

\bibitem[{Mai}99]{Maieff}
M.~E. {Maietti}.
\newblock About effective quotients in {C}onstructive {T}ype {T}heory.
\newblock In W.~Naraschewski T.~Altenkirch and B.~Reus, editors, {\em Types for
  proofs and programs. International workshop, TYPES '98. Kloster Irsee,
  Germany, March 27-31. 1999}, volume 1657 of {\em Lectures Notes in Computer
  Science}, pages 164--178. Springer Verlag, 1999.

\bibitem[{Mai}05]{tumscs}
M.E. {Maietti}.
\newblock Modular correspondence between dependent type theories and categories
  including pretopoi and topoi.
\newblock {\em Mathematical Structures in Computer Science}, 15(6):1089--1149,
  2005.

\bibitem[{Mai}07]{mai07}
M.E. {Maietti}.
\newblock Quotients over {M}inimal {T}ype {T}heory.
\newblock In {\em Computation and Logic in the Real World- CiE 2007, Siena},
  volume 4497 of {\em LNCS}, pages 517--531. Springer, 2007.

\bibitem[{Mar}75]{modmar}
P.~{Martin-L{\"o}f}.
\newblock About models for intuitionistic type theories and the notion of
  definitional equality.
\newblock In {\em Proceedings of the Third Scandinavian Logic Symposium (Univ.
  Uppsala, Uppsala, 1973)}, volume~82 of {\em Stud. Logic Found. Math.}, pages
  81--109. North-Holland, Amsterdam, 1975.

\bibitem[{Mar}84]{ML84}
P.~{Martin-L\"{o}f}.
\newblock {\em Intuitionistic Type Theory, notes by G. Sambin of a series of
  lectures given in Padua, June 1980}.
\newblock Bibliopolis, Naples, 1984.

\bibitem[{Mar}85]{sienlec}
P.~{Martin-L{\"o}f}.
\newblock On the meanings of the logical constants and the justifications of
  the logical laws.
\newblock In {\em Proceedings of the conference on mathematical logic (Siena,
  1983/1984)}, volume~2, pages 203--281, 1985.
\newblock reprinted in: Nordic J. Philosophical Logic 1 (1996), no. 1, pages
  11--60.

\bibitem[ML06]{MLac}
P.~Martin-L{\"o}f.
\newblock 100 years of {Z}ermelo's axiom of choice:what was the problem with
  it?
\newblock {\em The Computer Journal}, 49(3):10--37, 2006.

\bibitem[MP00]{mptt}
I.~{Moerdijk} and E.~{Palmgren}.
\newblock Wellfounded trees in categories.
\newblock {\em Annals of Pure and Applied Logic}, 104(1-3):189--218, 2000.
\newblock Proceedings of the Workshop on Proof Theory and Complexity, PTAC'98
  (Aarhus).

\bibitem[MP02]{swopos}
I.~{Moerdijk} and E.~{Palmgren}.
\newblock Type theories, toposes and constructive set theory: predicative
  aspects of {AST}.
\newblock {\em Annals of Pure and Applied Logic}, 114(1-3):155--201, 2002.

\bibitem[MS05]{mtt}
M.E. {Maietti} and G.~{Sambin}.
\newblock {Toward a minimalist foundation for constructive mathematics}.
\newblock In {L. Crosilla and P. Schuster}, editor, {\em From Sets and Types to
  Topology and Analysis: Practicable Foundations for Constructive Mathematics},
  number~48 in {Oxford Logic Guides}, pages 91--114. {Oxford University Press},
  2005.

\bibitem[MV99]{MV}
M.E. {Maietti} and S.~{Valentini}.
\newblock Can you add powersets to {M}artin-{L}{\"o}f intuitionistic type
  theory?
\newblock {\em Mathematical Logic Quarterly}, 45:521--532, 1999.

\bibitem[NPS90]{PMTT}
B.~{N}ordstr{\"o}m, K.~{P}etersson, and J.~{S}mith.
\newblock {\em Programming in Martin L{\"o}f's Type Theory.}
\newblock Clarendon Press, Oxford, 1990.

\bibitem[{Pal}05]{notepal}
E.~{Palmgren}.
\newblock Bishop's set theory.
\newblock Slides for lecture at the TYPES summer school, 2005.

\bibitem[{Pit}00]{Pitts}
A.M. {Pitts}.
\newblock Categorical logic.
\newblock In Oxford~University Press, editor, {\em Handbook of Logic in
  Computer Science}, volume~5, pages 39--128, 2000.

\bibitem[{Sam}03]{Sam2002}
G.~{Sambin}.
\newblock {Some points in formal topology}.
\newblock {\em Theoretical Computer Science}, 305:347--408, 2003.

\bibitem[{Sam}09]{gsmin}
G.~{Sambin}.
\newblock A minimalist foundation at work.
\newblock In P.~{Clark}, D.~{DeVidi}, and M.~{Hallett}, editors, {\em Essays in
  honour of John L. Bell}, Western Ontario Series in the Philosophy of Science.
  Springer, 2009.
\newblock To appear.

\bibitem[{Sam}ar]{BP}
G.~{Sambin}.
\newblock {\em The {B}asic {Picture}. {S}tructures for constructive topology.}
\newblock Oxford University Press, To appear.
\newblock (including two papers with Per Martin-Loef and with Venanzio
  Capretta). To appear.

\bibitem[{See}83]{hyper}
R.~A.~G. {Seely}.
\newblock Hyperdoctrines, natural deduction and the beck condition.
\newblock {\em Zeitschr. f. math. Logik. und Grundlagen d. Math.}, 29:505--542,
  1983.

\bibitem[{See}84]{Seel}
R.~{Seely}.
\newblock Locally cartesian closed categories and type theory.
\newblock {\em Math. Proc. Cambr. Phyl. Soc.}, 95:33--48, 1984.

\bibitem[Str91]{Streicher}
Th. Streicher.
\newblock {\em Semantics of type theory.}
\newblock Birkh{\"{a}}user, 1991.

\bibitem[SV98]{toolbox}
G.~{Sambin} and S.~{Valentini}.
\newblock Building up a toolbox for {M}artin-{L}{\"o}f's type theory: subset
  theory.
\newblock In G.~{Sambin} and J.~{Smith}, editors, {\em Twenty-five years of
  constructive type theory, {P}roceedings of a {C}ongress held in {V}enice,
  {O}ctober 1995}, pages 221--244. Oxford U. P., 1998.

\bibitem[{Tro}77]{tr}
A.S. {Troelstra}.
\newblock A note on non-extensional operations in connection with continuity
  and recursiveness.
\newblock {\em Indag. Math.}, 39(5):455--462, 1977.

\bibitem[Tv88]{DT88}
A.~S. {Troelstra} and D.~{van Dalen}.
\newblock Constructivism in mathematics, an introduction, vol. {I}.
\newblock In {\em Studies in logic and the foundations of mathematics}.
  North-Holland, 1988.

\bibitem[TvD88]{constrII}
A.~{Troelstra} and D.~{v}an {Dalen}.
\newblock {\em Constructivism in mathematics. Vol. {II} An Introduction.},
  volume 123 of {\em Studies in {L}ogic and the {F}oundations of
  {M}athematics}.
\newblock North-Holland Publishing Co. (Amsterdam), 1988.

\bibitem[{Val}95]{Val97}
S.~{Valentini}.
\newblock The forget-restore principle: a paradigmatic example.
\newblock In J.~Smith G.~Sambin, editor, {\em Twenty Five Years of Constructive
  Type Theory}, pages 275--283. Oxford Science Publications, Venice, 1995.

\bibitem[vdB06]{benphd}
B.~van~den {Berg}.
\newblock {\em Predicative topos theory and models for constructive set
  theory}.
\newblock PhD thesis, University of Utrecht, 2006.

\end{thebibliography}

\section{Appendix: The typed calculus \mtt\ }
\label{mttsyn}
 We present here the inference rules to build types
 in \mtt. 
The inference rules involve judgements  written in the style of Martin-L{\"o}f's type theory
\cite{ML84,PMTT} that may be of the form:\\

{\small  $$A \ type \ [\Gamma] \hspace{.5cm} A=B\ type\ [\Gamma] 
\hspace{.5cm} a \in A\ 
 [\Gamma] \hspace{.5cm} a=b \in A\ [\Gamma] $$}
where types include collections, sets, propositions and small propositions,
namely 
{\small $$type \in \{ col, set,prop,prop_s\, \}$$}

For easiness, the piece of context common to all judgements
involved in a rule
is omitted and
 typed variables appearing in a context
are meant to be added to the implicit
context as the last one.

\noindent
 Note that to write
 the elimination constructors of our types
 we adopt the higher-order syntax  in
\cite{PMTT}~\footnote{For example, note that the elimination constructor of
disjunction ${\bf \it El}_{\vee}(w,a_{B},a_{C})$  binds the open terms
$a_{B}(x)\in A\ [x\in B]$ and 
  $      a_{C}(y)\in A \ [y\in C]$.
Indeed, given that they are needed in the disjunction conversion rules,
it follows that these open terms
 must be  encoded into 
the elimination constructor.
To encode them we use the higher-order syntax as  in \cite{PMTT} (see also \cite{gui}).
According to this syntax
the open term $a_B(x)\in A\ [x\in B]$ yields to  $(x\in B)\,  a_B(x)$
of higher type $(x\in B)\, A$. Then, by $\eta$-conversion among higher types,
 it follows that  $(x\in B)\,  a_B(x)$ is equal to $a_B$. Hence, we often simply write the short expression
$a_B$ to recall the open term  where it comes from.}.

\noindent
We also have a form of judgement to build contexts:
{\small $$\Gamma\ cont$$}

\noindent
whose rules are the following

{\small
$$ \begin{array}{ll}
\emptyset\ cont\qquad  &\qquad \mbox{F-c}\displaystyle{
 \frac{\displaystyle\  A\ type\  [\Gamma]\ }
{\displaystyle \ \Gamma, x\in A \ cont\ }}\  (x \in A \not\in \Gamma)
\end{array}
$$}

\noindent 
Then, the first rule to build elements of type is 
the assumption of variables:

{\small
$$
\mbox{ var) }
\displaystyle{\frac{\displaystyle\  \Gamma, x\in A, \Delta \hspace{.3cm}\ cont\ }{\displaystyle\  x\in A\ [ \Gamma,x\in A , \Delta]\ }}
$$
}

\noindent
Among types there are the following embeddings: sets are collections and propositions are collections
\\
 
\noindent
{\small
$\begin{array}{l}
\mbox{\bf set-into-col) }\ \ \displaystyle{ \frac
       {\displaystyle\  A
        \hspace{.1cm} set\ }
{ \displaystyle\  A
        \hspace{.1cm} col\ }}
\end{array}
\qquad\qquad
\begin{array}{l}
    \mbox{\bf prop-into-col) }\ \ 
\displaystyle{ \frac
       {\displaystyle\  A\ prop\  }
      { \displaystyle\  A\ col  \ }}
\end{array}
$}
\\





\noindent
Moreover, collections are closed under strong indexed sums:\\

\noindent
{\small
$\begin{array}{l}
      \mbox{ \bf Strong Indexed Sum } \\[10pt]
      \mbox{F-}\Sigma ) \ \
\displaystyle{ \frac{\displaystyle   C(x)
         \hspace{.1cm} \ col \ [x\in B]}
         {\displaystyle \Sigma_{x\in B} C(x)\hspace{.1cm} col }}
         \qquad
      \mbox{I-}\Sigma )\ \
\displaystyle{ \frac
         {\displaystyle b\in  B \hspace{.3cm} c\in C(b)\qquad C(x)\ col\ [x\in B]}
         {\displaystyle \langle b,c\rangle\in  \Sigma_{x\in B} C(x)}}
\\[15pt]
\mbox{E-}\Sigma )\ \
\displaystyle{ \frac
         {\displaystyle \begin{array}{l}
M(z)\ col \ [ z\in \Sigma_{x\in B} C(x)]\\
d\in  \Sigma_{x\in B} C(x) \hspace{.3cm}
 m(x,y)\in M(\langle x, y \rangle)\ [x\in B,
         y\in C(x)]
\end{array}}
      {\displaystyle {\it  El}_{\Sigma}(d,m)\in  M(d)}}
      \\[15pt]
\mbox{C-}\Sigma ) \ \
\displaystyle{ \frac
         {\displaystyle \begin{array}{l}
M(z)\ col \ [ z\in \Sigma_{x\in B} C(x)]\\
b\in B\ \ \  c\in C(b)\hspace{.3cm} m(x,y)\in M(\langle x,y \rangle)\ [x\in B,
         y\in C(x)]\end{array}}
      {\displaystyle {\bf \it El}_{\Sigma}(\, \langle b, c\rangle
,m\, )=m(b,c)\in M(\langle b,c\rangle)} }
      \end{array}$
\\
\\

}

\noindent
Sets are generated  as follows:
\\
\\

\noindent
{\small
$\begin{array}{l}
      \mbox{ \bf Empty set} \\
      \mbox{ F-Em)}\ \ \mathsf{N_0} \hspace{.1cm} set \qquad
\mbox{ E-Em)}\ \
\displaystyle{ \frac
         {\displaystyle a\in  \mathsf{N_0} \hspace{.3cm} A(x)\hspace{.1cm} col \
[x\in \mathsf{N_0}] }
         {\displaystyle \orig{emp_{o}}(a)\in A(a)}}
      \end{array}$
\\
\\

\noindent
$\begin{array}{l}
\mbox {\bf Singleton}\\
 \mbox{\small S)}\ \orig{\mathsf{N_1}} \hspace{.1cm} set
 \qquad
 \mbox{\small I-S)}\ 
 \orig{\star} \in\orig{\mathsf{N_1}}
 \qquad 
 \mbox{\small E-S)}\
\displaystyle{  \frac
    {t\in \orig{\mathsf{N_1}}\quad  M(z)\ col\ 
[z\in \mathsf{\mathsf{N_1}}] \quad c\in M(\star)}
    {{ \it El}_{ \orig{\mathsf{N_1}} }(t,c)\in  M(t)}}
\qquad 
 \mbox{\small C-S)}\
\displaystyle{  \frac
    { M(z)\ col\ 
[z\in \mathsf{\mathsf{N_1}}] \quad c\in M(\star)}
    {{ \it El}_{ \orig{\mathsf{N_1}} }(\star,c)=c\in  M(\star)}}
 \end{array}$
\\
\\

\noindent
$\begin{array}{l}
      \mbox{ \bf Strong Indexed Sum set } \\[10pt]
      \mbox{F-}\Sigma_s ) \ \
\displaystyle{ \frac{\displaystyle   C(x)
         \hspace{.1cm} set\ \ [x\in B]\qquad B\ set}
         {\displaystyle \Sigma_{x\in B} C(x)\hspace{.1cm} set }}
      \end{array}$
\\
\\

\noindent
      $\begin{array}{l}
\mbox{\bf List set} \\[10pt]
      \mbox{F-list)}\
\displaystyle{ \frac
         { \displaystyle C \hspace{.1cm} set}
         {\displaystyle List(C) \hspace{.1cm} set }}
      \qquad
       \mbox{${\rm I}_{1}$-list)}\ \
\displaystyle{ \frac
         {\displaystyle \quad List(C) \hspace{.1cm} set\  }
         {\epsilon \in List(C)}}
      \qquad
      \mbox{${\rm I}_{2}$-list)}\ \
\displaystyle{ \frac
         {\displaystyle s\in List(C) \hspace{.3cm} c\in  C}
         {\displaystyle \orig{cons}(s,c)\in List(C)}}
      \end{array}$
\\
\\

\noindent
$\begin{array}{l}
\mbox{E-list)}\ \
\displaystyle{ \frac
         {\displaystyle \begin{array}{l}
L(z)\ col\  [z\in List(C)]
\hspace{.3cm}s\in List(C) \hspace{.3cm}\qquad  a\in L(\epsilon)\hspace{.3cm}\\
    l(x,y,z)\in
         L(\orig{cons}(x,y))  \ [x\in List(C),y\in C, z\in L(x)]
\end{array}}
      {\displaystyle {\bf \it El}_{List}(s,a, l)\in  L(s)}}
      \end{array}$
\\
\\

\noindent
$\begin{array}{l}
\mbox{${\rm C}_{1}$-list)}\ \
\displaystyle{  \frac
         {\displaystyle\begin{array}{l}
L(z)\ col\  [z\in List(C)] \hspace{.3cm}\qquad
 a\in L(\epsilon)\hspace{.3cm}\\
     l(x,y,z)\in
         L(\orig{cons}(x,y))  \ [x\in List(C),y\in C, z\in L(x)]
\end{array}}
      {\displaystyle {\bf \it El}_{List}( \epsilon, a,l)=a\in  L(\epsilon)}}
      \\[15pt]
\mbox{ ${\rm C}_{2}$-list)}\ \
\displaystyle{ \frac
         {\displaystyle\begin{array}{l}
L(z)\ col\  [z\in List(C)]
\hspace{.3cm} s\in List(C) \hspace{.3cm}c\in C \hspace{.3cm} a\in
L(\epsilon)\hspace{.3cm} \\
l(x,y,z)\in
         L(\orig{cons}(x,y))  \ [x\in List(C),y\in C, z\in L(x)]
\end{array}}
      {\displaystyle {\bf \it El}_{List}(\orig{cons}(s,c),a,l)=l(s,
c,{\bf \it El}_{List}(s,a,
      l))\in  L(\orig{cons}(s,c))}}
\end{array}$
\\
\\

\noindent
      $\begin{array}{l}
      \mbox{\bf Disjoint Sum set } \\[10pt]
      \mbox{F-+)} \ \
\displaystyle{  \frac
         { \displaystyle B \hspace{.1cm} set \hspace{.3cm}C
         \hspace{.1cm} set}
         {\displaystyle B+ C \hspace{.1cm} set }}
      \qquad
      \mbox{${\rm I}_{1}$-}+ ) \ \
\displaystyle{  \frac
         {\displaystyle b\in  B\qquad B\ set \qquad C\  set }
         {\displaystyle\orig{inl}(b)\in B+ C}}
      \qquad
\mbox{${\rm I}_{2}$-}+ ) \ \
\displaystyle{  \frac
         {\displaystyle c\in  C \qquad B\ set\qquad C\ set}
         {\displaystyle  \orig{inr}(c)\in B+ C}}
      \\[15pt]
\mbox{E-} + ) \ \
\displaystyle{  \frac
         {\displaystyle\begin{array}{l}
A(z)\ col\  [z\in B+C]\\
      w\in B+ C \hspace{.3cm}
      a_{B}(x)\in A(\orig{inl}(x))\ [x\in B]
      \hspace{.3cm}\
      a_{C}(y)\in A(\orig{inr}(y))\ [y\in C]
\end{array}}
         {\displaystyle { \it El}_{+}(w,a_{B},a_{C})\in A(w)}}
      \\[15pt]
      \mbox{${\rm C}_{1}$-}+ ) \ \
\displaystyle{ \frac
       {\displaystyle \begin{array}{l}
A(z)\ col\  [z\in B+C]\\
b\in B \hspace{.3cm}  a_{B}(x)\in A(\orig{inl}(x))\ [x\in B]
      \hspace{.3cm}\
      a_{C}(y)\in A(\orig{inr}(y))\ [y\in C]\end{array}}
      {\displaystyle { \it El}_{+}(\orig{inl}(b),a_{B},a_{C})=a_{B}(b)
\in A(\orig{inl}(c))}}
\\[15pt]
      \mbox{${\rm C}_{2}$-}+ ) \ \
\displaystyle{ \frac
       {\displaystyle\begin{array}{l}
A(z)\ col\  [z\in B+C]\\
      c\in C\hspace{.3cm} a_{B}(x)\in A(\orig{inl}(x))\ [x\in B]
      \hspace{.3cm}\
      a_{C}(y)\in A(\orig{inr}(y))\ [y\in C]\end{array}}
      {\displaystyle {\bf \it El}_{+}(\orig{inr}(c),
a_{B},a_{C})=a_{C}(c)\in A(\orig{inr}(c))}}
      \end{array}$
\\
\\

\noindent
      $\begin{array}{l}
\mbox{\bf Dependent Product set}
\\[10pt]
\mbox{F-$\Pi$)}\ \
\displaystyle{\frac{ \displaystyle  C(x)
        \hspace{.1cm} set\ [x\in B] \qquad B\ set }
{\displaystyle \Pi_{x\in B} C(x)\hspace{.1cm} set }}
 \qquad
 \mbox{ I-$\Pi$)}\ \
 \displaystyle{\frac{ \displaystyle c(x)\in  C(x)\ [x\in B] \qquad C(x)
        \hspace{.1cm} set\ [x\in B] \qquad B\ set }
 { \displaystyle \lambda x^{B}.c(x)\in \Pi_{x\in B} C(x)}}
       \\[15pt]
 \mbox{  E-$\Pi$)}\ \
 \displaystyle{\frac{ \displaystyle
 b\in B \hspace{.3cm} f\in \Pi_{x\in B} C(x) }
 {\displaystyle \mathsf{Ap}(f,b)\in C(b) }}
\\[15pt]
 \mbox{$\mathbf{\beta}$C-$\Pi$)}\ \
 \displaystyle{\frac{ \displaystyle b\in B \hspace{.3cm} 
c(x)\in  C(x)\ [x\in B]\qquad C(x)
        \hspace{.1cm} set\ [x\in B] \qquad B\ set}
 {\displaystyle \mathsf{Ap}(\lambda x^{B}.c(x),b)=c(b)\in   C(b) }}
      \end{array}
$}
\\
\\



\noindent
Propositions are generated as follows:
\\

\noindent {\small
      $\begin{array}{l}
       {\bf Falsum} \\
      \mbox{ F-Fs)}\ \ \bot \hspace{.1cm} prop \qquad
\mbox{E-Fs)}\ \
\displaystyle{ \frac
         {\displaystyle a\in  \bot \hspace{.3cm} A\hspace{.1cm} prop }
         { \displaystyle \orig{r_{o}}(a)\in A}}
      \end{array}$
\\
\\

\noindent
      $\begin{array}{l}
\mbox{\bf Disjunction}
      \\[10pt]
      \mbox{F-}\vee )\ \
\displaystyle{ \frac
         { \displaystyle B \hspace{.1cm} prop \hspace{.3cm}C
         \hspace{.1cm} prop}
         {\displaystyle B\vee C \hspace{.1cm} prop }}
      \quad
      \mbox{${\rm I}_{1}$-}\vee) \ \
\displaystyle{ \frac
         {\displaystyle b\in  B\qquad B\ prop\qquad C \ prop}
         {\displaystyle \mathsf{inl}_{\vee}(b)\in B\vee C}}
      \quad
\mbox{${\rm I}_{2}$-}\vee ) \ \
\displaystyle{ \frac
         {\displaystyle c\in  C \qquad B\ prop\qquad C\ prop}
         {\displaystyle \orig{inr}_{\vee}(c)\in B\vee C}}
\end{array}$\\[15pt]
$\begin{array}{l}
\mbox{E-} \vee ) \ \
\displaystyle{ \frac
         {\displaystyle \begin{array}{l}
A\  prop\ \\
      w\in B\vee C \hspace{.3cm}
      a_{B}(x)\in A\ [x\in B]
      \hspace{.3cm}\
      a_{C}(y)\in A \ [y\in C]\end{array}}
         {\displaystyle {\it El}_{\vee}(w,a_{B},a_{C})\in A}}
 \\[15pt]
      \mbox{${\rm C}_{1}$-}\vee ) \ \
\displaystyle{ \frac
       {\displaystyle \begin{array}{l}
A\  prop\qquad B\ prop\qquad C\ prop\\
b\in B \hspace{.3cm}  a_{B}(x)\in A\ [x\in B]
      \hspace{.3cm}\
      a_{C}(y)\in A \ [y\in C]
\end{array}}
      { \displaystyle {\bf \it El}_{\vee}
(\orig{inl}_{\vee}(b),a_{B},a_{C})=a_{B}(b)\in A}}
\\[15pt]
      \mbox{${\rm C}_{2}$-}\vee ) \ \
\displaystyle{ \frac
       {\displaystyle \begin{array}{l}
A\  prop\   \qquad B\ prop\qquad C\ prop  \\
c\in C\hspace{.3cm} a_{B}(x)\in A\ [x\in B]
      \hspace{.3cm}\
      a_{C}(y)\in A \ [y\in C]\end{array}}
      {\displaystyle {\bf \it El}_{\vee}
(\orig{inr}_{\vee}(c),a_{B},a_{C})=a_{C}(c)\in A}}
\end{array}
$
\\
\\

\noindent
$\begin{array}{l}
\mbox{\bf Conjunction}
\\[10pt]
      \mbox{F-}\wedge) \ \
\displaystyle{ \frac
         { \displaystyle  B   \hspace{.1cm} prop\qquad  C
         \hspace{.1cm} prop }
         {\displaystyle  B\wedge C\hspace{.1cm} prop }}
         \qquad
      \mbox{I-}\wedge )\ \
\displaystyle{ \frac
         {\displaystyle b\in  B \hspace{.3cm} c\in C\qquad B\ prop\qquad C\ prop}
         {\displaystyle \langle b,_{\wedge} c \rangle\in B\wedge C  }}
\\[15pt]
\mbox{$\mathbf{\textrm{ E}_{1}}$-$\wedge$) }\ \
\displaystyle{ \frac
       {\displaystyle d\in  B \wedge C}
{\displaystyle \pi_{1}^{B}(d)\in  B}}
\qquad \qquad
\mbox{$\mathbf{\textrm{ E}_{2}}$-$\wedge$) }\ \
\displaystyle{ \frac
       {\displaystyle d\in  B \wedge C}
{\displaystyle \pi_{2}^{C}(d)\in  C}}
\\[15pt]
\mbox{$\beta_{1}$ C-$ \wedge$) }\ \
\displaystyle{ \frac
       {\displaystyle b\in  B \hspace{.3cm} c\in C\qquad B\ prop\qquad C\ prop}
{\pi_{1}^{B}(\langle b,_{\wedge} c \rangle)=b\in B}}
\qquad\qquad
\mbox{$\beta_{2}$ C-$\wedge$) }\displaystyle{ \frac
       {\displaystyle b\in  B \hspace{.3cm} c\in C\qquad B\ prop\qquad C\ prop}
{\displaystyle \pi_{2}^{C}(\langle b,_{\wedge} c \rangle)=c\in C}}
\end{array}
$
\\
\\

\noindent
$
\begin{array}{l}
\mbox{\bf Implication}
\\[10pt]
\mbox{F-$\rightarrow$) }\ \
\displaystyle{ \frac
       {\displaystyle B
        \hspace{.1cm} prop \qquad C
        \hspace{.1cm} prop  }
{B \rightarrow C \hspace{.1cm} prop }}
\\[15pt]
\mbox{I-$\rightarrow$) }\ \
\displaystyle{ \frac
{\displaystyle
      c(x)\in  C\ [x\in B] \qquad B \hspace{.1cm} prop\qquad C\ prop\  }
{\displaystyle \lambda_\rightarrow x^{B}.c(x)\in B \rightarrow C}}
        \qquad
\mbox{ E-$\rightarrow$) }\ \
\displaystyle{ \frac
{\displaystyle b\in B \hspace{.3cm} f\in B\rightarrow C }
{\displaystyle \mathsf{Ap}_\rightarrow(f,b)\in C }}
       \\[15pt]
\mbox{$\mathbf{\beta}$C-$\rightarrow$) }\ \
\displaystyle{ \frac
{\displaystyle b\in B \hspace{.3cm} c(x)\in  C
 \
[x\in B]\qquad B \hspace{.1cm} prop\qquad C\ prop}
{ \displaystyle \mathsf{Ap}_\rightarrow
(\lambda_\rightarrow x^{B}.c(x),b)=c(b)\in   C }}
      \end{array}$
\\
\\

\noindent
      $\begin{array}{l}
      \mbox{\bf Existential quantification}
      \\[10pt]
      \mbox{F-}\exists ) \ \
\displaystyle{ \frac
         { \displaystyle   C(x)
         \hspace{.1cm} prop \ \ [x\in B]}
         {\displaystyle  \exists_{x\in B} C(x)\hspace{.1cm} prop }}
         \qquad
      \mbox{I-}\exists )\ \
\displaystyle{ \frac
         {\displaystyle b\in  B \hspace{.3cm} c\in C(b)\qquad C(x)\ prop\ [x\in B]}
         {\displaystyle \langle b,_{\exists} c \rangle\in  \exists_{x\in B} C(x)}}
      \\[15pt]
\mbox{E-}\exists )\ \
\displaystyle{  \frac
         {\displaystyle \begin{array}{l}
M\ prop\ \\
d\in  \exists_{x\in B} C(x) \hspace{.3cm} m(x,y)\in M\ [x\in B,
         y\in C(x)]
\end{array}}
      {\displaystyle { \it El}_{\exists}(d,m)\in  M}}
      \\[15pt]
\mbox{C-}\exists ) \ \
\displaystyle{ \frac
         {\displaystyle \begin{array}{l}
M\ prop\qquad C(x)\ prop\ [x\in B]\\
b\in B\ \ \  c\in C(b)\hspace{.3cm} m(x,y)\in M\ [x\in B,
         y\in C(x)]
\end{array}}
      {\displaystyle { \it El}_{\exists}(\langle b,_{\exists}c \rangle,m)=m(b,c)\in M} }
      \end{array}$
\\
\\

\noindent
      $\begin{array}{ll}
\mbox{\bf Universal quantification}
\\[10pt]
\mbox{F-$\forall$) }\ \
\displaystyle{ \frac
{\displaystyle C(x)\hspace{.1cm} prop\ [x\in B]  }
{\displaystyle \forall_{x\in B} C(x)\hspace{.1cm} prop }}
&
\mbox{I-$\forall$) }\ \
\displaystyle{ \frac
{\displaystyle c(x)\in  C(x)\ [x\in B]\qquad  C(x)\hspace{.1cm} prop\ [x\in B]  }
{\lambda_\forall x^{B}.c(x)\in \forall_{x\in B} C(x)}}
        \\[15pt]
\mbox{E-$\forall$) }\ \
\displaystyle{ \frac
{\displaystyle b\in B \hspace{.3cm} f\in \forall_{x\in B} C(x) }
{\displaystyle \mathsf{Ap}_{\forall}(f,b)\in C(b) }}
       &
\mbox{$\mathbf{\beta}$C-$\forall$) }\ \
\displaystyle{ \frac
{\displaystyle b\in B \hspace{.3cm} c(x)\in  C(x)\ [x\in B]\qquad
 C(x)\ prop\  [x\in B]
}
{\displaystyle \mathsf{Ap}_{\forall}(\lambda_\forall
      x^{B}.c(x),b)=c(b)\in   C(b) }}
      \end{array}$
\\
\\

\noindent
      $\begin{array}{l}
\mbox{\bf Propositional Equality}
      \\[15pt]
      \mbox{F-Id)}\ \
\displaystyle{ \frac
{\displaystyle  A \hspace{.1cm} col \hspace{.3cm}  a\in A \hspace{.3cm} b\in A}
         {\displaystyle \orig{Id}(A, a, b)\hspace{.1cm} prop }}
      \qquad
      \mbox{I-Id)}\ \
\displaystyle{ \frac
{\displaystyle  a \in A}
         {\displaystyle \orig{id_{A}}(a)\in \orig{Id}(A, a, a)}}
\\[15pt]
      \mbox{E-Id)}\ \
      \displaystyle{\frac
       {\displaystyle \begin{array}{l}C(x,y)\ prop \ [x\in A,y\in A]\\
a\in A\quad b\in A\quad p\in \orig{Id}(A, a, b)
\quad c(x)\in C(x,x)
\ [x\in A]\end{array}}
      {\displaystyle {\it El}_{\mathsf{Id}}(p,c)\in C(a,b) }}
      \\[15pt]
       \mbox{C-Id)} \ \
\displaystyle{\frac
       {\displaystyle \begin{array}{l}C(x,y)\ prop \ [x\in A,y\in A]\\
a\in A
\quad c(x)\in C(x,x)
\ [x\in A]\end{array}}
      {\displaystyle {\bf \it El}_{\mathsf{Id}}(\mathsf{id}_A(a),c)=c(a)
\in C(a,a) }}
      \end{array}$
\\
\\
}

\noindent
Then, small propositions are generated as follows:
\\

\noindent
{\small
$\begin{array}{l} 
\bot \ prop_s\qquad
\displaystyle{ \frac
         { \displaystyle B \ prop_s \hspace{.3cm}C
         \ prop_s}
         {\displaystyle B\vee C \ prop_s }}\qquad
\displaystyle{ \frac
       {\displaystyle B
        \ prop_s \qquad C
        \ prop_s }
{B \rightarrow C \ prop_s }}\qquad
\displaystyle{ \frac
         { \displaystyle  B    \ prop_s\qquad  C
          \ prop_s }
         {\displaystyle  B\wedge C \ prop_s }}
\\[15pt]
\displaystyle{ \frac
         { \displaystyle   C(x)
         \ prop_s\ \ [x\in B]\qquad B\  set}
         {\displaystyle  \exists_{x\in B} C(x)\
prop_s}}\qquad
\displaystyle{ \frac
{\displaystyle C(x) \ prop_s\ [x\in B]\qquad B\ set  }
{\displaystyle \forall_{x\in B} C(x) \ prop_s }}
\qquad
\displaystyle{ \frac
{\displaystyle  A \hspace{.1cm} set \hspace{.3cm}  a\in A \hspace{.3cm} b\in A}
         {\displaystyle \orig{Id}(A, a, b) \ prop_s }}
\end{array}
$}
\\
\\

\noindent
And we add rules saying that a small proposition is a proposition and
 that a small proposition is  a set:
\\

\noindent
{\small

$\begin{array}{l}
    \mbox{\bf prop$_s$-into-prop)}\ \displaystyle{ \frac
       {\displaystyle A\ prop_s }
      { \displaystyle A\ prop}}
\end{array}\qquad
\qquad
\begin{array}{l}
    \mbox{\bf prop$_s$-into-set)}\
\displaystyle{ \frac
       {\displaystyle A\ prop_s }
      { \displaystyle A\ set}}
\end{array}\qquad
$}
\\
\\

\noindent
Then, we also have  the  collection of small propositions and
function collections from a set  toward it:
\\

\noindent
{\small
$ 
\begin{array}{l}
 \mbox{\bf Collection of small propositions}\\[5pt]
\mbox{F-Pr) } \ \displaystyle{\mathsf{prop_s}\  col }\qquad \mbox{I-Pr) }\displaystyle{ \frac
     { \displaystyle B \ prop_s}{ \displaystyle B \in \mathsf{prop_s}}}
\qquad \mbox{E-Pr) }\displaystyle{ \frac { \displaystyle B \in \mathsf{prop_s}}
         { \displaystyle B \ prop_s}}
\end{array}
$}
\\
\\

 \noindent
 {\small
       $\begin{array}{ll}
 \mbox{\bf Function collection to $\mathsf{prop_s}$ }
 \\[10pt]
 \mbox{ F-Fun)}\ 
 \displaystyle{\frac{ \displaystyle B\ set\  }
 {\displaystyle B\rightarrow \mathsf{prop_s} \hspace{.1cm} col }}
 &
 \mbox{ I-Fun)}\ \
 \displaystyle{\frac{ \displaystyle c(x)\in \mathsf{prop_s}\ [x\in B]
     \qquad B\ set }
 { \displaystyle \lambda x^{B}.c(x)\in  B\rightarrow \mathsf{prop_s}}}
       \\[15pt]
 \mbox{  E-Fun)}\ \
 \displaystyle{\frac{ \displaystyle
 b\in B \hspace{.3cm} f\in  B\rightarrow \mathsf{prop_s} }
 {\displaystyle \mathsf{Ap}(f, b)\in \mathsf{prop_s}  }}
 &
 \mbox{$\mathbf{\beta}$C-Fun)}\ \
 \displaystyle{\frac{ \displaystyle b\in B \hspace{.3cm} c(x)\in\mathsf{prop_s}
 \  [x\in B]  \qquad B\ set }
 {\displaystyle  \mathsf{Ap}(\lambda x^{B}.c(x),b)=c(b)\in  \mathsf{prop_s} }}
 \end{array}
 $}
 \\
 \\


\noindent
Equality rules include those saying that  type equality is an equivalence relation and  substitution of equal terms in a type:
\\

\noindent
{\small
$\begin{array}{l}
    \mbox{ref)}\,
\displaystyle{ \frac
       {\displaystyle A\ type }
      { \displaystyle A=A\ type}}
\end{array}
\qquad
\begin{array}{l}
    \mbox{sym)}\,
\displaystyle{ \frac
       {\displaystyle A=B\ type }
      { \displaystyle B=A\ type}}
\end{array}
\qquad
\begin{array}{l}
    \mbox{tra)}\,
\displaystyle{ \frac
       {\displaystyle A=B\ type \quad \ B=C\ type }
      { \displaystyle A=C\ type}}
\end{array}
\\[15pt]
\begin{array}{l}
    \mbox{subT)}\,\displaystyle{ \frac
       {\displaystyle\begin{array}{l}
C(x_1,\dots,x_n)\ type\ 
 [\, x_1\in A_1,\,  \dots,\,  x_n\in A_n(x_1,\dots,x_{n-1})\, ]   \\[2pt]
a_1=b_1\in A_1\ \dots \ a_n=b_n\in A_n(a_1,\dots,a_{n-1})\end{array}}
         {\displaystyle 
 C(a_1,\dots,a_{n})=C( b_1,\dots, b_n) \ type}}
   \end{array} 
$}
\\
\\

\noindent
where {\small $type \in \{ col, set,prop, prop_s\, \}$} with 
the same choice both in the premise and in the conclusion.
\\

\noindent
For terms into sets we add the following equality rules:\\
\\

\noindent
{\small
$\begin{array}{l}
    \mbox{ref)}\,
\displaystyle{ \frac
       {\displaystyle \ a\in A\ }
      { \displaystyle\ a=a\in A\ }}
\end{array}
\qquad
\begin{array}{l}
    \mbox{sym)}\,
\displaystyle{ \frac
       {\displaystyle\  a=b\in A\ }
      { \displaystyle\  b=a\in A\ }}
\end{array}
\qquad
\begin{array}{l}
    \mbox{tra)}\,
\displaystyle{ \frac
       {\displaystyle \ a=b\in A \qquad b=c\in A\ }
      { \displaystyle\  a=c\in A\ }}
\end{array}
$
\\
\\

\noindent
$\begin{array}{l}
      \mbox{sub)} \ \
\displaystyle{ \frac
         { \displaystyle 
\begin{array}{l}
 c(x_1,\dots, x_n)\in C(x_1,\dots,x_n)\ \
 [\, x_1\in A_1,\,  \dots,\,  x_n\in A_n(x_1,\dots,x_{n-1})\, ]   \\[2pt]
a_1=b_1\in A_1\ \dots \ a_n=b_n\in A_n(a_1,\dots,a_{n-1})
\end{array}}
         {\displaystyle c(a_1,\dots,a_n)=c(b_1,\dots, b_n)\in
 C(a_1,\dots,a_{n})  }}
      \\[15pt]
\mbox{conv)} \ \
\displaystyle{\frac{\displaystyle a\in A\ \qquad A=B\ type }
{\displaystyle a\in B}}
\qquad   
\mbox{conv-eq)} \ \
\displaystyle{\frac{\displaystyle a=b\in A\ \qquad A=B\ type }
{\displaystyle a=b\in B\ }}
\end{array}
$}
\\
\\

We call {\small \mtts\ } the fragment of \mtt\  consisting only of  judgements forming
sets, small propositions with their elements and corresponding equalities.
\\

Finally, the calculus \mttdp\ is obtained by extending \mtt\ with 
generic dependent collections:
\\
\\

{\small \noindent
      $\begin{array}{ll}
\mbox{\bf Dependent Product collection}
\\[10pt]
\mbox{F-$\Pi_c$)}\ \
\displaystyle{\frac{ \displaystyle  C(x)
        \hspace{.1cm} col\ [x\in B]  }
{\displaystyle \Pi_{x\in B} C(x)\hspace{.1cm} col }}
&
\mbox{I-$\Pi_c$)}\ \
\displaystyle{\frac{ \displaystyle c(x)\in  C(x)\ [x\in B]  }
{ \displaystyle \lambda x^{B}.c(x)\in \Pi_{x\in B} C(x)}}
      \\[15pt]
\mbox{E-$\Pi_c$)}\ \
\displaystyle{\frac{ \displaystyle
b\in B \hspace{.3cm} f\in \Pi_{x\in B} C(x) }
{\displaystyle \mathsf{Ap}(f,b)\in C(b) }}
&
\mbox{$\mathbf{\beta}$C-$\Pi_c$)}\ \
\displaystyle{\frac{ \displaystyle b\in B \hspace{.3cm} c(x)\in  C(x)\ [x\in B]}
{\displaystyle \mathsf{Ap}(\lambda x^{B}.c(x),b)=c(b)\in   C(b) }}
      \end{array}
$}
\\
\\

Note that in \mttdp\ function collections toward $\mathsf{prop_s}$ are
clearly a special instance of dependent product collections.

\section{Appendix:  The typed calculus
  \emtt\  }\label{emtt}
Here we present the calculus \emtt.
To build its types and terms we use
the same kinds of judgements used in \mtt,
namely
{\small  $$A \ type \ [\Gamma] \hspace{.5cm} A=B\ type\ [\Gamma] 
\hspace{.5cm} a \in A\ 
 [\Gamma] \hspace{.5cm} a=b \in A\ [\Gamma] $$}

\noindent
where types include collections, sets, propositions and small propositions:
namely 
{\small $$type \in \{ col, set,prop,prop_s\, \}$$}

\noindent
Contexts are generated by the same context rules of \mtt.
Also here note that the piece of context common to all judgements
involved in a rule
is omitted and that
 typed variables appearing in a context
are meant to be added to the implicit
context as the last one.

\noindent
Among types there are the following embeddings: sets are collections and propositions are collections
\\
 
\noindent
{\small
$\begin{array}{l}
\mbox{\bf set-into-col) }\ \ \displaystyle{ \frac
       {\displaystyle\  A
        \hspace{.1cm} set\ }
{ \displaystyle\  A
        \hspace{.1cm} col\ }}
\end{array}
\qquad\qquad
\begin{array}{l}
    \mbox{\bf prop-into-col) }\ \ 
\displaystyle{ \frac
       {\displaystyle\  A\ prop\  }
      { \displaystyle\  A\ col  \ }}
\end{array}
$}
\\

\noindent
Collections are closed under strong indexed sums:
\\

\noindent
{\small
$\begin{array}{l}
      \mbox{ \bf Strong Indexed Sum } \\[10pt]
     \mbox{F-}\Sigma ) \ \
\displaystyle{ \frac{\displaystyle   C(x)
         \hspace{.1cm} \ col \ [x\in B]}
         {\displaystyle \Sigma_{x\in B} C(x)\hspace{.1cm} col }}
         \qquad
      \mbox{I-}\Sigma )\ \
\displaystyle{ \frac
         {\displaystyle b\in  B \hspace{.3cm} c\in C(b)
\qquad C(x)\ col\ [x\in B]}
         {\displaystyle \langle b,c \rangle\in  \Sigma_{x\in B} C(x)}}
\\[15pt]
\mbox{E-}\Sigma )\ \
\displaystyle{ \frac
         {\displaystyle \begin{array}{l}
M(z)\ col \ [ z\in \Sigma_{x\in B} C(x)]\\
d\in  \Sigma_{x\in B} C(x) \hspace{.3cm} m(x,y)\in M(\langle x,y \rangle)\ [x\in B,
         y\in C(x)]
\end{array}}
      {\displaystyle {\it  El}_{\Sigma}(d,m)\in  M(d)}}
      \\[15pt]
\mbox{C-}\Sigma ) \ \
\displaystyle{ \frac
         {\displaystyle \begin{array}{l}
M(z)\ col \ [ z\in \Sigma_{x\in B} C(x)]\\
b\in B\ \ \  c\in C(b)\hspace{.3cm} m(x,y)\in M(\langle x,y \rangle)\ [x\in B,
         y\in C(x)]\end{array}}
      {\displaystyle {\bf \it El}_{\Sigma}(\langle b,c \rangle,m)=m(b,c)\in M(\langle b,c \rangle)} }
      \end{array}
$
\\
\\}

\noindent
Sets are generated  as follows:
\\
\\

\noindent
{\small
$\begin{array}{l}
      \mbox{\bf Empty set} \\
      \mbox{\small F-Em)} \ \mathsf{N_0} \hspace{.1cm} set \qquad
\mbox{\small E-Em)}\ 
\displaystyle{ \frac
         {\displaystyle a\in  \mathsf{N_0} \hspace{.3cm}
 A(x)\hspace{.1cm} \ col\
[x\in \mathsf{N_0}] }
         {\displaystyle \orig{emp_{o}}(a)\in A(a)}}
      \end{array}$
\\

\noindent
$\begin{array}{l}
\mbox {\bf Singleton set}\\
 \mbox{\small S)}\ \orig{\mathsf{N_1}} \hspace{.1cm} set
 \qquad
 \mbox{\small I-S)}\ 
 \orig{\star} \in\orig{\mathsf{N_1}}
 \qquad 
 \mbox{\small E-S)}\
\displaystyle{  \frac
    {t\in \orig{\mathsf{N_1}}\quad  M(z)\ col\ 
[z\in \mathsf{\mathsf{N_1}}] \quad c\in M(\star)}
    {{ \it El}_{ \orig{\mathsf{N_1}} }(t,c)\in  M(t)}}
\qquad 
 \mbox{\small C-S)}\
\displaystyle{  \frac
    { M(z)\ col\ 
[z\in \mathsf{\mathsf{N_1}}] \quad c\in M(\star)}
    {{ \it El}_{ \orig{\mathsf{N_1}} }(\star,c)=c\in  M(\star)}}
 \end{array}$
\\
\\

\noindent
$\begin{array}{l}
\mbox {\bf Strong Indexed Sum set} \\[10pt]
 \mbox{\small F-}\Sigma_s ) \ 
\displaystyle{  \frac
    { \displaystyle   C(x)
    \hspace{.1cm} set \ [x\in B]\qquad B\ set}
    {\displaystyle \Sigma_{x\in B} C(x)\hspace{.1cm} set }}
  \end{array}$
\\
\\

\noindent
      $\begin{array}{l}
\mbox{\bf List set} \\[10pt]
      \mbox{F-list)}\
\displaystyle{ \frac
         { \displaystyle C \hspace{.1cm} set}
         {\displaystyle List(C) \hspace{.1cm} set }}
      \qquad
       \mbox{${\rm I}_{1}$-list)}\ \
\displaystyle{ \frac
         {\displaystyle\quad  List(C) \hspace{.1cm} set\ }
         {\epsilon \in List(C)}}
      \qquad
      \mbox{${\rm I}_{2}$-list)}\ \
\displaystyle{ \frac
         {\displaystyle s\in List(C) \hspace{.3cm} c\in  C}
         {\displaystyle \orig{cons}(s,c)\in List(C)}}
      \end{array}$
\\
\\

\noindent
$\begin{array}{l}
\mbox{E-list)}\ \
\displaystyle{ \frac
         {\displaystyle \begin{array}{l}
L(z)\ col\  [z\in List(C)]
\hspace{.3cm}s\in List(C) \hspace{.3cm}\qquad  a\in L(\epsilon)\hspace{.3cm}\\
    l(x,y,z)\in
         L(\orig{cons}(x,y))  \ [x\in List(C),y\in C, z\in L(x)]
\end{array}}
      {\displaystyle {\bf \it El}_{List}(s,a, l)\in  L(s)}}
      \end{array}$
\\
\\

\noindent
$\begin{array}{l}
\mbox{${\rm C}_{1}$-list)}\ \
\displaystyle{  \frac
         {\displaystyle\begin{array}{l}
L(z)\ col\  [z\in List(C)] \hspace{.3cm}\qquad
 a\in L(\epsilon)\hspace{.3cm}\\
     l(x,y,z)\in
         L(\orig{cons}(x,y))  \ [x\in List(C),y\in C, z\in L(x)]
\end{array}}
      {\displaystyle {\bf \it El}_{List}( \epsilon, a,l)=a\in  L(\epsilon)}}
      \\[15pt]
\mbox{ ${\rm C}_{2}$-list)}\ \
\displaystyle{ \frac
         {\displaystyle\begin{array}{l}
L(z)\ col\  [z\in List(C)]
\hspace{.3cm} s\in List(C) \hspace{.3cm}c\in C \hspace{.3cm} a\in
L(\epsilon)\hspace{.3cm} \\
l(x,y,z)\in
         L(\orig{cons}(x,y))  \ [x\in List(C),y\in C, z\in L(x)]
\end{array}}
      {\displaystyle {\bf \it El}_{List}(\orig{cons}(s,c),a,l)=l(s,
c,{\bf \it El}_{List}(s,a,
      l))\in  L(\orig{cons}(s,c))}}
\end{array}$
\\
\\

\noindent
      $\begin{array}{l}
      \mbox{\bf Disjoint Sum set } \\[10pt]
      \mbox{F-+ )}\ \
\displaystyle{  \frac
         { \displaystyle B \hspace{.1cm} set \hspace{.3cm}C
         \hspace{.1cm} set}
         {\displaystyle B+ C \hspace{.1cm} set }}
      \qquad
      \mbox{${\rm I}_{1}$-}+ ) \ \
\displaystyle{  \frac
         {\displaystyle b\in  B \qquad B\ set\qquad C\ set }
         {\displaystyle\orig{inl}(b)\in B+ C}}
      \qquad
\mbox{${\rm I}_{2}$-}+ ) \ \
\displaystyle{  \frac
         {\displaystyle c\in  C\qquad B\ set\qquad C\ set }
         {\displaystyle  \orig{inr}(c)\in B+ C}}
      \\[15pt]
\mbox{E-} + ) \ \
\displaystyle{  \frac
         {\displaystyle\begin{array}{l}
A(z)\ col\  [z\in B+C]\\
      w\in B+ C \hspace{.3cm}
      a_{B}(x)\in A(\orig{inl}(x))\ [x\in B]
      \hspace{.3cm}\
      a_{C}(y)\in A(\orig{inr}(y))\ [y\in C]
\end{array}}
         {\displaystyle { \it El}_{+}(w,a_{B},a_{C})\in A(w)}}
      \\[15pt]
      \mbox{${\rm C}_{1}$-}+ ) \ \
\displaystyle{ \frac
       {\displaystyle \begin{array}{l}
A(z)\ col\  [z\in B+C]\qquad\\
b\in B \hspace{.3cm}  a_{B}(x)\in A(\orig{inl}(x))\ [x\in B]
      \hspace{.3cm}\
      a_{C}(y)\in A(\orig{inr}(y))\ [y\in C]\end{array}}
      {\displaystyle { \it El}_{+}(\orig{inl}(b),a_{B},a_{C})=a_{B}(b)
\in A(\orig{inl}(c))}}
\\[15pt]
      \mbox{${\rm C}_{2}$-}+ ) \ \
\displaystyle{ \frac
       {\displaystyle\begin{array}{l}
A(z)\ col\  [z\in B+C]\qquad \\
      c\in C\hspace{.3cm} a_{B}(x)\in A(\orig{inl}(x))\ [x\in B]
      \hspace{.3cm}\
      a_{C}(y)\in A(\orig{inr}(y))\ [y\in C]\end{array}}
      {\displaystyle {\bf \it El}_{+}(\orig{inr}(c),
a_{B},a_{C})=a_{C}(c)\in A(\orig{inr}(c))}}
      \end{array}$
\\
\\

\noindent
   $\begin{array}{l}
\mbox{\bf Dependent Product set}
\\[10pt]
\mbox{\small F-$\Pi$)} \
\displaystyle{\frac{ \displaystyle  C(x)
        \hspace{.1cm} set\ [x\in B]\qquad B\ set  }
{\displaystyle \Pi_{x\in B} C(x)\hspace{.1cm} set }}
\qquad
 \mbox{ I-$\Pi$)}\ \
 \displaystyle{\frac{ \displaystyle c(x)\in  C(x)\ [x\in B] \qquad  C(x)
        \hspace{.1cm} set\ [x\in B]\qquad B\ set }
 { \displaystyle \lambda x^{B}.c(x)\in \Pi_{x\in B} C(x)}}
       \\[15pt]
 \mbox{  E-$\Pi$)}\ \
 \displaystyle{\frac{ \displaystyle
 b\in B \hspace{.3cm} f\in \Pi_{x\in B} C(x) }
 {\displaystyle \mathsf{Ap}(f,b)\in C(b) }}
 \\[15pt]
 \mbox{$\mathbf{\beta}$C-$\Pi$)}\ \
 \displaystyle{\frac{ \displaystyle b\in B \hspace{.3cm} c(x)\in  C(x)\ [x\in B]\qquad  C(x)
        \hspace{.1cm} set\ [x\in B]\qquad B\ set}
 {\displaystyle \mathsf{Ap}(\lambda x^{B}.c(x),b)=c(b)\in   C(b) }}
\\[15pt]
 \mbox{$\mathbf{\eta}$C-$\Pi$}\ \ \displaystyle{\frac{ \displaystyle
 f\in \Pi_{x\in B} C(x) }
 {\lambda x^{B}.\mathsf{Ap}(f,x)=f\in \Pi_{x\in B} C(x)  }}
 (x\ not\ free\  in \ f)
      \end{array}
$
\\
\\

\noindent
 $\begin{array}{l}
\mbox{\bf Quotient set}
\\[10pt]
 \mbox{\small Q)} \ \
 \displaystyle{\frac
    { \displaystyle\begin{array}{ll}
 A\ set &
 R(x,y) \hspace{.1cm} prop_s\ [x\in A, y\in A] \\[5pt]  
& \begin{array}{l}
\mathsf{true}\in  R(x,x)\ [x\in A]
\\
 \mathsf{true}\in  R(y,x)\ [x\in A,\, y\in A,\, u\in R(x,y)]\\
\mathsf{true} \in  R(x,z)\ [x\in A,\, y\in A,\,  z\in A,\, \\
\qquad\qquad\qquad \qquad\qquad \quad\,  u\in 
  R(x,y),\,  v\in R(y,z)]\end{array}\end{array}}
    { \displaystyle A/R \hspace{.1cm} set }}
 \\[10pt]
  \mbox{\small I-Q)} \ 
\displaystyle{ \frac
    {\displaystyle a\in A  \hspace{.1cm}\qquad A/R \ set}
     {\displaystyle [a]\in A/R }}
    \qquad
   \mbox{ eq-}Q ) \ \  
 \displaystyle{\frac
     {\displaystyle a\in A \hspace{.3cm}b\in A \hspace{.3cm}
 \mathsf{true}\in R(a,b)\qquad
 A/R\ set}
    {\displaystyle [a]=[b]\in A/R }}
 \\[15pt]
 \mbox{\small  E-Q)} \ 
 \displaystyle{\frac
      {\displaystyle
 \begin{array}{l}
 L(z)\ col\ [z\in  A/R]\\[2pt]
  p\in A/R \hspace{.3cm} l(x)\in L([x])\ [x\in A]\hspace{.3cm} 
      l(x)=l(y)\in L([x])\ [x\in A, y\in A,
   d\in R(x,y)]\end{array}}
     {\displaystyle {\bf \it El}_{Q}(p,l)\in L(p)}}
 \\[12pt]
  \mbox{\small  C-Q)}  \ 
 \displaystyle{ \frac
      {\displaystyle \begin{array}{l}
 L(z)\ col\ [z\in  A/R]\\[2pt]
 a \in A \hspace{.3cm}
  l(x)\in L([x])\ [x\in A]\hspace{.3cm} l(x)=l(y)\in L([x])\ [x\in A, y\in A,
   d\in R(x,y)]
 \end{array}}
     {\displaystyle {\bf \it El}_{Q}([a],l)=l(a)\in L([a])}}
 \\[10pt]
 \mbox{\bf Effectiveness}\\[5pt]
\mbox{eff) }
 \displaystyle{  \frac
      {\displaystyle a\in A \hspace{.3cm} b\in A
  \hspace{.3cm} [a]=[b]\in  A/R\qquad A/R\ set}
     {\displaystyle \mathsf{true}\in R(a,b)}}
\end{array}
 $
\\
\\

}





\noindent
\emtt\ propositions are mono, namely
they are inhabited by at most a canonical proof-term:
\\

\noindent
{\small
$
\begin{array}{l}
\mbox{{\bf prop-mono}) } \,\displaystyle{ \frac{ \displaystyle A\ prop\ \qquad p\in A\  \quad q\in A\ }
{\displaystyle p=q\in A}} 
\end{array}
\qquad
\begin{array}{l}
\mbox{{\bf prop-true}) } \,\displaystyle{ \frac{ \displaystyle A\ prop\ \qquad p\in A\  }
{\displaystyle \mathsf{true}\in A}} 
\end{array}
$
}
\\

\noindent
Propositions are generated as follows:
\\

\noindent 
{\small $\begin{array}{lcc}
    \mbox{\bf Falsum }&\\
 \mbox{\small F-Fs)}\ \ \bot \hspace{.1cm} prop&\qquad 
\mbox{\small E-Fs)}\ 
\displaystyle{  \frac
    {\displaystyle \mathsf{true}\in  \bot \hspace{.3cm} A\hspace{.1cm} prop}
    { \displaystyle \mathsf{true}\in A}}
\end{array} 
$
\\
\\

\noindent
$
\begin{array}{l}
\mbox{\bf Extensional Propositional Equality}
\\[10pt]
 \mbox{\small F-Eq)}\ \
\displaystyle{  \frac
    {\displaystyle C \hspace{.1cm} \ col
 \hspace{.3cm}  c\in C \hspace{.3cm} d\in C}
    {\displaystyle \orig{Eq}(C, c, d)\hspace{.1cm} prop}}
 \qquad
 \mbox{\small I-Eq)}\ \
 \displaystyle{ \frac{\displaystyle  c \in C}
    {\displaystyle \mathsf{true}\in \orig{Eq}(C, c, c)}}
\qquad
 \mbox{\small E-Eq)}\ \
\displaystyle{  \frac
  {\displaystyle \mathsf{true}\in \orig{Eq}(C, c, d)}
 { \displaystyle c=d \in C }}
\end{array} $
\\
\\

\noindent
$
\begin{array}{l}
\mbox{\bf Implication}
\\[10pt]
\mbox{\small F-Im) }\ 
\displaystyle{ \frac
       {\displaystyle B
        \ prop\quad C
        \ prop }
{B \rightarrow C \ prop}}
\qquad
\mbox{\small I-Im) }\ 
\displaystyle{ \frac
{\displaystyle
 \mathsf{true} \in  C\ [x\in B]\qquad B \ prop\qquad  C\ prop  }
{\displaystyle \mathsf{true} \in B \rightarrow C}}
        \\[10pt]
\mbox{\small  E-Im) } \
\displaystyle{ \frac
{\displaystyle  \mathsf{true} \in B \quad  \mathsf{true} \in B\rightarrow C }
{\displaystyle \mathsf{true} \in C }}
      \end{array}$
\\
\\

\noindent
$\begin{array}{l}
\mbox{\bf Conjunction}
\\[10pt]
     \mbox{\small F-}\wedge)  \
\displaystyle{ \frac
         { \displaystyle  B   \ prop\quad  C
         \ prop}
         {\displaystyle  B\wedge C\hspace{.1cm} prop}}
         \quad
      \mbox{I-}\wedge )\ \
\displaystyle{ \frac
         {\displaystyle  \mathsf{true}\in  B \hspace{.3cm} \mathsf{true}\in  C
\qquad B\ prop\qquad C\ prop}
         {\displaystyle \mathsf{true}\in B\wedge C  }}
\\[10pt]
\mbox{\small $\mathbf{ \textrm{ E}_{1}}$-$\wedge$) }\ \
\displaystyle{ \frac
       {\displaystyle \mathsf{true}\in  B \wedge C}
{\displaystyle \mathsf{true}\in  B}}
\quad
\mbox{\small $\mathbf{ \textrm{ E}_{2}}$-$\wedge$) }\ \
\displaystyle{ \frac
       {\displaystyle \mathsf{true}\in B \wedge C}
{\displaystyle \mathsf{true}\in   C}}
\end{array}$
\\
\\

\noindent
      $\begin{array}{l}
\mbox{\bf Disjunction}
      \\[10pt]
     \mbox{\small F-}\vee )\ 
\displaystyle{ \frac
         { \displaystyle B \hspace{.1cm} prop\hspace{.3cm}C
         \hspace{.1cm} prop}
         {\displaystyle B\vee C \hspace{.1cm} prop}}
      \quad
      \mbox{\small ${\rm I}_{1}$-}\vee) \ 
\displaystyle{ \frac
         {\displaystyle\mathsf{true}\in  B\qquad B\ prop\qquad C\ prop}
         {\displaystyle  \mathsf{true}\in B\vee C}}
      \quad
\mbox{\small ${\rm I}_{2}$-}\vee ) \
\displaystyle{ \frac
         {\displaystyle   \mathsf{true}\in  C \qquad B\ prop\qquad C\ prop}
         {\displaystyle \mathsf{true}\in B\vee C}}
\end{array}$
\\[10pt]

\noindent
$\begin{array}{l}
\mbox{E-}\vee ) \ 
\displaystyle{ \frac
         {\displaystyle \begin{array}{l}
A\  prop\qquad      \mathsf{true}\in B\vee C \hspace{.3cm}
     \mathsf{true}\in A\ [x\in B]
      \hspace{.3cm}\
      \mathsf{true}\in A \ [y\in C]\end{array}}
         {\displaystyle  \mathsf{true}\in A}}
      \end{array}$
\\
\\



\noindent
   $\begin{array}{l}
      \mbox{\bf Existential quantification}
      \\[10pt]
     \mbox{\small F-}\exists ) \
\displaystyle{ \frac
         { \displaystyle   C(x)
         \hspace{.1cm} prop\ \ [x\in B]}
         {\displaystyle  \exists_{x\in B} C(x)\hspace{.1cm} prop}}
         \qquad
      \mbox{\small I-}\exists )\ \
\displaystyle{ \frac
         {\displaystyle   b\in  B \hspace{.3cm}
  \mathsf{true}\in  C(b)\qquad C(x)
         \hspace{.1cm} prop\ \ [x\in B]}
         {\displaystyle \mathsf{true} \in  \exists_{x\in B} C(x)}}
      \\[10pt]
\mbox{\small E-}\exists )\ 
\displaystyle{  \frac
         {\displaystyle \begin{array}{l}
M\ prop\qquad \mathsf{true}\in  \exists_{x\in B} C(x) \hspace{.3cm}
\mathsf{true}\in M\ [x\in B,
         y\in C(x)]
\end{array}}
      {\displaystyle \mathsf{true}\in  M}}
      \end{array}$
\\
\\

\noindent
      $\begin{array}{ll}
\mbox{\bf Universal quantification}
\\[10pt]
\mbox{\small F-$\forall$) } \
\displaystyle{ \frac
{\displaystyle C(x)\hspace{.1cm} prop\ [x\in B]  }
{\displaystyle \forall_{x\in B} C(x)\hspace{.1cm} prop}}
&
\mbox{\small I-$\forall$) } \
\displaystyle{ \frac
{\displaystyle \mathsf{true}\in   C(x)\ [x\in B] \qquad C(x)
         \hspace{.1cm} prop\ \ [x\in B] }
{\mathsf{true}\in  \forall_{x\in B} C(x)}}
        \\[10pt]
\mbox{\small  E-$\forall$) }\ 
\displaystyle{ \frac
{\displaystyle b\in B \hspace{.3cm} \mathsf{true}\in \forall_{x\in B} C(x) }
{\displaystyle \mathsf{true}\in C(b) }}
      \end{array}$
}
\\
\\

\noindent
As in \mtt, small propositions are generated as follows:
\\

\noindent
{\small
$
\begin{array}{l}
\bot \ prop_s\qquad
\displaystyle{ \frac
         { \displaystyle B \ prop_s \hspace{.3cm}C
         \ prop_s}
         {\displaystyle B\vee C \ prop_s}}\qquad
\displaystyle{ \frac
       {\displaystyle B
        \ prop_s\qquad C
        \ prop_s}
{B \rightarrow C \ prop_s}}\qquad
\displaystyle{ \frac
         { \displaystyle  B \   prop_s \qquad  C
          \ prop_s}
         {\displaystyle  B\wedge C \ prop_s}}
\\[15pt]
\displaystyle{ \frac
         { \displaystyle   C(x)
          \  prop_s \ [x\in B]\qquad B\  set}
         {\displaystyle  \exists_{x\in B} C(x)\in
prop_s}}\qquad
\displaystyle{ \frac
{\displaystyle C(x) \ prop_s\ [x\in B]\qquad B\ set  }
{\displaystyle \forall_{x\in B} C(x) \ prop_s}}
\qquad
\displaystyle{ \frac
{\displaystyle  A \hspace{.1cm} set \hspace{.3cm}  a\in A \hspace{.3cm} b\in A}
         {\displaystyle \orig{Eq}(A, a, b) \ prop_s}}
\end{array}
$}
\\
\\

\noindent
And we add rules saying that a small proposition is a proposition and
 that a small proposition is  a set:
\\

\noindent
{\small

$\begin{array}{l}
    \mbox{\bf prop$_s$-into-prop)}\ \displaystyle{ \frac
       {\displaystyle A\ prop_s }
      { \displaystyle A\ prop}}
\end{array}\qquad
\qquad
\begin{array}{l}
    \mbox{\bf prop$_s$-into-set)}\
\displaystyle{ \frac
       {\displaystyle A\ prop_s }
      { \displaystyle A\ set}}
\end{array}\qquad
$}
\\
\\

\noindent
Contrary to \mtt, in \emtt\ we do not have the intensional collection
of small propositions but 
 the  quotient of the collection of small propositions under equiprovability
representing the power collection of the singleton:
\\

\noindent
{\small
$
\begin{array}{l}
 \mbox{\bf Power collection of the singleton}\\[5pt]
 \mbox{F-{\cal P})}\ \ \displaystyle{{\cal P}(1)\  col }\qquad \mbox{I-{\cal P})}\ \displaystyle{ \frac
     { \displaystyle B \ prop_s}{ \displaystyle [B] \in {\cal P}(1)}}
\qquad \mbox{eq-{\cal P})}\ \displaystyle{ \frac { \displaystyle\mathsf{true}\in  B\leftrightarrow C\ }
         { \displaystyle [B]=[C]\in {\cal P}(1)  }}\qquad
 \mbox{eff-{\cal P})}\ \displaystyle{ \frac 
         { \displaystyle [B]=[C]\in {\cal P}(1) } { \displaystyle\mathsf{true}\in  B\leftrightarrow C\ }}\qquad\\[15pt]
 \displaystyle{ \frac
     { \displaystyle  U\in {\cal P}(1)\qquad  V\in {\cal P}(1)}{ \displaystyle 
\mathsf{Eq}( {\cal P}(1),\, U, \, V\, ) \  prop_s}}
\qquad\qquad
\mbox{$\eta$-{\cal P})}\ \displaystyle{ \frac 
         { \displaystyle U\in {\cal P}(1) } 
{ \displaystyle U= [\, 
\mathsf{Eq}( {\cal P}(1),\, U, \, [\mathsf{tt}]\, )\, ]\ }}$$
\end{array}
$}
\\

\noindent
 where {\small $\mathsf{tt}\, \equiv\, \bot \rightarrow \bot$} represents
the truth constant.
\\
\\

\noindent
Then, we have also function collections from a set toward ${\cal P}(1)$:
\\

 \noindent
 {\small
       $\begin{array}{l}
 \mbox{\bf Function collection to ${\cal P}(1)$}
 \\[10pt]
\begin{array}{ll}
 \mbox{ F-Fc)}\ \
 \displaystyle{\frac{ \displaystyle B\ set\  }
 {\displaystyle B\rightarrow {\cal P}(1) \hspace{.1cm} col }}
 \qquad 
&
 \mbox{ I-Fc)}\ \
 \displaystyle{\frac{ \displaystyle c(x)\in {\cal P}(1)\ [x\in
       B]\qquad B\ set  }
 { \displaystyle \lambda x^{B}.c(x)\in  B\rightarrow  {\cal P}(1)}}
       \\[15pt]
 \mbox{  E-Fc)}\ \
 \displaystyle{\frac{ \displaystyle
 b\in B \hspace{.3cm} f\in  B\rightarrow {\cal P}(1)\qquad }
 {\displaystyle  \mathsf{Ap}(f,b)\in {\cal P}(1) }}
 \qquad &
 \mbox{$\mathbf{\beta}$C-Fc)}\ \
 \displaystyle{\frac{ \displaystyle b\in B \hspace{.3cm} c(x)\in {\cal P}(1)
 \  [x\in B] \qquad B\ set  }
 {\displaystyle \mathsf{Ap}( \lambda x^{B}.c(x),b)=c(b)\in  {\cal P}(1) }}
 \end{array} 
\\
\\
 \quad\mbox{$\mathbf{\eta}$C-Fc)}\ \ \displaystyle{\frac{ \displaystyle
 f\in B\rightarrow   {\cal P}(1)   }
 {\lambda x^{B}.\mathsf{Ap}(f,x)=f\in B\rightarrow    {\cal P}(1)  }}
 ( x\ not\ free\  in \ f)
 \end{array}
 $}
 \\
 \\



 \noindent
 Then, as in \mtt\  we add the embedding rules of sets into collections {\bf set-into-col}, of propositions into collections {\bf prop-into-col},
 of small propositions into sets  {\bf prop$_s$-into-set} and of small propositions into propositions {\bf prop$_s$-into-prop}.

\noindent
Moreover, we also add the equality rules {\small ref), sym), tra)} both for types and for
terms saying that type and term equalities
are equivalence relations, and the rules {\small conv), conv-eq)}.

\noindent
Contrary to \mtt, we add all the  equality rules about collections and sets
saying that their constructors preserve  type equality as follows:
\\
\\

{\small
\noindent
$\begin{array}{l}
      \mbox{ \bf Strong Indexed Sum-eq} \\[10pt]
      \mbox{eq-}\Sigma ) \ \
\displaystyle{ \frac{\displaystyle   C(x)=D(x)
         \hspace{.1cm} col \ \ [x\in B]\qquad B=E\ col}
         {\displaystyle \Sigma_{x\in B} C(x)= \Sigma_{x\in E} D(x)
\hspace{.1cm} col }}
\end{array}
 \qquad
 \begin{array}{l}
 \mbox{\bf Function collection-eq}
 \\[10pt]
 \mbox{eq-Fc)}\ \
 \displaystyle{\frac{ \displaystyle  B=E\ set }
 {\displaystyle B\rightarrow {\cal P}(1)=E\rightarrow {\cal P}(1)
\hspace{.1cm} col }}
\end{array}$}
\\
\\



\noindent
{\small
$\begin{array}{l}
\mbox{\bf Lists-eq} \\[10pt]
      \mbox{eq-list)}\
\displaystyle{ \frac
         { \displaystyle C=D \hspace{.1cm} set}
         {\displaystyle List(C)=List(D) \hspace{.1cm} set }}
\end{array}
\qquad
\begin{array}{l}
      \mbox{ \bf Strong Indexed Sum set-eq} \\[10pt]
      \mbox{eq-}\Sigma_s ) \ \
\displaystyle{ \frac{\displaystyle   C(x)=D(x)
         \hspace{.1cm} set \ \ [x\in B]\qquad B=E\ set}
         {\displaystyle \Sigma_{x\in B} C(x)= \Sigma_{x\in E} D(x)
\hspace{.1cm} set }}\end{array}
$
\\
\\    

\noindent
      $\begin{array}{l}
      \mbox{\bf Disjoint Sum-eq} \\[10pt]
      \mbox{eq-+)} \  \
\displaystyle{  \frac
         { \displaystyle B=D \hspace{.1cm} set \hspace{.3cm}C=E
         \hspace{.1cm} set}
         {\displaystyle B+ C =D+E\hspace{.1cm} set }}
\end{array}\qquad\begin{array}{l}
\mbox{\bf Dependent Product-eq}
\\[10pt]
\mbox{eq-$\Pi$)}\ \
\displaystyle{\frac{ \displaystyle  C(x)=D(x)
        \hspace{.1cm} set\ [x\in B] \qquad B=E\ set }
{\displaystyle \Pi_{x\in B} C(x)=\Pi_{x\in E} D(x)\hspace{.1cm} set }}
      \end{array}$
\\
\\

\noindent
$
\begin{array}{l}
\mbox{\bf Quotient set-eq}\\[10pt]
\mbox{\small eq-Q)} \ \
 \displaystyle{\frac
    { \displaystyle A=B\ set\qquad
 R(x,y)=S(x,y) \hspace{.1cm}\ prop_s\ [x\in A, y\in A] \qquad  
\mathsf{Equiv}(R)\qquad \mathsf{Equiv}(S)}
    { \displaystyle A/R=B/S \hspace{.1cm} \ set}}
\end{array}$}
\\
\\


\noindent
Then, \emtt\ includes the following equality rules about propositions:
\\
\\

\noindent
{\small
      $\begin{array}{l}
\mbox{\bf Disjunction-eq}
      \\[10pt]
      \mbox{eq-}\vee )\ \
\displaystyle{ \frac
         { \displaystyle B=D \hspace{.1cm} prop \hspace{.3cm}C=E
         \hspace{.1cm} prop}
         {\displaystyle B\vee C= D\vee E \hspace{.1cm} prop }}
\end{array}
\qquad \begin{array}{l}
\mbox{\bf Implication-eq}
\\[10pt]
\mbox{eq-$\rightarrow$)}\ \
\displaystyle{ \frac
       {\displaystyle B=D
        \hspace{.1cm} prop \qquad C=E
        \hspace{.1cm} prop  }
{B \rightarrow C= D\rightarrow E \hspace{.1cm} prop }}
\end{array}
$
\\
\\

\noindent
$
 \begin{array}{l}
\mbox{\bf Conjunction-eq}
\\[10pt]
      \mbox{eq-}\wedge) \ \
\displaystyle{ \frac
         { \displaystyle  B=D   \hspace{.1cm} prop\qquad  C=E
         \hspace{.1cm} prop }
         {\displaystyle  B\wedge C=D\wedge E \hspace{.1cm} prop }}
      \end{array}
\qquad
\begin{array}{l}
\mbox{\bf Propositional equality-eq}
      \\[15pt]
      \mbox{eq-Eq)}\ \
\displaystyle{ \frac
{\displaystyle  A=E \hspace{.1cm} col \hspace{.3cm}  a=e\in A \hspace{.3cm} 
b=c\in A}
         {\displaystyle \orig{Eq}(A, a, b)= \orig{Eq}(E, e, c)\hspace{.1cm} prop }}
      \end{array}$
\\
\\

\noindent
      $\begin{array}{l}
      \mbox{\bf Existential quantification-eq}
      \\[10pt]
      \mbox{eq-}\exists ) \ \
\displaystyle{ \frac
         { \displaystyle   C(x)=D(x)
         \hspace{.1cm} prop \ \ [x\in B]\qquad B=E\ col}
         {\displaystyle  \exists_{x\in B} C(x)=\exists_{x\in E} D(x)
\hspace{.1cm} prop }}
\end{array}     \qquad
\begin{array}{l}
\mbox{\bf Universal quantification-eq}
\\[10pt]
\mbox{eq-$\forall$)}\ \
\displaystyle{ \frac
{\displaystyle C(x)=D(x)\hspace{.1cm} prop\ [x\in B]\qquad B=E \hspace{.1cm}
 col  }
{\displaystyle \forall_{x\in B} C(x)=\forall_{x\in E} D(x)\hspace{.1cm} prop }}
\end{array}
$}
\\
\\
Analogously, we add
 \mbox{eq-$\vee$)}, \mbox{eq-$\rightarrow$)},   
\mbox{eq-$\wedge$)}, \mbox{eq-Eq)}, \mbox{eq-$\exists$)}, \mbox{eq-$\forall$)}
restricted to small propositions.

\noindent
Moreover,  equality of propositions is that of collections, that of small propositions coincides with that
 of  propositions
and that of sets:
\\

\noindent
{\small
$\begin{array}{l}
    \mbox{\bf prop-into-col eq})\ 
\displaystyle{ \frac
       {\displaystyle A=B\ prop }
      { \displaystyle A=B\ col}}\end{array}
\qquad
\\[15pt]
\begin{array}{l}
 \mbox{\bf prop$_{s}$-into-prop eq})\
\displaystyle{ \frac
       {\displaystyle A=B\ prop_s }
      { \displaystyle A=B\ prop}}\end{array}
\qquad \begin{array}{l}
    \mbox{\bf prop$_s$-into-set eq})\
\displaystyle{ \frac
       {\displaystyle A=B\ prop_s }
      { \displaystyle A=B\ set}}
\end{array}
$}
\\
\\

\noindent
Equality of sets is that of collections:
\\

\noindent
{\small $
\begin{array}{l}
    \mbox{\bf set-into-col eq})\
\displaystyle{ \frac
       {\displaystyle A=B\ set }
      { \displaystyle A=B\ col}}
\end{array}
$}
\\
\\

\noindent
Contrary to \mtt, also for terms we add  equality rules 
saying that all the constructors preserve equality as in \cite{PMTT}:
\\

\noindent
{\small
$\begin{array}{l}
      \mbox{I-eq }\Sigma )\ \
\displaystyle{ \frac
         {\displaystyle b=b'\in  B \hspace{.3cm} c=c'\in C(b)
\qquad C(x)\ col\ [x\in B]}
         {\displaystyle \langle b,c\rangle=\langle b',c'\rangle \in  \Sigma_{x\in B} C(x)}}
\\[15pt]
\mbox{E-eq }\Sigma )\ \
\displaystyle{ \frac
         {\displaystyle \begin{array}{l}
M(z)\ col \ [ z\in \Sigma_{x\in B} C(x)]\\
d=d'\in  \Sigma_{x\in B} C(x) \hspace{.3cm}
 m(x,y)=m'(x,y)\in M(\langle x, y \rangle)\ [x\in B,
         y\in C(x)]
\end{array}}
      {\displaystyle {\it  El}_{\Sigma}(d,m)={\it  El}_{\Sigma}(d',m') \in  M(d)}}
      \end{array}$
\\
\\

\noindent
      $\begin{array}{ll}
\mbox{E-eq Em)}\ \
\displaystyle{ \frac
         {\displaystyle a=a'\in  \mathsf{N_0} \hspace{.3cm} A(x)\hspace{.1cm} col \
[x\in \mathsf{N_0}] } {\displaystyle \orig{emp_{o}}(a)= \orig{emp_{o}}(a')\in A(a)}}
&
\mbox{E-eq S)}\ \displaystyle{  \frac
    {t=t'\in \orig{\mathsf{N_1}}\quad  M(z)\ col\ 
[z\in \mathsf{\mathsf{N_1}}] \quad c=c'\in M(\star)}
    {{\it El}_{ \orig{\mathsf{N_1}} }(t,c)={\it El}_{ \orig{\mathsf{N_1}} }(t',c')\in
 M(t)
}}
 \end{array}$
\\
\\

\noindent
      $\begin{array}{l}
      \mbox{${\rm I}_{2}$-eq list)}\ \
\displaystyle{ \frac
         {\displaystyle s=s'\in List(C) \hspace{.3cm} c=c'\in  C}
         {\displaystyle \orig{cons}(s,c)=\mathsf{cons}(s',c')\in List(C)}}\\[12pt]
\mbox{E-eq list)}\ \
\displaystyle{ \frac
         {\displaystyle \begin{array}{l}
L(z)\ col\  [z\in List(C)]
\hspace{.3cm}s=s'\in List(C) \hspace{.3cm} a=a'\in L(\epsilon)\hspace{.3cm}\\[1pt]
    l(x,y,z)=l'(x,y,z)\in
         L(\orig{cons}(x,y))  \ [x\in List(C),y\in C, z\in L(x)]
\end{array}}
      {\displaystyle {\it El}_{List}(s,a, l)={\it El}_{List}(s',a', l')\in  L(s)}}
      \end{array}$
\\
\\

\noindent
      $\begin{array}{l}
\mbox{\small I-eq Q)} \ 
\displaystyle{ \frac
     {\displaystyle a=a'\in A\qquad  \hspace{.1cm} A/R \ set}
    {\displaystyle [a]=[a']\in A/R }}
\\[15pt]
\mbox{\small E-eq Q)} \ 
\displaystyle{\frac
     {\displaystyle
\begin{array}{l}
L(z)\ col\ [z\in  A/R]\\[1pt]
 p=p'\in A/R \hspace{.3cm} l(x)=l'(x)\in L([x])\ [x\in A]\hspace{.3cm} 
     l(x)=l(y)\in L([x])\ [x\in A, y\in A,
  d\in R(x,y)]\end{array}}
    {\displaystyle  {\it El}_{Q}(p,l)={\it El}_{Q}(p',l')\in L(p)}}
\end{array}$
\\
\\

\noindent
      $\begin{array}{l}
\mbox{\small ${\rm I}_{1}$-}+ ) \ 
\displaystyle{  \frac
    {\displaystyle b=b'\in  B\qquad B\ set\qquad C\ set }
    {\displaystyle \orig{inl}(b)=\orig{inl}(b')\in B+ C}}
 \qquad
 \mbox{\small ${\rm I}_{2}$-}+ ) \ 
 \displaystyle{ \frac
    {\displaystyle c=c'\in  C \qquad B\ set\qquad C\ set}
    {\displaystyle  \orig{inr}(c)=\orig{inr}(c')\in B+ C}}
\\[12pt]

\mbox{\small E-} + ) \ 
 \displaystyle{ \frac
    {\displaystyle 
\begin{array}{l}
A(z)\  col\ [z\in  B+ C]\\[1pt]
d=d'\in B+ C \hspace{.3cm}
 a_{B}(x)=a_{B}'(x) \in A(\orig{inl}(x))\ [x\in B]
 \hspace{.3cm}\ 
 a_{C}(y)= a_{C}'(y) \in A(\orig{inr}(y))\ [y\in C]
\end{array}}
    {\displaystyle {\it El}_{+}(d,a_{B},a_{C})={\it El}_{+}(d',a_{B}',a_{C}')\in A(w)}}
\end{array}$
\\
\\

\noindent
      $\begin{array}{ll}
\mbox{I-eq $\Pi$)}\ \
\displaystyle{\frac{ \displaystyle c(x)=c'(x)\in  C(x)\ [x\in B]\qquad
C(x)\ set\ [x\in B]\quad B\ set  }
{ \displaystyle \lambda x^{B}.c(x)= \lambda x^{B}.c'(x)\in \Pi_{x\in B} C(x)}}
      &\ \ 
\mbox{E-eq $\Pi$)}\ \
\displaystyle{\frac{ \displaystyle
b=b'\in B \hspace{.3cm} f=f'\in \Pi_{x\in B} C(x) }
{\displaystyle \mathsf{Ap}(f,b)= \mathsf{Ap}(f',b') \in C(b) }}
  \end{array}
$
\\
\\

\noindent
      $\begin{array}{ll}
\mbox{I-eq Fc)}\ \
\displaystyle{\frac{ \displaystyle c(x)=c'(x)\in {\cal P}(1)\ [x\in B]
   \quad B\ set }
{ \displaystyle \lambda x^{B}.c(x)= \lambda x^{B}.c'(x)\in
  B\rightarrow {\cal P}(1)}}
      &\ \ 
\mbox{E-eq Fc)}\ \
\displaystyle{\frac{ \displaystyle
b=b'\in B \hspace{.3cm} f=f'\in B\rightarrow {\cal P}(1) }
{\displaystyle \mathsf{Ap}(f,b)= \mathsf{Ap}(f',b') \in {\cal P}(1) }}
  \end{array}
$
\\
\\
}

\noindent
Note that \mbox{I-eq $\Pi$)} is the so-called $\xi$-rule in \cite{modmar}.
\\





We call \emtts\  the fragment of \emtt\ consisting only of  judgements forming
sets, small propositions with their elements and corresponding equalities.
\\

Finally, the calculus \emttdp\ is obtained by extending \emtt\ with 
generic dependent collections and quotient collections:
\\
\\

{\small
\noindent
      $\begin{array}{l}
\mbox{\bf Dependent Product Collection}
\\[10pt]
\mbox{F-$\Pi_c$)}\ \
\displaystyle{\frac{ \displaystyle  C(x)
        \hspace{.1cm} col\ [x\in B]  }
{\displaystyle \Pi_{x\in B} C(x)\hspace{.1cm} col }}
\qquad
\mbox{I-$\Pi_c$)}\ \
\displaystyle{\frac{ \displaystyle c(x)\in  C(x)\ [x\in B]  }
{ \displaystyle \lambda x^{B}.c(x)\in \Pi_{x\in B} C(x)}}
      \qquad
\mbox{E-$\Pi_c$)}\ \
\displaystyle{\frac{ \displaystyle
b\in B \hspace{.3cm} f\in \Pi_{x\in B} C(x) }
{\displaystyle \mathsf{Ap}(f,b)\in C(b) }}
\\[15pt]
\mbox{$\mathbf{\beta}$C-$\Pi_c$)}\ \
\displaystyle{\frac{ \displaystyle b\in B \hspace{.3cm} c(x)\in  C(x)\ [x\in B]}
{\displaystyle \mathsf{Ap}(\lambda x^{B}.c(x),b)=c(b)\in   C(b) }}
 \qquad
 \mbox{$\mathbf{\eta}$C-$\Pi_c$)}\ \ \displaystyle{\frac{ \displaystyle
 f\in \Pi_{x\in B} C(x) }
 {\lambda x^{B}.\mathsf{Ap}(f,x)=f\in \Pi_{x\in B} C(x)  }}
 (x\ not\ free\  in \ f)
      \end{array}
$
\\
\\

\noindent
$
 \begin{array}{l}
\mbox{\bf Quotient collection}
\\[10pt]
 \mbox{\small Q$_c$)} \ \
 \displaystyle{\frac
    { \displaystyle 
 R(x,y) \hspace{.1cm}\ prop\ [x\in A, y\in A] \quad  
 \begin{array}{l}
\mathsf{true}\in  R(x,x)\ [x\in A]
\\
 \mathsf{true}\in  R(y,x)\ [x\in A,\, y\in A,\, u\in R(x,y)]\\
\mathsf{true} \in  R(x,z)\ [x\in A,\, y\in A,\,  z\in A,\, \\
\qquad\qquad\qquad \qquad\qquad \quad\,  u\in 
  R(x,y),\,  v\in R(y,z)]\end{array}}
    { \displaystyle A/R \hspace{.1cm} \ col}}
 \\[15pt]
 \mbox{\small I-Q$_c$)} \ 
\displaystyle{ \frac
     {\displaystyle a\in A \qquad A/R \ col}
    {\displaystyle [a]\in A/R }}
   \qquad
  \mbox{eq-Q$_c$)} \ \  
\displaystyle{\frac
     {\displaystyle a\in A \hspace{.3cm}b\in A \hspace{.3cm} \mathsf{true}\in R(a,b)\qquad A/R\ col}
    {\displaystyle [a]=[b]\in A/R }}
\\[15pt]
\mbox{\small E-Q$_c$) } \ 
\displaystyle{\frac
     {\displaystyle
\begin{array}{l}
L(z)\ col\ [z\in  A/R]\\[2pt]
 p\in A/R \hspace{.3cm} l(x)\in L([x])\ [x\in A]\hspace{.3cm} 
     l(x)=l(y)\in L([x])\ [x\in A, y\in A,
  d\in R(x,y)]\end{array}}
    {\displaystyle  {\it El}_{Q}(p,l)\in L(p)}}
\\[10pt]
 \mbox{\small C-Q$_c$)}  \ 
\displaystyle{ \frac
     {\displaystyle \begin{array}{l}
L(z)\ col\ [z\in  A/R]\\[2pt]
a \in A \hspace{.3cm}
 l(x)\in L([x])\ [x\in A]\hspace{.3cm} l(x)=l(y)\in L([x])\ [x\in A, y\in A,
  d\in R(x,y)]
\end{array}}
    {\displaystyle {\bf \it El}_{Q}([a],l)=l(a)\in L([a])}}
\\[10pt]
\mbox{\bf Effectiveness}\\[5pt]
\mbox{eff$_c$) }\displaystyle{  \frac
     {\displaystyle a\in A \hspace{.3cm} b\in A
 \hspace{.3cm} [a]=[b]\in  A/R\qquad A/R\  col}
    {\displaystyle \mathsf{true}\in R(a,b)}}
\end{array}
 $}
\\
\\

Then, we also add the corresponding equality rules about
dependent product collections 
and about quotient collections as {\small eq-$\Pi$)},
{\small I-eq $\Pi$)}, {\small E-eq $\Pi$)}
 {\small eq-Q)},
{\small I-eq Q)}, {\small E-eq Q)}.
\\





Note that in \emttdp\ function collections toward {\small ${\cal P}(1)$} are
clearly a special instance of dependent product collections.

\section{Appendix: Interpretation of \emtt\   into \mtt}
\label{emcq}
Here we define
the interpretation of  \emtt-type and term signatures as
\mtt-extensional
dependent types and  terms, respectively, by using  canonical isomorphisms.
This interprets  the so-called ``raw syntax'' in
\cite{tumscs}, namely the signatures of types and terms in \emtt,
 in a partial way. Indeed, as it is  well explained in   \cite{Streicher},
we can not define a total interpretation  by induction on the derivation of 
types and typed terms in \emtt, because
term equalities are involved  in the formation of types and typed terms,
and hence the interpretation would depend on the validity of equality
which should instead follow as a consequence of the chosen interpretation.

\noindent
\begin{definition}[Interpretation of \emtt\ into \mtt\ and of \emttdp\
into \mttdp]
\label{intext}\em

\noindent
In the following we define  supports and related equalities
of \mtt\  (\mttdp)-extensional dependent types
 interpreting \emtt\  (\emttdp)-type signatures,
and  \mtt\  (\mttdp)-extensional terms interpreting  \emtt\ (\emttdp)-term signatures, all 
by induction on their formation.
They are  both described under a context, since
free variables are assumed
to be typed.
We also  warn that,
 when we interpret a type or a term signatures depending
on more than one  term, we assumed  to have matched  the types of the already interpreted
terms via canonical isomorphisms.

\noindent
The {\bf empty context} is interpreted as follows:\\
{\small $\begin{array}{l}
 [ \ ]^I  \equiv\, [x: \mathsf{N}_1]   \mbox{ and }\
  [ \mathsf{N}_{1_{=}}]\ \equiv\ ( \mathsf{N}_1, =_{\mathsf{N}_1})\\
\mbox{with }\ z=_{\mathsf{N}_{1}}z'\, \equiv \mathsf{Id}(\mathsf{N}_1, z,z')\ \mbox{ for }
z,z'\in \mathsf{N}_1 
\end{array}$}\\

\noindent
The {\bf assumption of variable} is interpreted as follows:\\
{\small $ (\, x\in A \  ) [\Gamma^I, x\in A, \Delta^I]\, ) \equiv\,
x\in A^I\ [\Gamma^I, x\in A^I, \Delta^I]$}
\\

\noindent
Collection constructors are interpreted as follows:\\

\noindent
{\small
$\begin{array}{l}\mbox{\bf Power collection of the singleton}:\\ 
{\cal P}(1)^I\ col\  [\Gamma^I]\, \equiv\, 
\mathsf{prop_s}\ [\Gamma^I]\\[3pt]
\mbox{and }z=_{\mathsf{prop_s}^I}z'\, \equiv\,
 ( z \rightarrow z')
\wedge  ( z' \rightarrow z)$ for $z,z'\in \mathsf{prop_s}\end{array}$
\\
$\sigma_{\overline{x}}^{\overline{x}'}(w)\, \equiv\, w$ for
$\overline{x}, \overline{x'}\in \Gamma^I, w\in \mathsf{prop_s}$.\\
$([A])^I\, \equiv\, A^I$ for $A$ small proposition.
\\
$(\, \mathsf{true}\in  B\leftrightarrow C\  \ [\Gamma]\, )^I\, \equiv\,  r\in (B\leftrightarrow C)^I \  [\Gamma^I] $ provided that $r\in  (B\leftrightarrow C)^I \  [\Gamma^I] $  is derivable in  \mtt.\\

\noindent
$\begin{array}{l}
\mbox{\bf Strong Indexed Sum}:\\
  (\, \Sigma_{y\in B} C(y)\,)^I \ col\ [\Gamma^I]\, \equiv\, 
\Sigma_{y\in B^{I}} C^{I}(y)\ col \ [\Gamma^I] \\[3pt]
\mbox{and }
 z=_{\Sigma_{y\in B} C(y)^{I}}z'\, \equiv\quad \exists_{d\in \pi_1(z)=_{B^{I}}
\pi_1(z')}\ \
\sigma_{\pi_1(z)}^{\pi_1(z')}( \pi_2(z))=_{C^I(\pi_1(z'))}\pi_2(z')\quad
\mbox{ for }z,z'\in (\, \Sigma_{y\in B} C(y)\, )^{I}.
\end{array}$\\
$ (\, \langle b,d\rangle\, )^{I}\, \equiv\, \langle b^{\widetilde{I}}, d^{\widetilde{I}}\rangle $ and 
$  {\it  El}_{\Sigma}(d,m)^I \, \equiv\,
 {\it  El}_{\Sigma}(d^{\widetilde{I}},(w_1,w_2).m^{\widetilde{I}})$\\[3pt]
$\sigma_{\overline{x}}^{\overline{x'}}(w)\, \equiv\,
\mathsf{El}_{\Sigma}(w, (w_1,w_2). \langle\,  \sigma_{\overline{x}}^{ \overline{x'}}(w_1)\, ,\,
\sigma_{\overline{x},w_1}^{ \overline{x'},
 \sigma_{\overline{x}}^{ \overline{x'}}(w_1)}(w_2)\, \rangle\, )$ for
$\overline{x}, \overline{x'}\in \Gamma^I$ and $
w\in (\, \Sigma_{y\in B} C(y)\, )^{\widetilde{I}}(\overline{x})$.\\

\noindent
$\begin{array}{l}
\mbox{\bf Dependent Product collection}:\\
(\, \Pi_{y\in B} C(y)\, )^I\ col \ [\Gamma^I]\, \equiv\,\ 
 \Sigma_{h\in \Pi_{y\in B^I} C^{I}(y)}\quad\ 
  \forall_{y_1\in B^I}\ \forall_{y_2\in B^I}\ \ \forall_{d\in  y_1=_{B^I}y_2}\ \ 
\sigma_{\overline{x},y_1}^{\overline{x},y_2}\, (\, \mathsf{Ap}(h,y_1)\, )=_{
C^I(y_2)}\mathsf{Ap}( h,y_2) \\
\mbox{and }z=_{\Pi_{y\in B^I} C(y)^I}z'\, 
\equiv\, \forall_{y\in B^I}\ \   \mathsf{Ap}(\, \pi_1(z)\, ,y)=_{C^I(y)}
\mathsf{Ap}(\, \pi_1(z')\, ,y)\ 
\mbox{ for }z,z'\in (\, \Pi_{y\in B}\, C(y)\, )^I.\end{array}$
\\
$(\, \lambda y^{B}.c\, )^I\, \equiv\, \langle\,
\lambda y^{B^I}.c^{\widetilde{I}}\, , \, p\, \rangle\quad $ where $p\in
\forall_{y_1\in B^I}\ \forall_{y_2\in B^I}\  \forall_{d\in   y_1=_{B^I}y_2}\ \ 
\sigma_{\overline{x},y_1}^{\overline{x},y_2}(\, c^{\widetilde{I}}(y_1)\, )
=_{C^I(y_2)} c^{\widetilde{I}}(y_2)$\\
$(\, \mathsf{Ap}(f,b) \, )^I\, \equiv\, 
\mathsf{Ap}(\, \pi_1(f^{\widetilde{I}})\, , b^{\widetilde{I}})$\\
$\sigma_{\overline{x}}^{\overline{x'}}(w)\, \equiv\,\langle\, 
\lambda y'^{B^I(\overline{x'})}.\  \sigma_{\overline{x},
\sigma_{\overline{x'}}^{\overline{x}}(y')}^{ \overline{x'},y'}
(\, \mathsf{Ap}(\, \pi_1(w)\, ,\, \sigma_{\overline{x'}}^{\overline{x}}(y') \,  )\,)\, ,\,  p\, \rangle$ for
$\overline{x}, \overline{x'}\in \Gamma^I$ and $w\in (\, \Pi_{y\in B}\, C(y)\, )^{I}(\overline{x})$\\
where $p$ is the proof-term witnessing the preservation of equalities obtained from $\pi_2(w)$.\\

\noindent
$\begin{array}{l}
\mbox{\bf Function collection to  ${\cal P}(1)$}:\\
 (\, B\rightarrow {\cal P}(1)\  col\ [\Gamma]\, )^I\, \equiv\,
\Sigma_{h\in B^I\rightarrow \mathsf{prop_s}} \ 
\quad\ 
  \forall_{y_1\in B^I}\ \forall_{y_2\in B^I}\ \  
\  y_1=_{B^I}y_2\, \rightarrow \,
(\, \mathsf{Ap}(h,y_1)\, \leftrightarrow\,  \mathsf{Ap}(h,y_2)\,) \\
\mbox{and }z=_{{\cal P}}z'\, 
\equiv\, \forall_{y\in B^{I}}\ \ \mathsf{Ap}( \,\pi_1(z)\, ,y)\, 
\leftrightarrow  \, \mathsf{Ap}( \, \pi_1(z')\, ,y)\ 
\mbox{ for }z,z'\in (\, B\rightarrow {\cal P}(1)\, )^{I}\end{array}$\\
$(\, \lambda y^{B}.c\, )^I\, \equiv\, \langle\,
\lambda y^{B^I}. c^{\widetilde{I}}\, , \, p \, \rangle$ where $p\in
\forall_{y_1\in B^{I}}\ \forall_{y_2\in B^{I}}\ \  y_1=_{B^I}y_2\, \rightarrow\, (\, 
c^{\widetilde{I}}(y_1)\, \leftrightarrow \, c^{\widetilde{I}}(y_2)\, )$\\
$(\, \mathsf{Ap}(f,b) \, )^I\, \equiv\, \mathsf{Ap}(\,
\pi_1(f^{\widetilde{I}})\, , b^{\widetilde{I}})$\\
$\sigma_{\overline{x}}^{\overline{x'}}(w)\, \equiv\,\langle\, 
\lambda y'^{B^{I}(\overline{x'})}.\,  \sigma_{\overline{x},
\sigma_{\overline{x'}}^{\overline{x}}(y')}^{ \overline{x'},y'}
(\, \mathsf{Ap}(\, \pi_1(w)\, ,\,  \sigma_{\overline{x'}}^{\overline{x}}(y') \,  )\, ),\, p\, \rangle$ for
$\overline{x}, \overline{x'}\in \Gamma^I$ and
 $w\in (\, B\rightarrow {\cal P}(1)\, )^I(\overline{x})$\\
where $p$ is the proof-term witnessing the preservation of equalities obtained from $\pi_2(w)$.
 \\

\noindent
$\begin{array}{l}
\mbox{\bf Quotient collection}:\\
(\, A/R\ col\ [\Gamma] )^I\, \equiv\, A^I \ col\ [\Gamma^I] \\[3pt]
\mbox{and }z=_{A/R^I}z' \, \equiv\,  R^I(z,z')\mbox{ for } 
z,z'\in A^I.\end{array}$
\\
$(\, [a]\, )\, \equiv\, a^I$ and
 ${\it El}_{Q}(p,l)^I\, \equiv\, l^{\widetilde{I}}(p^{\widetilde{I}}) $\\
$\sigma_{\overline{x}}^{\overline{x'}}(w)$ is defined
as the substitution isomorphism of $A^{I}_=\, [\, \Gamma^I_=]$.\\
$(\, \mathsf{true}\in R(a, b)\ [\Gamma]\, )^I\, \equiv\,  r\in R(a, b)^I \  [\Gamma^I] $ provided that $r\in R(a, b)^I \ [\Gamma^I] $ is derivable in  \mtt.\\}

\noindent
Now we give the interpretation of sets:
\\

{\small
\noindent
$\mbox{\bf Empty set}:
\ (\, \mathsf{N_0}\ set\ [\Gamma] )^I \, \equiv\, \mathsf{N_0}\ set\ [\Gamma^I] \\
  \mbox{and }
z=_{\mathsf{N_0}^I}z'\, \equiv \mathsf{Id}(\mathsf{N_0}, z,z')\ \mbox{ for }
z,z'\in \mathsf{N_0}.$\\
$(\,  \orig{emp_{o}}(a)\, )^I\, \equiv\, \orig{emp_{o}}(a^{\widetilde{I}})$\\
 $\sigma_{\overline{x}}^{\overline{x'}}(w)\, \equiv\, w$
for $\overline{x}, \overline{x'}\in \Gamma^I$ and   $w\in \mathsf{N_0}$.\\

\noindent
$\mbox{\bf Singleton set}:
\ (\, \mathsf{N_1}\ set\   [\Gamma] )^I \, \equiv\, \mathsf{N_1}\ set\ [\Gamma^I] \\
\mbox{and } 
z=_{\mathsf{N_1}^I}z'\, \equiv \mathsf{Id}(\mathsf{N_1}, z,z')\ \mbox{ for }
z,z'\in \mathsf{N_1}.$\\
$(\, \star\, )^I\, \equiv\, \star\quad $ and
$\quad (\, { \it El}_{ \mathsf{N_1} }(t,c)\, )^I\, \equiv\,
{ \it El}_{ \mathsf{N_1} }(t^{\widetilde{I}},c^{\widetilde{I}})$\\
 $\sigma_{\overline{x}}^{\overline{x'}}(w)\, \equiv\, w$
for  $\overline{x}, \overline{x'}\in \Gamma^I$ and 
 $w\in \mathsf{N_1}$.\\

\noindent
$
\mbox{\bf List set}:
\ (\, List(C)\ set\  [\Gamma] )^I\, \equiv\, List(C^I)\ set\  [\Gamma^I]$\\
\mbox{and } $ z =_{List(C^I)}z'$ defined as in theorem~\ref{mainth1}.
\\
$(\, \epsilon\, )^I\, \equiv\, \epsilon\ $ and 
$(\, \orig{cons}(s,c)\, )^I\, \equiv\, \mathsf{cons}(s^{\widetilde{I}},c^{\widetilde{I}})$\\
$(\, { \it El}_{List}(s,a,l)\, )^I\, \equiv\, { \it El}_{List}(\,
s^{\widetilde{I}},\,
a^{\widetilde{I}},\, l^{\widetilde{I}}\,)$\\
 $\sigma_{\overline{x}}^{\overline{x'}}(w)\, \equiv\, 
{\it El}_{List}(\, w\, ,\, \epsilon ,\,  (y_1,y_2,z).\mathsf{cons}( \,
z\, ,\, 
\sigma_{\overline{x}}^{\overline{x'}}
(y_2)\, )\,)$
for $\overline{x}, \overline{x'}\in \Gamma^I$ and 
 $w\in (\,  List(C)\, )^{I}(\overline{x})$.\\

\noindent
$\mbox{\bf Disjoint Sum set}:
\ (\, B+C\ set\ [\Gamma] )^I\, \equiv\, B^I+C^I\ set\ [\Gamma^I]$\\
\mbox{and }$ z =_{B^I+C^I}z'$
is defined as in theorem~\ref{mainth1}.\\
$( \, \orig{inl}(b)\, )^I\, \equiv\, \orig{inl}(b^I)\ $ and
$\ ( \, \orig{inl}(c)\, )^I\, \equiv\, \orig{inl}(c^I)\ $ and $\
( \, { \it El}_{+}(d,a_{B},a_{C})\, )^I\, \equiv\, 
{ \it El}_{+}(\, d^{\widetilde{I}},\, a_{B}^{\widetilde{I}},\, a_{C}^{\widetilde{I}})\, )$
\\
 $\sigma_{\overline{x}}^{\overline{x'}}(w)\, \equiv\, 
{\it El}_{+}(w, (y_1).\sigma_{\overline{x}}^{\overline{x'}}(y_1)
\, ,\, (y_2).\sigma_{\overline{x}}^{\overline{x'}}
(y_2)\, )$
for  $\overline{x}, \overline{x'}\in \Gamma^I$ and 
 $w\in (\, B+C\, )^I(\overline{x})$}.\\

\noindent
Finally, $\mbox{\bf Strong Indexed Sum set}$, $\mbox{\bf Dependent Product set}$ and
  $\mbox{\bf Quotient set}$  constructors with their terms are
 interpreted
analogously to Strong Indexed Sums, Dependent Product collections and Quotient collections
with their terms, respectively.\\
\\



\noindent
Lastly,
\emtt-propositions are interpreted as \mtt-extensional propositions
whose support is similar,
except for extensional propositional equality, and whose
equality  is trivial, namely if  {\small $A\ prop\ [\Gamma]$} is
 a proposition then
{\small $z=_{A^I}z'\, \equiv\, \mathsf{tt}$} for all {\small $\overline{x}\in \Gamma^I,
z,z'\in A^I(\overline{x}) $}.
Therefore in the following we just specify the support of their
 interpretation.
 \\

\noindent
{\small  {\bf Falsum}: $\bot^I \  prop\ [\Gamma^I])\, 
\equiv\, \bot \  prop\ [\Gamma^I]$\\
$(\, \mathsf{true}\in A\, )^I
\, \equiv\, \mathsf{r_0}(p)\in A^I$ provided that
 $p\in \bot$ is derived in \mtt.\\
$\sigma_{\overline{x}}^{\overline{x'}}(w)\, \equiv\,
w$ for
 $\overline{x}, \overline{x'}\in \Gamma^I$ and
 $w\in \bot$.\\

\noindent
{\bf Extensional Propositional Equality}:
 $ \orig{Eq}(B, b_1, b_2)^I\ prop\
[\Gamma^I] \, \equiv\, b_1^{\widetilde{I}}=_{B^I} b_2^{\widetilde{I}}\ prop\ [\Gamma^I]$\\
$(\, \mathsf{true}\in \orig{Eq}(B, b, b)\ [\Gamma]\, )^I\,
\equiv\, \mathsf{rfl}(b^{\widetilde{I}})\in  b^{\widetilde{I}}=_{B^I} b^{\widetilde{I}}\ [\Gamma^I]$
provided that $b^{\widetilde{I}}\in B^I\ [\Gamma^I]$  is derived in \mtt.\\
where $\sigma_{\overline{x}}^{\overline{x'}}$ is that of $b_1^{\widetilde{I}}=_{B^I} b_2^{\widetilde{I}}\
 prop\ [\Gamma^I]$\\

\noindent
{\bf Implication}:
$( B \rightarrow C)^I\ prop\ [\Gamma^I] \, \equiv\,
 B^I  \rightarrow C^I\  prop\ [\Gamma^I] $\\
$(\, \mathsf{true} \in B \rightarrow C\ [\Gamma]\, )^I
\, \equiv\,  \lambda_\rightarrow x^{B}.c\in B^I \rightarrow C^I\
[\Gamma^I]$ provided that $c\in C^I\ [\Gamma^I, y\in B^I]$ is
derived in \mtt.\\
$(\, \mathsf{true} \in C \ [\Gamma]\, )^I
\, \equiv\, \mathsf{Ap}_\rightarrow (f,b)\in C^I\ [\Gamma^I]$ provided that
$f\in B^I\rightarrow C^I\ [\Gamma^I]$ and $  b\in B^I\ [\Gamma^I]$ are
derived in \mtt.\\
$\sigma_{\overline{x}}^{\overline{x'}}(w)\, \equiv\,
\lambda_\rightarrow y\in B^I(\overline{x'}).\,
 \sigma_{\overline{x}}^{\overline{x'}}
(\, \mathsf{Ap}_\rightarrow ( w,\,  \sigma_{\overline{x'}}^{\overline{x}}(y)\, )\, ) $ for
 $\overline{x}, \overline{x'}\in \Gamma^I$ and
 $w\in ( B\rightarrow C)^I(\overline{x})$.
\\

\noindent
{\bf Conjunction}
$( B \wedge C)^I\ prop\ [\Gamma^I] \, \equiv\,
 B^I  \wedge C^I\  prop\ [\Gamma^I] $\\
$(\, \mathsf{true}\in B\wedge C \ [\Gamma]\, )^I
\, \equiv\, \langle b,_{\wedge} c \rangle\in B^I\wedge C^I\
  [\Gamma^I]$ provided that
$b\in B^I\ [\Gamma^I]$ and $  c\in C^I\ [\Gamma^I]$ are
derived in \mtt.\\
$(\, \mathsf{true}\in  B\ [\Gamma]\, )^I
\, \equiv\,
\pi_{1}^{B^I}(d)\in  B^I \ [\Gamma^I]$  provided that
$d\in B^I\wedge C^I\ [\Gamma^I]$ is derived in \mtt.\\
$(\, \mathsf{true}\in  C\ [\Gamma]\, )^I
\, \equiv\,
\pi_{2}^{C^I}(d)\in  C^I \ [\Gamma^I]$  provided that
$d\in B^I\wedge C^I\ [\Gamma^I]$ is derived in \mtt.\\
$\sigma_{\overline{x}}^{\overline{x'}}(w)\, \equiv\,
\langle\sigma_{\overline{x}}^{\overline{x'}}(\, 
\pi_{1}^{B^I}(w)\, )\, 
,_{\wedge}\, \sigma_{\overline{x}}^{\overline{x'}}(\, \pi_{2}^{C^I}(w)\, )
 \rangle$ for
 $\overline{x}, \overline{x'}\in \Gamma^I$ and
 $w\in (B\wedge C)^I(\overline{x})$.
\\

\noindent
{\bf Disjunction}
$( B \vee C)^I\ prop\ [\Gamma^I] \, \equiv\,
 B^I  \vee C^I\  prop\ [\Gamma^I] $\\
$(\, \mathsf{true}\in B\vee C \ [\Gamma]\, )^I
\, \equiv\,  \mathsf{inl}_{\vee}(b) \in B^I\vee C^I\
  [\Gamma^I]$ provided that
$b\in B^I\ [\Gamma^I]$ is
derived in \mtt.\\
$(\, \mathsf{true}\in B\vee C \ [\Gamma]\, )^I
\, \equiv\,  \mathsf{inr}_{\vee}(c) \in B^I\vee C^I\
  [\Gamma^I]$ provided that
$c\in C^I\ [\Gamma^I]$ is
derived in \mtt.\\
$(\, \mathsf{true}\in A\ [\Gamma]\, )^I\, \equiv\,
{\it El}_{\vee}(d,a_{B},a_{C})\in A^I\ [\Gamma^I]$
provided that $d\in B^I\vee C^I\ [\Gamma^I]$,
$ a_{B}\in A^I\ [\Gamma^I, y\in B^I]$ and
$ a_{C}\in A^I\ [\Gamma^I, y\in C^I]$ are
derived in \mtt.\\
$\sigma_{\overline{x}}^{\overline{x'}}(w)\, \equiv\,
{\it El}_{\vee}(w,\, (y_1). \, \mathsf{inl}_{\vee}(\,
 \sigma_{\overline{x}}^{\overline{x'}}(y_1)\, )\, 
, \, (y_2). \mathsf{inr}_{\vee}(\,
 \sigma_{\overline{x}}^{\overline{x'}}(y_2)\,)\, )$ for
 $\overline{x}, \overline{x'}\in \Gamma^I$ and
 $w\in (B\vee C)^I(\overline{x})$.
\\

\noindent
 {\bf Existential quantifier}:
$(  \exists_{y\in B} C(y))^I \ prop\ [\Gamma^I] )\, \equiv\, 
\exists_{y\in B^I} C^I(y)\  prop\ [\Gamma^I]$\\
$(\, \mathsf{true} \in  \exists_{y\in B} C(y)\ [\Gamma]\, )^I
\, \equiv\, \langle b^{\widetilde{I}},_{\exists} c \rangle\in  \exists_{y\in B^I} C^I(y)\
[\Gamma^I]$ provided that
 $c\in C^I(b^{\widetilde{I}})\ [\Gamma^I]$
 is
derived in \mtt.\\
$(\, \mathsf{true}\in  M\ [\Gamma]\, )^I
\, \equiv\, { \it El}_{\exists}(d,m)\in  M^I\ [\Gamma^I]$
provided that $d\in \exists_{y\in B^I} C^I(y)\ [\Gamma^I]$
and $m\in M^I\ [\Gamma^I, y\in B^I, z\in C^I(y)]$ are 
derived in \mtt.\\
$\sigma_{\overline{x}}^{\overline{x'}}(w)\, \equiv\,
{\it El}_{\vee}(w,\, (y,z). \, \langle \,
 \sigma_{\overline{x}}^{\overline{x'}}(y)\, \, 
, \,
 \sigma_{\overline{x},y}^{\overline{x'}, 
\sigma_{\overline{x}}^{\overline{x'}(y)}}(z)\,\rangle\, ) $ for
 $\overline{x}, \overline{x'}\in \Gamma^I$ and
 $w\in (\,\exists_{y\in B}\,  C(y)\, )^I(\overline{x})$.
\\

 \noindent
{\bf Universal quantifier}:
 $(\,   \forall_{y\in B} C(y)^I\ prop\ [\Gamma^I]  )\, 
\equiv\, \forall_{y\in B^I} C^I(y)\  prop\ [\Gamma^I]) $\\
$(\, \mathsf{true}\in  \forall_{y\in B} C(y)\ [\Gamma]\, )^I
 \, \equiv\, \lambda_\forall y^{B^I}.c\in \forall_{y\in B^I} C^I(y)
$ provided that
$c\in C^I(y)\ [\Gamma^I, y\in B^I]$ is derived in \mtt.\\
$(\, \mathsf{true}\in C(b) \ [\Gamma]\, )^I
 \, \equiv\,\mathsf{Ap}_\forall( f,b^{\widetilde{I}})
\in C^I(b^{\widetilde{I}})\ [\Gamma^I]
$ provided that
$f\in \forall_{y\in B^I} C^I(y)\ [\Gamma^I]$ is
derived in \mtt.\\
$\sigma_{\overline{x}}^{\overline{x'}}(w)\, \equiv\,
\lambda_\forall y^{B^I(\overline{x'})}.\, 
\sigma_{\overline{x},\sigma_{\overline{x'}}^{\overline{x}}(y)}^{\overline{x'},y}
(\, \mathsf{Ap}_\forall(  w, \, \sigma_{\overline{x'}}^{\overline{x}}(y))\, ) $ for
 $\overline{x}, \overline{x'}\in \Gamma^I$ and
 $w\in (\, \forall_{y\in B}\, C(y)\,)^I(\overline{x})$.
\\}

\end{definition}

\nocite{finmin}

\end{document}